\documentclass[11pt]{article}
\usepackage{amsfonts}
\usepackage{eucal}
\usepackage{amssymb,epsfig,latexsym,xcolor}
\usepackage{graphicx,amsmath,amssymb,latexsym,psfrag}
\usepackage{graphics,graphicx,comment,bm,empheq}
\usepackage{multirow,changepage}
\usepackage[subnum]{cases}

\newcommand{\cG}{\mathcal{G}}
\newcommand{\cA}{\mathcal{A}}

\newcommand{\cU}{\mathcal{U}}

\newcommand{\N}{\mathbb{N}}

\newcommand{\diff}{{\rm\,d}}

\newtheorem{Theorem}{Theorem}[section]
\newtheorem{Proposition}[Theorem]{Proposition}

\newtheorem{Corollary}[Theorem]{Corollary}
\newtheorem{Lemma}[Theorem]{Lemma}
\newtheorem{Remark}[Theorem]{Remark}

\newcommand{\nac}[2]{g_{#2}^{(#1)}}
\newcommand{\nan}{g_n}
\newcommand{\oD}{\mathcal{D}}
\newcommand{\ocalD}{\tilde{\mathcal{D}}}
\newcommand{\oE}{\mathcal{E}}
\newcommand{\oc}[1]{c_{#1}}
\newcommand{\ov}{\bm{\phi}}
\newcommand{\otaun}{v_n}
\newcommand{\otaumn}{v_n^{-1}}
\newcommand{\ozeta}{w}
\newcommand{\ote}{y_{\mathcal E}}
\newcommand{\oan}[1]{a_n^{(#1)}}
\newcommand{\on}[1]{n_n^{(#1)}}
\newcommand{\osn}[1]{n_{#1}}
\newcommand{\op}[1]{p_n^{(#1)}}
\newcommand{\oq}[1]{q_n^{(#1)}}
\newcommand{\obq}[1]{\bold{q}_n^{(#1)}}

\newcommand{\oTheta}{Q}
\newcommand{\oboldm}{\bold{u}}
\newcommand{\oldz}{\zeta}
\newcommand{\oldb}{\rho}
\newcommand{\oldB}{R}
\newcommand{\oldpi}{b}

\newcommand{\tgifeps}[3]{
\begin{figure}[htb]
\centering
\includegraphics[width=#1cm]{#2.eps}
 \vspace{-1mm}
\caption{#3\label{fig:#2}}
\vspace{-2mm}
\end{figure}
}

\newcommand{\sidebyside}[4]{
\begin{figure}
\begin{minipage}[t]{8cm}
\includegraphics[width=0.95\textwidth]{#1.eps}
\caption{#2\label{fig:#1}}
\end{minipage}
\hfill \hspace{0.4cm}
\begin{minipage}[t]{8cm}
\includegraphics[width=0.95\textwidth]{#3.eps}
\caption{#4\label{fig:#3}}
\end{minipage}
\hfill
\vspace{-0mm}
\end{figure}
}

\begin{document}
\title{Bootstrap percolation on the stochastic block model \\ with $k$ communities}
\author{Giovanni Luca Torrisi\thanks{
CNR, Roma, Italy. e-mail: \tt{giovanniluca.torrisi@cnr.it}} \and
Michele Garetto\thanks{
Universit\`{a} di Torino, Italy. e-mail: \tt{michele.garetto@unito.it}} \and
Emilio Leonardi\thanks{ 
Politecnico di Torino, Italy. e-mail: \tt{emilio.leonardi@polito.it}}}

\date{}
\maketitle

\begin{abstract}
We analyze the bootstrap percolation process on the stochastic
block model (SBM), a natural extension of the Erd\"{o}s--R\'{e}nyi random graph
that allows representing the \lq\lq community structure" observed in many real systems.
In the SBM,  nodes are partitioned into subsets, which represent different communities,
and pairs of nodes are independently connected with a probability that depends on the communities they belong to.
Under mild assumptions on  system parameters, we prove the existence of a sharp phase transition
for the final number of active nodes and characterize sub-critical and  super-critical regimes
in terms of the number of initially active nodes, which are selected uniformly at random in each community.
\end{abstract}


\section{Introduction}

Bootstrap percolation on a graph is a simple activation process
that starts with a given number of initially active nodes (called seeds) and evolves as follows.
Every inactive node that has at least $r\geq 2$ active neighbors is activated, and remains so forever.
The process stops when no more nodes can be  activated. There are two main cases of interest: one in which the seeds are selected uniformly at random
among the nodes, and one in which the seeds are arbitrarily chosen.
In both cases, the main question concerns the final size of the set of active nodes.

Bootstrap percolation was introduced in \cite{chalupa} on a Bethe lattice, and successively investigated
on regular grids and trees \cite{bollobas,BPP}.
More recently, bootstrap percolation has been studied
on random graphs and random trees \cite{fountolakis2,amini1,amini2,angel,pittel,galton,rgg,fountolakis,JLTV,kozma16,turova15}, motivated by the increasing interest in large-scale complex systems
such as technological, biological and social networks. For example, in the case of social networks,
bootstrap percolation may serve as a primitive model for the spread of ideas, rumors and
trends among individuals. Indeed, in this context one can assume that a person will adopt
an idea after receiving sufficient influence by friends who have already adopted it \cite{kempe,munik,watts}.

More in detail, bootstrap percolation has been studied on random regular graphs \cite{pittel},
on random graphs with given vertex degrees \cite{amini1}, on Galton-Watson random trees \cite{galton}, on random geometric graphs \cite{rgg}, on Chung-Lu random graphs \cite{amini2,fountolakis},
on small-world random graphs \cite{kozma16,turova15} and on Barabasi-Albert random graphs \cite{fountolakis2}.
Particularly relevant to our work is reference \cite{JLTV}, where the authors have provided a detailed analysis
of the bootstrap percolation process on the Erd\"{o}s-R\'{e}nyi random graph. We emphasize that in \cite{JLTV}  seeds are chosen uniformly at random among the nodes.
However, the critical number of seeds triggering percolation can be significantly reduced if  seed selection procedure  is
optimized~\cite{feige}.

Over the years, several variants of classical bootstrap percolation on a graph have
been considered. In majority bootstrap percolation,
a node becomes active if at least half of its neighbors are active.
In jigsaw percolation, introduced in \cite{jigsaw}, there are two types of edges among nodes, one
representing \lq\lq social links" and one representing
\lq\lq compatibility of ideas". Two clusters of nodes merge each other
if there exists at least one edge of each type between them.
Majority and jigsaw bootstrap percolation have been
analyzed on the Erd\"{o}s--R\'{e}nyi random graph in \cite{cecilia17} and
\cite{belajigsaw}, respectively.

An important characteristic of real graphs that is not captured by random graphs
models on which bootstrap percolation (and its variants) has been studied so far,
is the \lq\lq community structure".
Informally, one says that a graph has a \lq\lq community structure"
if nodes are partitioned into clusters in such a way that many edges join nodes of the same cluster and comparatively
few edges join nodes of different clusters \cite{girvanpnas02}. Many methods have been proposed for community detection in real networks
(see review article \cite{fortunato}). In the development of theoretical foundations of community detection, it has recently attracted considerable attention the so-called stochastic block model (SBM),
which is basically defined as the superposition of different Erd\"{o}s--R\'{e}nyi random graphs. In particular, detection of
two symmetric communities has been studied in \cite{mass14}, while partial or exact recovery of the community
membership has been investigated in \cite{abbe15} and \cite{abbe16}.

In this paper we study classical bootstrap percolation on the SBM, assuming
that  seeds are selected uniformly at random within each community and allowing
a different number of seeds for different communities.
We prove the existence of a sharp phase transition for the final size of the set of active nodes,
identifying a sub-critical regime, in which the bootstrap percolation process essentially does not evolve,
and a super-critical regime, in which
the activation process essentially percolates the whole graph.
Our results generalize and strengthen the main achievements in \cite{JLTV}.
We emphasize that our techniques significantly differ from those  employed in \cite{JLTV}.
In particular, we devise a suitable extension of the classical \lq\lq binomial chain" construction
originally proposed in \cite{Scalia} (which was successively applied in \cite{JLTV} to Erd\"{o}s--R\'{e}nyi
graphs), adapting it to the SBM.
Furthermore, differently from \cite{JLTV} where martingales concentration inequalities are exploited,
we resort on concentration inequalities for the binomial distribution to prove that bootstrap percolation on the SBM concentrates around its average. Our approach provides exponential bounds on the related tail probabilities,
which allow us, under a mild additional assumption, to strengthen  convergence in \lq\lq probability" for the final size
of active nodes (as obtained in \cite{JLTV}) to the level of \lq\lq almost sure" convergence.
To better understand the main difficulties in the analysis of bootstrap percolation on the SBM, we recall
that in the classical \lq\lq binomial chain" construction a (virtual) discrete time is introduced and at each time step
a single active node is \lq\lq explored" by revealing its neighbors. Nodes become active as soon as the number of their \lq\lq explored" neighbors reaches percolation threshold $r$.
In the SBM stochastic properties of the set of active nodes at time step $t$ heavily depend on
the number of nodes that have been \lq\lq explored" in each community up to time $t$, and this makes the analysis
of the bootstrap percolation process on the SBM
significantly more complex.
In particular, it
requires the identification of an appropriate \lq\lq strategy" to select the community in which a new node is \lq\lq explored" at every time step.

We acknowledge that the SBM, in spite of its flexibility to represent a wide variety of community-based
systems and its mathematical tractability, does not accurately describe most real-world graphs. Indeed, it does not model the possible
heterogeneity among nodes belonging to the same community.
Different variants of the SBM have been proposed to better fit real network data, such as making nodes
to follow a given degree sequence \cite{amin10,newman11} or considering overlapping communities and mixed membership models \cite{airoldi08,gopalan13}.
The investigation of bootstrap percolation  on the SBM is certainly a first step
towards the analysis of this process on more sophisticated community-based models.

The paper is organized as follows. In Section \ref{sec:SBM} we introduce the model, we describe
the extension to the SBM of the classical \lq\lq binomial chain" representation of bootstrap percolation, and we state
the model assumptions. The main results of the paper are stated in Section \ref{sec:main} where, to better convey the ideas which
lead us to identify critical conditions for the bootstrap percolation process on the SBM,
we briefly recall also the main achievements in \cite{JLTV}, and the intuition behind them.
In Section  \ref{sec:discussion} we discuss some consequences of our results with the help of numerical illustrations.
All proofs are presented in Section \ref{sec:proofs}.

\section{The stochastic block model}\label{sec:SBM}

\subsection{Model's description}

The SBM $G=G(\{\osn{i}\}_{1\leq i\leq k},\{q_{ij}\}_{1\leq i,j\leq k})$, $k,\osn{i}\in\mathbb N:=\{1,2,\ldots\}$, $\sum_{i=1}^{k}\osn{i}=n$,
$q_{ij}\in [0,1)$, $q_{ij}=q_{ji}$, is a random graph formed by the superposition of $k$ Erd\"os-R\'enyi's random graphs $G_i=G(\osn{i},p_i)$, $p_i:=q_{ii}$,
called hereafter communities, where edges joining nodes of communities $G_{i}$ and $G_j$, $i\neq j$, are independently
added with probability $q_{ij}$.

The bootstrap percolation process on the SBM is a nodes' activation process which obeys to the following rules:
\begin{itemize}
\item Nodes can be active or inactive.
\item At the beginning, an arbitrary number $a_i$ ($0 \leq a_i \leq \osn{i}$) of nodes, called seeds, are chosen
uniformly at random among the nodes of $G_i$. Seeds are declared to be active.
\item Nodes not belonging to the set of seeds are initially inactive.
\item An inactive node becomes active as soon as at least $r\geq 2$ of its neighbors are active.
\item Active nodes never become inactive and so the set of active nodes grows monotonically.
\item The process stops when no more nodes can be activated.
\end{itemize}
Bootstrap percolation naturally evolves through generations
of nodes that are sequentially activated. Zeroth generation $\cG_0$ is the set of seeds; first generation $\cG_1$ is composed by all the nodes which are activated by seeds;
second generation $\cG_2$ is composed by all the nodes that are activated
by both seeds and   nodes of the first generation, and so on. Bootstrap percolation  stops when either
an empty generation is obtained or all the nodes are active. The final set of active nodes is clearly given by
\[
\cG\equiv\bigcup_{h\in\N\cup\{0\}}\cG_h.
\]

We propose an extension of the classical \lq\lq binomial chain" representation of the system dynamics,
which makes easier
the analysis of the final size of active nodes $|\cG|$.
We introduce  a (virtual) discrete time $t\in\N\cup\{0\}$
and we assign a marks counter $M_{v}(t)$, $M_v(0):=0$, to every node $v$ which is not a seed. Seeds are activated at time $t=0$.
At every time step $t$ a single active node from one of the $k$ communities is \lq\lq used", i.e., \lq\lq explored"
by revealing its neighbors and by adding a mark to each of them. Nodes, which are not seeds, become active as soon as they collect $r$ marks.
For $i=1,\ldots,k$, we denote by $\mathcal{A}_{i}(t)$ and $\mathcal{U}_{i}(t)$
the set of active nodes at time $t$ in community $G_i$ and the set of \lq\lq used" nodes up to time time $t$ in community $G_i$,
respectively. We set $\mathcal{U}_{i}(0):=\emptyset$ and denote by $\mathcal{A}_{i}(0)$ the set of seeds in community $G_i$.

The process evolves according to the following recursive procedure. At time $t\in\mathbb N$:
\begin{itemize}
\item We select arbitrarily a community, in which active and not yet \lq\lq explored" nodes
are available, i.e., we select arbitrarily a community $G_j$ such that $\mathcal{A}_{j}(t-1)\setminus\mathcal{U}_{j}(t-1)\neq\emptyset$.
More formally, we select a community $G_j$ according to probability
\[
\frac{C_j(t)}{\sum_{i=1}^{k}C_i(t)},
\]
where $\{C_{i}(t)\}_{1\leq i\leq k}$ are arbitrary non-negative random variables such that:\\
\noindent $(i)$ $C_i(t)=0$,  if $\mathcal{A}_{i}(t-1)\setminus\mathcal{U}_{i}(t-1)=\emptyset$;\\
\noindent $(ii)$ $\sum_{i=1}^{k}C_i(t)\in\mathbb{R}_+:=(0,\infty)$ if there exists $i'\in\{1,\ldots,k\}$ such that $\mathcal{A}_{i'}(t-1)\setminus\mathcal{U}_{i'}(t-1)\neq\emptyset$.
\item We choose uniformly at random a node $v_t^{(j)}\in\mathcal{A}_{j}(t-1)\setminus\mathcal{U}_{j}(t-1)$.
\item 
We \lq\lq use" chosen node $v_t^{(j)}$, i.e., we \lq\lq explore" the node by revealing its neighbors and by adding a mark to each of them.
\item We set $\mathcal{A}_{i}(t):=\mathcal{A}_{i}(t-1)\cup\Delta\mathcal{A}_{i}(t)$, for $1 \le i\le k$, where
$\Delta\mathcal{A}_{i}(t)$ is the set of nodes in community $G_i$ that become active exactly at time $t$, i.e., the set of nodes in community $G_i$
that collected the $r$th mark exactly at time $t$. Note that, by construction, $\Delta\mathcal{A}_{i}(t)=\emptyset$ if $t<r$.
We also set $\mathcal{U}_{j}(t):=\mathcal{U}_{j}(t-1)\cup\{v_t^{(j)}\}$  and $\mathcal{U}_{i}(t):=\mathcal{U}_{i}(t-1)$, for $1 \le i\le k$ and  $i\neq j$.
\item The process terminates as soon as
$\mathcal{A}_{i}(t)=\mathcal{U}_{i}(t)$, $\forall$ $i=1,\ldots,k$, i.e., at time
\begin{equation}\label{eq:T}
T:=\min\{t\in\mathbb N:\,\,\mathcal{A}_{i}(t)=\mathcal{U}_{i}(t),\,\forall i=1,\ldots,k\}.
\end{equation}
\end{itemize}
Throughout this  paper  random variables $\{C_i(t)\}_{1\leq i\leq k,\,t<T}$   specify  the \lq\lq strategy"
(or \lq\lq policy") followed by the bootstrap percolation process. We emphasize that random variables $C_i(t)$ are introduced to formalize in mathematical terms
the selection mechanism of the community in which a new node is \lq\lq used" at a given time step.
The very weak assumptions that we impose on such random variables allow us to represent essentially
any possible way to select communities.

Note that, by construction, for any $t\leq T$,
\begin{equation*}
|\mathcal{U}(t)|=t,\quad\text{where}\quad\mathcal{U}(t):=\mathcal{U}_{1}(t)\cup\ldots\cup\mathcal{U}_{k}(t).
\end{equation*}
Indeed, at every time $t\leq T$, the process \lq\lq explores" a node belonging only to one of the $k$ communities.
Let $\mathcal{A}(t):=\bigcup_{i=1}^{k}\mathcal{A}_i(t)$ denote the set of active nodes at time $t\leq T$ and consider a node $v\not\in\mathcal{A}(0) $. We clearly have
\begin{equation}\label{eq:AMv}
v\in\mathcal{A}(t)\quad\text{if and only if}\quad M_v(t)\geq r,\quad 1\leq t\leq T.
\end{equation}
We also have
\begin{equation}\label{eq:Mv}
M_{v}(t)=\sum_{i=1}^{k}\sum_{s=1}^{U_{i}(t)}I^{(i)}_{v}(s),\quad\text{$\forall$ $v\not\in\mathcal{A}(t)$}
\end{equation}
where $U_{i}(t):=|\mathcal{U}_{i}(t)|$ and random variables $\{I_{v}^{(i)}(s)\}_{v\notin\cA(t),1\leq i\leq k,1\leq s\leq T}$
are independent, with $I_{v}^{(i)}(s)$  distributed as $\mathrm{Be}(q_{ij})$ when $v$ is a node of community $G_j$.
Here $\mathrm{Be}(p)$ denotes a Bernoulli distributed random variable with mean $p\in [0,1]$.
The following proposition, whose proof is given in Subsection \ref{subsec:prop21}, holds.
\begin{Proposition}\label{prop:equiv}
For any \lq\lq strategy" $\{C_i(t)\}_{1\leq i\leq k,\,t<T}$,  we have
\[
\cG\equiv\mathcal{A}(T).
\]
\end{Proposition}
We have defined random marks $I_v^{(i)}(s)$
for $v\notin\cA(t)$, $1\leq i\leq k$ and $1\leq s\leq T$, but, similarly to \cite{JLTV},
it is possible to introduce additional, redundant random marks, which are independent and Bernoulli distributed with mean $q_{ij}$ if $v$ is a node of community $G_j$,
in such a way that $I_v^{(i)}(s)$ is defined for all $v\in G$, $1\leq i\leq k$ and  $s\in\mathbb N$.
Such additional random marks are added, for any $1\leq s\le T$,
to already active nodes and so they have no effect on the
underlying bootstrap percolation process.
Throughout this paper,
we denote by $\mathrm{Bin}(u,p)$, $u\in\mathbb N$, $p\in [0,1]$, a random variable following the binomial distribution with parameters $(u,p)$.
For fixed $i\in\{1,\ldots,k\}$ and $t\in\mathbb N$, we have
\begin{equation}\label{eq:sumbin}
\text{$M_{v}(t)\,|\,\bold{U}(t)=\oboldm(t)\overset{\mathcal L}=\sum_{j=1}^{k}\mathrm{Bin}(u_j(t),q_{ij})$,\quad $v\in G_i$}
\end{equation}
where $\bold{U}(t):=(U_1(t),\ldots,U_k(t))$, $\oboldm(t):=(u_1(t),\ldots,u_k(t))$, symbol $\overset{\mathcal L}=$ denotes equality in law and random variables
$\mathrm{Bin}(u_j(t),q_{ij})$, $j=1,\ldots,k$, are independent. The number of active nodes in community $G_i$ at time $t\in\mathbb N\cup\{0\}$ is given by
\begin{equation}\label{eq:Ai}
A_{i}(t):=|\mathcal{A}_{i}(t)|=a_{i}+S_{i}(t),
\end{equation}
where
\begin{equation}\label{eq:Si}
S_{i}(t):=\sum_{v\in G_i\setminus\mathcal{A}_{i}(0)}\bold{1}\{Y_{v}\leq t\},\qquad Y_{v}:=\min\{s\in\mathbb N:\,\,M_{v}(s)\geq r\}.
\end{equation}
Since random variables
\[
\{M_{v}(t)\,|\,\bold{U}(t)=\oboldm(t)\}_{v\in G_i}
\]
are independent and identically distributed with law specified by \eqref{eq:sumbin}, we have
\begin{equation}\label{eq:bin}
S_{i}(t)\,|\,\bold{U}(t)=\oboldm(t)\overset{\mathcal L}=\mathrm{Bin}(\osn{i}-a_{i},\oldpi(\oboldm(t),\bold{q}_i)),
\end{equation}
where
\begin{equation}\label{eq:pi}
\oldpi(\oboldm(t),\bold{q}_i):=P\left(\sum_{j=1}^{k}\mathrm{Bin}(u_j(t),q_{ij})\geq r\right),
\end{equation}
and $\bold{q}_i:=(q_{ij})_{1\leq j\leq k}$. Hereafter, we denote by $A(t):=|\mathcal{A}(t)|=\sum_{i=1}^{k}A_i(t)$, the number of active nodes in the SBM at time $t$,
and by $A_*:=A(T)=T$ the final number of active nodes.

\begin{Remark}\label{rem1}
The analysis of the bootstrap percolation process is significantly more complex on the SBM
than on the Erd\"os-R\'enyi random graph due to the following two reasons. $(i)$ On
the SBM, at each time step $t\le T$, we select, according to the chosen \lq\lq strategy", a community
in which  \lq\lq exploring" an active node.
In particular, note that, for any $1\leq i\leq k$ and $t<T$,
random variables $A_i(t)$ heavily depend on  quantities $\{U_i(t)\}_{1\leq i\leq k}$, which are themselves constrained by the availability of \lq\lq usable" nodes in  communities,
and therefore
on the adopted  \lq\lq strategy". In contrast, to analyze the bootstrap percolation process on
the Erd\"os-R\'enyi random graph
there is clearly no need to introduce any \lq\lq policy". $(ii)$
As a consequence of $(i)$,  for any $i=1,\ldots,k$ and $t<T$, the law of random variables
$S_i(t)$ is binomial only given event $\{\bold{U}(t)=\oboldm(t)\}$. Therefore, for any $1\leq i\leq k$ and $t<T$,
the probabilistic structure of random variables $S_i(t)$ is
significantly more complex on the SBM than on the Erd\"os-R\'enyi random graph, i.e., for $k=1$. On the latter graph, indeed,  $U_1(t)=t$
and  the law of $S_1(t)$ is binomial with parameters $(\osn{1}-a_1,\oldpi(t,p_1))$, where $\oldpi(t,p_1):=P(\mathrm{Bin}(t,p_1)\geq r)$.

Finally, we remark that thanks to Proposition \ref{prop:equiv}, differently from quantities $A_i$ and $U_i$, the final number of active nodes
$A^*$ does not depend on the chosen \lq\lq strategy".
\end{Remark}

\begin{Remark}
A possible alternative way to extend the classical \lq\lq binomial chain" construction to the SBM could be the following.
At each time step $t\le T$, a node is \lq \lq used"
in each community $G_i$ in which at least one \lq\lq usable" node can be found, i.e.,
in each community $G_i$ in which ${\mathcal A}_i(t-1)\setminus{\mathcal  U}_i(t-1)\neq \emptyset$. Although this alternative extension of the classical
\lq\lq binomial chain" construction on the SBM makes the system dynamics independent on the \lq\lq strategy",
it appears pretty difficult to analyze due to the complex
probabilistic structure of random variables $\{A_i(t)\}_{1\leq i\leq k,\,t<T}$ and $\{U_i (t)\}_{1\leq i\leq k, t<T}$.
It will appear from our investigation that the opportunity to \lq\lq arbitrarily" define the \lq\lq policy", according to which communities are
selected, provides a degree of flexibility that comes in handy to analyze the process
(see the proofs of Theorems \ref{teo-blocksub} and \ref{teosupcrit}).
\end{Remark}

\begin{table}
\begin{center}
\caption{Main notation}\label{main-notation}
\begin{adjustwidth} {-1.1 cm}{}
 \begin{tabular}{|l|l|l|}
 \hline
 \multicolumn{3}{|c|}{Graph parameters}\\
 \hline
 Symbol & Mathematical definition & Description\\
 \hline
  $n$ &  & number of nodes in the SBM \\
  $k$ & &   numbers of communities \\
  $G_i$ & & $i$th community \\
  $\on{i}$ & &  number of nodes in $G_i$ \\
  $\oq{{ij}}$ & &  edge prob. between nodes in $G_i$ and $G_j$ \\
  $\op{i}$ & \small{$\oq{{ii}}$}  & edge prob. between nodes in $G_i$\\
  $\nu_{ij}$ & $\lim_{n\to \infty} \on{i}/\on{j}\in\mathbb{R}_+$ \,\, (see \eqref{eq:a1a3})) & \\
   $\gamma_{ij}$ & \small{$\lim_{n\to \infty} \oq{{ij}}/\op{i} $}  \,\, (see \eqref{eq:gratiog}))& \\
    $\mu_{ij}$ & \small{ $\lim_{n\to \infty} \op{i}/\op{j}$ } \,\,  (see \eqref{eq:gratiog})& \\
    $\chi_{ij}$ & \small{$\gamma_{ij}(\nu_{ij}(\mu_{ij})^r)^{1/(r-1)}$} (see \eqref{eq:chi-add}) & \\
    $\widetilde G$ &   & community-level graph\\
    $\obq{i}$ & $(\oq{{ij}})_{1\le j\le k}$ & \\
    \hline
     \multicolumn{3}{|c|}{Bootstrap percolation parameters}\\
     \hline
    Symbol & Mathematical definition & Description\\
     \hline
     $r$ & & bootstrap percolation threshold\\
  $\nac{i}{n}$ & $(1-r^{-1})\Big(\frac{(r-1)!}{\on{i} (\op{{i}})^r}\Big)^{1/(r-1)}$ \,\, (see \eqref{eq:acrit}) & critical number of seeds in $G_i$ \\
  $\oan{i}$ & & number of seeds in $G_i$\\
   $\alpha_i$ & $\lim_{n\to \infty} \oan{i}/\nac{i}{n} $ \;\; (see \eqref{eq:trivial})& \\
 $ \lfloor \bold{x}\nan\rfloor $ & $ (\lfloor x_1 \nac{1}{n}\rfloor , \ldots  ,\lfloor x_k \nac{k}{n}\rfloor )$  &   \\
    \hline
     \multicolumn{3}{|c|}{Bootstrap percolation main dynamical variables}\\
     \hline
 Symbol & Mathematical definition & Description\\
 \hline
   $\mathcal A_i(t)$ &  & set of active nodes in  $G_i$ at time $t$\\
   $\mathcal U_i(t) $  &  & set of \lq\lq used" nodes in $G_i$ at time $t$\\
   ${\color{black}A_n^{(i)}}(t)$ & see \eqref{eq:Ai} & number of active nodes in  $G_i$ at time $t$\\
  ${\color{black}U_n^{(i)}}(t)$ &  & number of \lq\lq used" nodes in $G_i$ at time $t$\\
  $A_n(t)$ & $\sum_{i=1}^k {\color{black}A_n^{(i)}}(t)$  & \\
  $T_n$ & $ \min\{t: \bold{A}_n(t)=\bold{U}_n(t)\}$ \,\, (see \eqref{eq:T}) & termination time of the bootstrap process \\
  $A_n^*$ & $A_n(T_n)$ & final size of active nodes\\
  $M_v(t)$ & see \eqref{eq:Mv}  &  marks counter of node $v$ at time $t$ \\
  $\{{\color{black}C_n^{(i)}}(t)\}_{1\le i\le k}$ & & variables defining the \lq\lq strategy" \\
  ${\color{black}S_n^{(i)}}(t)$ & ${\color{black}A_n^{(i)}}(t)-\oan{i}$ \,\, (see \eqref{eq:Ai} and \eqref{eq:Si}) & \\
  $Y_v$ & $\min\{t\in \mathbb{N}: M_v(t)\ge r\}$ \,\, (see \eqref{eq:Si}) & activation time of node $v$\\
\hline
 \multicolumn{3}{|c|}{Variables representing average dynamics}\\
 \hline
    Symbol & Mathematical definition & Description\\
     \hline
$\oldpi(\oboldm_n(t),\obq{i})$ & \small{$P(\sum_{i=1}^k \text{Bin}(u_n^{(j)}(t), \oq{{ij}})\ge r)$ \;\, (see \eqref{eq:sumbin} and \eqref{eq:pi})}  &  prob. that a node in $G_i$ is active at time $t$  \\
${\color{black}\oldB_n^{(i)}}(\oboldm_n(t), \obq{i} ) $  & \small{$E[{\color{black}A_n^{(i)}}(t)-{\color{black}U_n^{(i)}}(t)\mid\bold{U}_n(t)=\oboldm_n(t)]$} \,\, (see \eqref{eq:R})
& mean of \lq\lq usable" nodes in $G_i$ at time $t$\\
$\oldb_i(\bold{x})$ & \small{$\lim_{n\to\infty}\frac{{\color{black}\oldB_n^{(i)}}(\lfloor \bold{x}\nan\rfloor,\obq{i} )}{\nac{i}{n}}$} \,\, (see \eqref{eq:Btob})  & \\
${{\bm\oldb}}$ & $(\oldb_1,\ldots,\oldb_k)$ & \\
 \hline
 \end{tabular}
 \end{adjustwidth}
\end{center}
\end{table}

\begin{table}
 \begin{center}
\caption{Further notation}\label{main-notation2}
 \begin{tabular}{|l|l|l|}
  \hline
  \multicolumn{3}{|c|}{Further parameters and  variables}\\
     \hline
   Symbol & Mathematical definition & Description/Remarks\\
 \hline
  \multirow{2}{*}{$\mathcal{K}_h$} &  \multirow{2}{*}{see \eqref{eq:mathcalk}} & set of communities \\
     & &  at dist. $h$ from $G_1$ on $\widetilde{G}$\\
 ${\color{black}d_i}$ & &  dist. $G_i \leftrightarrow G_1$ on $\widetilde{G}$ \\
 ${\color{black}\overline{d}}$ & see \eqref{eq:maxdist} & $\max_{i} {\color{black}d_i}$\\
 $\chi$ & $\min\{\chi_{ij}: \chi_{ij}>0\}$ \,\, (see \eqref{eq:chimin}) &  \\
  $\beta_h$ & $\frac{r^{-1}(1-r^{-1})^{r-1}}{2}(\chi \beta_{h-1})^{r-1}$, $\beta_0=\alpha_1/2$ \, (see \eqref{eq:beta})& \\
  $J_{\bm\oldb}(\bold x)$ &Jacobian of ${{\bm\oldb}}$ at $\bold x$ \; (see Sect. \ref{sect:jacobian}) & \\
  $\lambda_{PF}(\bold{x})$ &   eigenvalue of $J_{\bm\oldb}(\bold x)$ associated to $\bold \ov_{PF}(\bold{x})$  (see Sect. \ref{sect:jacobian} ) & relevant eigenvalue  \\
  $\bold \ov_{PF}(\bold{x})$ & eigenvector of $J_{\bm\oldb}(\bold x)$ associated to $\lambda_{PF}(\bold{x})$ \, (see Sect. \ref{sect:jacobian}) & $\bold \ov_{PF}(\bold{x})\ge\bold 0$\\
  \hline
     \multicolumn{3}{|c|}{Regions/Sets}\\
     \hline
   Symbol & Mathematical definition & Description/Remarks\\
 \hline
$\mathcal{R}_{Sub}$ &  $\{ \bm{\alpha}: ({\mathcal Sub}) \text{ holds} \}$ \, (see \eqref{Rcond-sub-crit}) & sub-critical region \\
  $\mathcal{R}_{Crit}$ &  $\{ \bm{\alpha}: ({\mathcal Crit}) \text{ holds}  \}$ \, (see \eqref{Rcond-sub-crit}) & critical region \\
 $\mathcal{R}_{Sup}$ & $\{ \bm{\alpha}:  ({\mathcal Sup}) \text{ holds} \}$ \, (see \eqref{Rcond-sup}) & super-critical region  \\
  $\oD$ & $\{\bold{x}\in [0,r/(r-1)]^k : x_i+ \sum_{j\neq i}^{1,k} x_j\chi_{ij}\le r/(r-1) \}$ \, (see \eqref{eq:D}) &  \\
  $\oD_{\bm\oldb}$ & $ \{\bold{x}\in\oD :\;\ \oldb_1(\bold{x})=\oldb_2(\bold{x})=\ldots=\oldb_k(\bold{x})\}$\: (see below \eqref{eq:D}) &  \\
  $\ocalD$  & $\cup_{i=1}^k \{ \bold{x}\in \oD : x_i + \sum_{j\neq i}^{1,k} x_j\chi_{ij}= r/(r-1) \}$ \,(see Sect. \ref{sect:jacobian})  & \\
  $\oE'$ & $\{ \bold{x}\in\mathcal D : {\bm\oldb}(\bold{x})\ge\bold 0 \}$ \, (see Sect. \ref{cauchy-prob}) &  \\
  $\oE$ & closure largest connect. comp. of $\oE'$: $\bold{x}_0\in\oE$ \, (see Sect. \ref{cauchy-prob}) & $\bold{0}\in \oE$\\
 $\mathcal S_h$ & $\{ \theta \bold x_0^{(h)} + (1-\theta ) \bold x_0^{(h-1)}:\, \theta\in [0,1] \}$ \, (see Sect. \ref{sect:equiv}) & \\
  $\mathcal S $ & $\cup_{1\le h\le  {\color{black}\overline{d}}} \mathcal S_h$ \,\, (see Sect. \ref{sect:equiv})& \\
  $\tilde{\mathcal S} $ & $ \{ \bold x \in \oD : \bold x=\bold x(\ote)+ \theta \bold \ov_{PF} (\bold x(\ote)), \theta\ge 0 \}$ \,\, (see Sect. \ref{sect:equiv}) & \\
$\mathcal C_{\oE}$ & $\{\bold x(y) :  y \in [y_0, \ote]$ \} \,(see Sect. \ref{sect:equiv})&  \\
$\mathcal C$ & $\mathcal S \cup \mathcal C_{\oE} \cup \tilde{\mathcal S}$ \, (see Sect. \ref{sect:equiv})&  $\mathcal C= \mathcal S \cup \mathcal{C}_{\oE} $ \mbox{if} $\ote<\infty$ \\
$\mathcal Z$ & set of distinct zeros of ${\bm\oldb}$ within $\oD$\, (see Sect. \ref{sect:equiv}) & \\
  \hline
     \multicolumn{3}{|c|}{Special points}\\
     \hline
 Symbol & Mathematical definition & Description/Remarks\\
 \hline
  $\bold{x}_0^{(h)}$ & $\sum_{s=0}^{h}\beta_{s}\bold{1}_{\mathcal K_s}$ \, (see \eqref{eq:x0h})&  \\
 $\bold{x}_0$ & $\bold x_0^{({\color{black}\overline{d}})}$ \, (see \eqref{eq:x0}) & i.c. of the Cauchy problem  \\
 $\ote$ & $\inf \{ y\in (y_0,y_+): \bold{x}(y)\not \in \overset\circ{\oE}\}$ \,(see \eqref{eq:tE}) &
 $\ote=\infty$ if $\{\ldots \}=\emptyset$ \\
    $\bold{x}(\ote)$ & $\lim_{y\uparrow \ote} \bold{x}(y)$ \, (see \eqref{eq:x(tE)}) & $\bold{x}(y_{\mathcal E})\in\partial\mathcal E$ \\
\hline
 \end{tabular}
 \end{center}
\end{table}

\subsection{Model assumptions}\label{model-assumptions}

Hereafter, given two functions $f_1$ and $f_2$ we write $f_1(n)\ll f_2(n)$ (or equivalently $f_1(n)=o(f_2(n))$), $f_1(n)\sim f_2(n)$, $f_1(n)\sim_e f_2(n)$, $f_1(n)=O(f_2(n))$
and $f_1(n)\lesssim f_2(n)$ if, as $n\to\infty$, $f_1(n)/f_2(n)\to 0$, $f_1(n)/f_2(n)\to c\in\mathbb{R}\setminus\{0\}$,
$f_1(n)/f_2(n)\to 1$, $\limsup_{n\to\infty}|f_1(n)/f_2(n)|<\infty$ and either $f_1(n)=o(f_2(n))$ or $f_1(n)\sim_e f_2(n)$, respectively.

In the following we consider a sequence of SBMs with a growing number of nodes $n$ and we explicit the dependence on $n$ of the related variables.
We warn the reader that, unless explicitly written, all the limits in this paper are taken as $n\to\infty$.

For any $i,j\in\{1,\ldots,k\}$, we assume
\begin{equation}\label{eq:a1a3}
\lim_{n\to\infty}\on{i}/\on{j}=\nu_{ij}\in\mathbb{R}_+,
\end{equation}
where $n_n^{(i)}:=n_i$ (i.e., $n_n^{(i)}$  replaces $n_i$).  Since $n=\sum_{i=1}^{k}\on{i}$ it follows
\[
\on{i}\sim_e \left(1+\sum_{j\neq i}^{1,k}\nu_{ji}\right)^{-1}n.
\]
For any $i,j\in\{1,\ldots,k\}$, we assume
\begin{equation}\label{eq:hyp2bis}
1/\on{i}\ll \op{{i}}\ll 1/(\on{i})^\frac{1}{r},\quad\text{where $\op{i}:=p_i$}
\end{equation}
\begin{equation}\label{eq:hyp4}
\oq{{ij}}:=q_{ij}\leq\min\{\op{{i}},\op{{j}}\},
\end{equation}
\begin{equation}\label{eq:gratiog}
\text{
$\gamma_{ij}:=\lim_{n\to\infty}\oq{{ij}}/\op{{i}}
\in [0,1]$, $\mu_{ij}:=\lim_{n\to\infty}\op{{i}}/\op{{j}}
\in\mathbb{R}_+$.}
\end{equation}
For any $i\in\{1,\ldots,k\}$, we define
\begin{equation}\label{eq:acrit}
\nac{i}{n}:=\left(1-\frac{1}{r}\right)\left(\frac{(r-1)!}{\on{i} (\op{{i}})^r}\right)^{\frac{1}{r-1}},
\end{equation}
and assume
\begin{equation}\label{eq:trivial}
\oan{{i}}/\nac{i}{n} \to \alpha_i\geq 0\quad\text{with $\alpha_i>0$ for some $i\in\{1,\ldots,k\}$}
\end{equation}
where $\oan{i}:=a_i$. By \eqref{eq:a1a3}, \eqref{eq:gratiog} and \eqref{eq:acrit}, for any $i,j\in\{1,\ldots,k\}$, we have
\begin{equation}\label{eq:chi-add}
\chi_{ij}:=\lim_{n\to\infty}\frac{\oq{{ij}}\nac{j}{n}}{\op{{i}}\nac{i}{n}}=\gamma_{ij}(\nu_{ij}(\mu_{ij})^r)^{1/(r-1)}.
\end{equation}
Note that $0\leq\chi_{ij}<\infty$, for all $i,j\in\{1,\ldots,k\}$. In particular, $\gamma_{ij}=0$ implies $\chi_{ij}=\chi_{ji}=0$, while
$\gamma_{ij}>0$ implies $\chi_{ij},\chi_{ji}\in\mathbb{R}_+$. Throughout this paper, we assume that
\begin{equation}\label{eq:ala4}
\text{Matrix $\bm{\chi}:=(\chi_{ij})_{1\leq i,j\leq k}$ is irreducible.}
\end{equation}
Hereon, in addition to \eqref{eq:a1a3}, \eqref{eq:hyp2bis}, \eqref{eq:hyp4}, \eqref{eq:gratiog}, \eqref{eq:trivial} and \eqref{eq:ala4},
without loss of generality, we assume that communities are numbered so as to guarantee
\begin{equation}\label{eq:alfaorder}
\alpha_1\ge\alpha_2\ge\ldots\ge\alpha_k\quad\text{with $\alpha_1>0$.}
\end{equation}
Finally, we note that, as proved in \cite{JLTV}, under
\eqref{eq:a1a3} and \eqref{eq:hyp2bis}, we have, for any $i\in\{1,\ldots,k\}$,
\begin{equation}\label{eq:ptcto0bm}
\nac{i}{n}\to \infty, \quad \nac{i}{n}/\on{i}\to 0,\quad \op{{i}}\nac{i}{n}\to 0.
\end{equation}

Let us briefly discuss model assumptions. Condition \eqref{eq:a1a3} states that the $k$ different communities have sizes which are asymptotically of the same order.
Condition \eqref{eq:hyp2bis} guarantees that the average degree of nodes in each community $G_i$ tends to infinity,
as $n\to\infty$. Under such condition a sharp phase transition occurs with a negligible number of seeds, i.e., a number of seeds that is $o(n)$, as shown in \cite{JLTV} for the case $k=1$.

Condition \eqref{eq:hyp4} means that the SBM is \lq\lq weakly assortative", indeed,
for any community $G_i$, the \lq\lq intra-community" edge probability $\op{i}$ is
required to be not smaller than each \lq\lq extra-community" edge probability $\oq{{ij}}$ ($i\neq j$).
Condition \eqref{eq:gratiog} guarantees that, asymptotically, the \lq\lq intra-community" edge probabilities are
comparable with each other (i.e. are of the same order). Note that, although two different communities $G_i$ and $G_j$ ($i\neq j$) may be
completely disconnected (in the sense that no edges are established between them,
which occurs when $\oq{{ij}}=0$), thanks to condition \eqref{eq:ala4} there are no isolated communities. Indeed, condition \eqref{eq:ala4} guarantees that
graph $\widetilde{G}=\widetilde G(k,\bm{\chi})$, whose $k$ nodes represent communities and an edge is established between  node $i$ and  node $j$ if $\chi_{ij}>0$,
is  connected. We remark once again that the number of seeds in each community can be arbitrarily chosen (i.e.,
quantities $a_n^{(i)}$, $i=1,\ldots,k$, are arbitrary)
while seeds within each community must be selected uniformly at random. We also remark that, although throughout this paper it is assumed
$\bm{\alpha}:=(\alpha_1,\ldots,\alpha_k)\neq\bold 0$, it is easily checked that for $\bm{\alpha\equiv\bold{0}}$ the system is in  sub-critical conditions.

From now on, we explicit the dependence of $C_i$, $T$, $U_{i}$, $u_i$, $\bold{U}$, $\bold u$, $A_{i}$, $S_{i}$, $\bold{q}_i$, $A$ and $A_*$ on $n$
writing $C_n^{(i)}$, $T_n$, ${\color{black}U_n^{({i})}}$, $u_n^{(i)}$, $\bold{U}_n$, $\oboldm_n$, ${\color{black}A_n^{({i})}}$, ${\color{black}S_n^{({i})}}$, $\obq{i}$, $A_n$ and $A_n^*$, respectively.

For reader's convenience, we summarize the main notation of the paper in Tables \ref{main-notation} and \ref{main-notation2}.
With a few exceptions, in our notation we have adopted the following general rules.
$(i)$ System's parameters are denoted by small latin letters;
$(ii)$ dynamical variables are denoted by capital latin letters;
$(iii)$ sets (including curves) are denoted by calligraphic letters;
$(iv)$ asymptotic limits (as $n \rightarrow \infty$) are denoted by small greek letters.
Moreover, we have used the following rule when adding indexes to symbols: index $n$ is always
put as pedex; other indexes are preferably put as pedex, unless they conflict with $n$,
in which case they are put as apex (in round brackets in order to avoid confusion with exponentiation).

\section{Main results}\label{sec:main}

\subsection{Bootstrap percolation on the Erd\"os-R\'enyi random graph: a quick review}

To better understand the ideas which led us to identify sub-critical and super-critical conditions for
the bootstrap percolation process on the SBM,
we briefly recall the main results in \cite{JLTV}. Note that the Erd\"os-R\'enyi random graph corresponds to a special case of the the SBM, (i.e., when
$k=1$). It has been proved in \cite{JLTV} that:\\
\noindent $(i)$ If \eqref{eq:hyp2bis} and \eqref{eq:trivial} hold (with $k=1$) and $\alpha_1<1$, then
\[
A_n^*/\nac{1}{n}\to\frac{r\varphi(\alpha_1)}{(r-1)\alpha_1},\quad\text{in probability}
\]
where $\varphi(\alpha_1)$ is the unique solution in $[0,1]$ of equation $rx-x^r=(r-1)\alpha_1$ (see Theorem 3.1$(i)$ in \cite{JLTV}).\\
\noindent $(ii)$ If \eqref{eq:hyp2bis} and \eqref{eq:trivial} hold (with $k=1$) and $\alpha_1>1$, then
\[
A_n^*/n\to 1,\quad\text{in probability}
\]
(see Theorem 3.1$(ii)$ in \cite{JLTV}).

In words, below the critical number of seeds $\nac{1}{n}$, the bootstrap percolation process essentially does not evolve, reaching, as $n\to\infty$, a
final size of active nodes which is of the same order as $\oan{1}$
(sub-critical case); above critical number of seeds, $\nac{1}{n}$, the process percolates
through the entire graph, reaching, as $n\to\infty$,
a final size of active nodes which is of the same order as $n$
(super-critical case).

We briefly and informally describe the intuition behind these results. First of all note that, since $k=1$, we have
${\color{black}U_n^{(1)}}(t)=t$ for any $t\in\mathbb N$, and so
the number of \lq\lq active and not yet used" nodes at time $t$ is given by $\oan{1}+{\color{black}S_n^{(1)}}(t)-t$. Exploiting that
${\color{black}S_n^{(1)}}(t)\overset{\mathcal L}=\mathrm{Bin}(n-\oan{{1}},\oldpi(t,\op{1}))$, where
$\oldpi(t,\op{1}):=P(\mathrm{Bin}(t,\op{{1}})\geq r)$, one proves that, under assumptions
\eqref{eq:hyp2bis} and \eqref{eq:trivial} (with $k=1$), the process $\{\oan{1}+{\color{black}S_n^{(1)}}(t)-t\}_{t\in\mathbb N}$ is \lq\lq concentrated" around its mean, as $n\to\infty$, and
\begin{equation}\label{eq:b}
\lim_{n\to\infty}\frac{{\color{black}\oldB_n^{(1)}}(\lfloor x\nac{1}{n}\rfloor,\op{1})}{\nac{1}{n}}=\oldb_1(x):=\alpha_1-x+r^{-1}(1-r^{-1})^{r-1}x^r,\quad\text{$\forall$ $x\in [0,r/(r-1)]$}
\end{equation}
(see \cite{JLTV} for details), where
\begin{align*}
{\color{black}\oldB_n^{(1)}}(t,\op{1}):=E[\oan{1}+{\color{black}S_n^{(1)}}(t)-t]=\oan{1}+(n-\oan{{1}})\oldpi(t,\op{1})-t,
\end{align*}
and $\lfloor x\rfloor$ denotes the greatest integer less than or equal to $x\in\mathbb R$. A simple computation shows that function $\oldb_1$ has a unique point of minimum at $x=r/(r-1)$ and $\oldb_1(r/(r-1))=\alpha_1-1$.
Therefore, since bootstrap percolation  stops the first time at which the number of \lq\lq active and not yet used"
nodes equals zero we have that $(i)$ if $\alpha_1<1$, then bootstrap percolation stops at a time which is, asymptotically in $n$, of the same order as $\nac{1}{n}$,
$(ii)$ if $\alpha_1>1$, then bootstrap percolation does not stop at
times which are, asymptotically in $n$, of the same order as $\nac{1}{n}$. In the super-critical case,
 further analysis of  function $t\mapsto {\color{black}\oldB_n^{(1)}}(t,\op{1})$ at times bigger than $\nac{1}{n}$ reveals that
${\color{black}\oldB_n^{(1)}}(\cdot,\op{1})$ quickly increases up to time $t=t_n\sim_e cn$, where $c>0$ is an arbitrarily small positive constant, and
${\color{black}\oldB_n^{(1)}}(t_n,\op{1})\approx n-t_n$. Then ${\color{black}\oldB_n^{(1)}}(\cdot,\op{1})$ decreases linearly and hits zero at time $t_n\sim_e n$. As a consequence, for $\alpha_1>1$,
bootstrap percolation terminates at a time which is of the same order as $n$.

\subsection{Phase transition of the bootstrap percolation process on the SBM}

As recalled in the previous section, bootstrap percolation  on the Erd\"os-R\'enyi random graph exhibits a sharp phase transition;
the reader may be wondering whether more complex phenomena, such as selective percolation of communities,
can be observed on the SBM. We shall show that this is not the case. Indeed, under the assumptions described in Subsection \ref{model-assumptions}, in analogy with the case $k=1$,
the bootstrap percolation process either stops at time scales $\nac{1}{n}$ (sub-critical case),
or it percolates the whole graph (super-critical case).

To present our results, we start introducing the  asymptotic normalized
mean number of \lq\lq active and not yet used" nodes (i.e., function corresponding to $\oldb_1$). For $n,t\in\mathbb N$ and $i\in\{1,\ldots,k\}$, we set
\begin{align}
{\color{black}\oldB_n^{(i)}}(\oboldm_{n}(t),\obq{i}):&= E[{\color{black}A_n^{({i})}}(t)-{\color{black}U_n^{({i})}}(t)\,|\,\bold{U}_{n}(t)=\oboldm_{n}(t)]\label{eq:R}\\
&=\oan{{i}}+(\on{i}-\oan{{i}})\oldpi(\oboldm_{n}(t),\obq{i})-u_n^{(i)}(t).\nonumber
\end{align}
Hereon, for $\bold{x}:=(x_1,\ldots,x_k)\in [0,\infty)^k$, we set
\[
\lfloor \bold{x}\nan\rfloor:=(\lfloor x_1\nac{1}{n}\rfloor,\ldots,\lfloor x_k \nac{k}{n}\rfloor).
\]
The following lemma holds.

\begin{Lemma}\label{le:bas2}
Assume \eqref{eq:a1a3}, \eqref{eq:hyp2bis}, \eqref{eq:hyp4}, \eqref{eq:gratiog}, \eqref{eq:trivial} and
let $i\in\{1,\ldots,k\}$ be fixed.
Then
\begin{equation}\label{eq:Btob}
\lim_{n\to\infty}\frac{{\color{black}\oldB_n^{(i)}}(\lfloor\bold{x}\nan\rfloor,\obq{{i}})}{\nac{i}{n}}=\oldb_i(\bold x),
\quad\text{$\forall$ $\bold{x}\in [0,\infty)^k$}
\end{equation}
where
\[
\oldb_i(\bold x):=\alpha_i-x_i+r^{-1}(1-r^{-1})^{r-1}\left(\sum_{j=1}^{k}x_j\chi_{ij}\right)^r.
\]
\end{Lemma}

Setting
\begin{equation}\label{eq:D}
\oD:=\left\{\bold{x}\in [0,r/(r-1)]^k:\,\,x_i+\sum_{j\neq i}^{1,k}\chi_{ij}x_j\leq\frac{r}{r-1},\,\,\forall\,\,i=1,\ldots,k\right\}
\end{equation}
and
\[
\oD_{{\bm\oldb}}:=\{\bold x\in\oD:\,\,\oldb_1(\bold x)=\ldots=\oldb_k(\bold x)\},
\]
we shall distinguish among the following  conditions:
\\
\\
\noindent$({\mathcal Sub})$: $\min_{\bold x\in\oD_{\bm\oldb}}\oldb_1(\bold{x})<0$.\\
\noindent$({\mathcal Crit})$: $\min_{\bold x\in\oD_{\bm\oldb}}\oldb_1(\bold{x})=0$.\\
\noindent$({\mathcal Sup})$: $\min_{\bold x\in\oD_{\bm\oldb}}\oldb_1(\bold x)>0$.\\
\\
respectively called sub-critical, critical and super-critical condtions.
Throughout this  paper, we consider $\bm{\alpha}$ to be a fixed (given) parameter. Therefore,
unless strictly necessary (such as in Section \ref{sec:discussion} and in the proofs of Proposition \ref{prop:alpharegion} and Theorem \ref{prop:sub}),
we do not explicit the dependence of $\oldb_i$ on $\alpha_i$ and of $\oD_{\bm\oldb}$ on $\bm{\alpha}$.
Note that, for $k=1$, $\oD$ reduces to $[0,r/(r-1)]$ and vectorial function ${\bm\oldb}:=(\oldb_1,\ldots,\oldb_k)$ reduces to function $\oldb_1$
defined in \eqref{eq:b}. Therefore, in the case of one community, $({\mathcal Sub})$ reduces to $\alpha_1<1$, $({\mathcal Crit})$ reduces to $\alpha_1=1$ and $({\mathcal Sup})$ reduces to
$\alpha_1>1$.

The next theorems provide the main results of this paper.

\begin{Theorem}\label{teo-blocksub}
Assume \eqref{eq:a1a3}, \eqref{eq:hyp2bis}, \eqref{eq:hyp4}, \eqref{eq:gratiog}, \eqref{eq:trivial}, \eqref{eq:ala4},
\eqref{eq:alfaorder} and
$({\mathcal Sub})$. Then, for any $\varepsilon>0$ there exists $c(\varepsilon)\in\mathbb R_+$ such that
\begin{equation}\label{eq:limsub}
P\left(\Big|\frac{A_n^*}{\nac{1}{n}}-x_*\Big|>\varepsilon\right)=O(\mathrm{e}^{-c(\varepsilon)g_n^{(1)}}),
\end{equation}
where the explicit expression for constant $x_*$ is given by \eqref{eq:xstar}.
\end{Theorem}

\begin{Theorem}\label{teosupcrit}
Assume \eqref{eq:a1a3}, \eqref{eq:hyp2bis}, \eqref{eq:hyp4}, \eqref{eq:gratiog}, \eqref{eq:trivial}, \eqref{eq:ala4}, \eqref{eq:alfaorder}
and $({\mathcal Sup})$. Then, for any $\varepsilon>0$ there exists $c(\varepsilon)\in\mathbb R_+$ such that
\[
P\left(\Big|\frac{A_n^*}{n}-1\Big|>\varepsilon\right)=O(\mathrm{e}^{-c(\varepsilon)g_n^{(1)}}).
\]
\end{Theorem}
These results can be roughly summarized as follows:
$(i)$ under $({\mathcal Sub})$, the bootstrap percolation process on the SBM
basically does not evolve as $n\to\infty$ (note that, due to  conditions \eqref{eq:trivial}, \eqref{eq:ala4} and \eqref{eq:alfaorder},
$\sum_{i=1}^{k}a_n^{(i)}\sim\nac{1}{n}$ as $n\to\infty$), $(ii)$ under $({\mathcal Sup})$, the bootstrap percolation process on the SBM basically percolates the whole graph as $n\to\infty$.

Replacing assumption \eqref{eq:hyp2bis} with  (slightly) stronger condition
\begin{equation}\label{eq:hyp2bisstr}
1/\on{i}\ll \op{{i}}\lesssim 1/(\on{i})^\frac{1}{r'},\quad\text{$r'\in (r,\infty)$}
\end{equation}
the convergence in probability provided by Theorems \ref{teo-blocksub} and \ref{teosupcrit} can be strengthened to an almost sure convergence by a standard application of Borel-Cantelli lemma.
For completeness, we state such refinements in the next two corollaries (for which we omit the proofs).
\begin{Corollary}
Assume \eqref{eq:a1a3}, \eqref{eq:hyp4}, \eqref{eq:gratiog}, \eqref{eq:trivial}, \eqref{eq:ala4},
\eqref{eq:alfaorder}, \eqref{eq:hyp2bisstr} and $({\mathcal Sub})$. Then
\begin{equation*}
\lim_{n\to\infty}\frac{A_n^*}{\nac{1}{n}}=x_*\in\mathbb{R}_+\quad\text{almost surely.}
\end{equation*}
\end{Corollary}

\begin{Corollary}
Assume \eqref{eq:a1a3}, \eqref{eq:hyp4}, \eqref{eq:gratiog}, \eqref{eq:trivial}, \eqref{eq:ala4}, \eqref{eq:alfaorder}, \eqref{eq:hyp2bisstr}
and $({\mathcal Sup})$. Then
\[
\lim_{n\to\infty}\frac{A_n^*}{n}=1\quad\text{almost surely.}
\]
\end{Corollary}

\subsubsection{An informal description of some basic ideas of the proofs}

In broad terms, proofs of Theorems \ref{teo-blocksub} and \ref{teosupcrit} adopt the following
scheme: first, we analyze the average dynamics of the number of
\lq\lq active and not yet used" nodes; second, we show that the \lq\lq true" random processes
are sufficiently concentrated around their averages. As already noticed in Remark \ref{rem1}, the following issue makes the analysis of the bootstrap percolation
process significantly more complex  on the SBM than  on the Erd\"{o}s--R\'{e}nyi random graph: for any $1\leq i\leq k$ and $t<T_n$,
random variables $A_i(t)$ heavily depend on quantities $\{U_i(t)\}_{1\leq i\leq k}$, and therefore
on the considered \lq\lq strategy". In turn, the chosen \lq\lq strategy"
$\{C_n^{(i)}(t)\}_{1\le i\le k,t <T_n}$, according to which the nodes are selected and \lq\lq used"
in different communities,
is constrained by the availability of \lq\lq active and not yet used" nodes in the different communities,
indeed ${\color{black}A_n^{(i)}}(t-1)-{\color{black}U_n^{(i)}}(t-1)>0$ is a necessary condition for $C_n^{(i)}(t)>0$.
In this discussion, we refer to such constraints as \lq\lq feasibility\rq \rq constraints, and to
\lq\lq strategies" satisfying such constraints as \lq\lq feasible strategies".

We remark once again that, due to Proposition \ref{prop:equiv}, $A_n^*$ does not depend on the considered \lq\lq feasible" strategy.
Therefore, in order to study $A_n^*$ we have the freedom to arbitrarily select  the \lq\lq feasible strategy" according to our convenience.
A first crucial step in the proofs of our main results consists in the  identification  of a suitable \lq\lq feasible strategy",
which allows us to analyze the bootstrap percolation process. A first natural candidate is the \lq\lq strategy" defined, for $1\leq t\leq T_n$, by:
\begin{equation}\label{eq:strategy}
C_n^{(i)}(t):=\frac{{\color{black}A_n^{(i)}}(t-1)-{\color{black}U_n^{(i)}}(t-1)}{\sum_{j=1}^{k}({\color{black}A_n^{(j)}}(t-1)-{\color{black}U_n^{(j)}}(t-1))},
\quad\text{$i=1,\ldots,k$}
\end{equation}
which corresponds to select (and then to \lq\lq use") a node uniformly at random,
among all the \lq\lq active and not yet used" nodes.
The main drawback of this \lq\lq strategy" is that the analysis of corresponding process $\bold{U}_n(t)$,
$t\leq T_n$,
appears prohibitive due to its complex correlation structure.
To circumvent these difficulties we introduce a \lq \lq fluid" version of \lq\lq strategy" \eqref{eq:strategy}, under which $\bold{U}_n(t)$ is deterministic.
The construction of such \lq\lq fluid strategy" is related to the solution of a Cauchy problem that will be introduced in Subsection \ref{cauchy-prob}.
Unfortunately, since process $\bold{U}_{n}(t)$ is deterministic, regardless of the evolution of process $\bold{A}_n(t)$, one can not guarantee that
such \lq\lq fluid policy" is \lq\lq feasible" up to time $T_n$. Therefore, we limit the application of  such a  \lq\lq fluid strategy" up to a certain random time
$T'_n$,
where its \lq\lq feasibility" is guaranteed, and we switch to  \lq\lq policy" described by  \eqref{eq:strategy}  in the time interval $(T'_n,T_n]$.

The asymptotic analysis of the average dynamics of the above \lq\lq fluid
strategy" permits us to identify three regimes, which  are shown to be equivalent to
$({\mathcal Sub})$,  $({\mathcal Sup})$ and   $({\mathcal Crit})$.
Then we prove that the number of \lq\lq active and not yet used" nodes
concentrates around its average. To accomplish this step, differently from \cite{JLTV}, where martigales' theory has been employed,
we exploit classical concentration inequalities for the binomial distribution.

\section{Analytical description of the critical regions in terms of the parameter $\bm{\alpha}$}\label{sec:discussion}

By Theorems \ref{teo-blocksub} and \ref{teosupcrit} we immediately determine the sub-critical and the super-critical regions,
i.e.,  set of $\bm{\alpha}$ for which the sub-critical and the super-critical behavior is experienced.
Indeed, by expliciting the dependence of $\oldb_i$ on $\alpha_i$ and the dependence of $\oD_{\bm\oldb}$ on $\bm{\alpha}$,
we can associate to conditions $({\mathcal Sub})$, $({\mathcal Crit})$ and $ ({\mathcal Sup})$ the regions:\\
\begin{equation}
\mathcal{R}_{Sub}:=\{\bm{\alpha}:\,\,\min_{\bold x\in\oD_{{\bm\oldb}, \bm{\alpha}}}\oldb_1(\bold{x}, \alpha_1)<0\},\quad
\mathcal{R}_{Crit}:=\{\bm{\alpha}:\,\,\min_{\bold x\in\oD_{{\bm\oldb}, \bm{\alpha}}}\oldb_1(\bold{x}, \alpha_1)=0\}, \label{Rcond-sub-crit}
\end{equation}
\begin{equation}
\mathcal{R}_{Sup}:=\{\bm{\alpha}:\,\,\min_{\bold x\in\oD_{{\bm\oldb}, \bm{\alpha}}}\oldb_1(\bold x,\alpha_1)>0\},\label{Rcond-sup}
\end{equation}
respectively. Since $\mathcal{R}_{Sub}$, $\mathcal{R}_{Crit}$ and $\mathcal{R}_{Sup}$
are disjoint and exhaustive, and $\mathcal{R}_{Crit}$ contains only those points which are in the boundary
of $\mathcal{R}_{Sup}$ and in the boundary of $\mathcal{R}_{Sub}$, our description of regions $\mathcal{R}_{Sup}$ and $\mathcal{R}_{Sub}$ is tight.

We have already noticed that in the case of one community (i.e., for the Erd\"{o}s--R\'{e}nyi random graph)
$\mathcal{R}_{Sub}\equiv [0,1)$, $\mathcal{R}_{Crit}\equiv\{1\}$ and $\mathcal{R}_{Sup}\equiv (1,\infty)$.
Unfortunately, in the more general case of $k\ge 2$ communities, sub-critical and super-critical regions
can not be always described in terms of  parameter $\bm{\alpha}$ through simple closed form expressions.
However, as we shall see in the next subsection, when matrix $\bm{\chi}$ is invertible
 critical region can be determined by a simple computational procedure.

\subsection{General results}

The following propositions hold.

\begin{Proposition}\label{prop:alphagreaterthan1}
Assume $k\geq 2$, \eqref{eq:a1a3}, \eqref{eq:hyp2bis}, \eqref{eq:hyp4}, \eqref{eq:gratiog}, \eqref{eq:trivial}, \eqref{eq:ala4}, and
\eqref{eq:alfaorder} with $\alpha_1\geq 1$. Then $({\mathcal Sup})$ holds.
\end{Proposition}

\begin{Proposition}\label{prop:alpharegion}
Assume $k\geq 2$, \eqref{eq:a1a3}, \eqref{eq:hyp2bis}, \eqref{eq:hyp4}, \eqref{eq:gratiog}, \eqref{eq:trivial}, \eqref{eq:ala4} and
\eqref{eq:alfaorder}.  Under the additional assumption that matrix $\bm{\chi}$ is invertible, we have
\begin{align*}
\mathcal{R}_{Crit} &=\{\bm{\alpha}\in [0,1]^k:\,\,\alpha_i=(\bold{x}_{\bm{\theta}}(\bm{\chi}^{-1})^t)_i-r^{-1}(1-r^{-1})^{r-1}(\bold{x}_{\bm{\theta}})_i^{r},
\,\,\bm{\theta}:=(1,\theta_2,\ldots,\theta_{k})\in\mathbb{R}_{+}^{k}\},
\end{align*}
where $\bold{x}_{\bm{\theta}}$ is the (row) vector with components $\frac{r}{r-1}\left(\bm{\theta}\bm{\chi}^{-1}\mathrm{diag}(\bm{\theta}^{-1})
\right)_i^{1/(r-1)}$, $i=1,\ldots,k$, $\mathrm{diag}(\bm{\theta}^{-1})$ is the diagonal $k\times k$ matrix with diagonal elements $1,\theta_2^{-1},\ldots,\theta_k^{-1}$,
$(\bm{\chi}^{-1})^t$ is the transpose of matrix $\bm{\chi}^{-1}$, and
\begin{align*}
\mathcal{R}_{Sub}&=\{\bm{\alpha}\in [0,1]^k:\,\,\{\theta\bm{\alpha}:\,\,\theta\in [0,1]\}\cap \mathcal{R}_{Crit}=\emptyset\},\\
\mathcal{R}_{Sup}&=\{\bm{\alpha}\in [0,\infty)^k:\,\,\{\theta\bm{\alpha}:\,\,\theta\in [0,1)\}\cap\mathcal{R}_{Crit}\neq\emptyset\},
\end{align*}
\end{Proposition}

\subsection{The case of $k$ \lq\lq identical" communities: numerical illustrations and informal discussion}

We illustrate Proposition \ref{prop:alpharegion} in the special case of $k\geq 2$ \lq\lq identical" communities, i.e., $\on{1}=\ldots=\on{k}$,
$\oq{{ij}}=q_n$, for $i\neq j$, and $\op{i}=p_n$ for any $i$.

We assume: $k/n\ll p_n\ll (k/n)^{1/r}$, $q_n<p_n$, $q_n/p_n\to\psi\in (0,1)$, \eqref{eq:trivial} and
\eqref{eq:alfaorder}. Then $\chi_{ii}=1$, $\chi_{ij}=\chi_{ji}=\psi$ for $i\neq j$ and the corresponding assumptions
\eqref{eq:a1a3}, \eqref{eq:hyp2bis}, \eqref{eq:hyp4}, \eqref{eq:gratiog} and \eqref{eq:ala4} are satisfied. A straightforward computation shows that
\begin{equation*}
{\bm \chi^{-1} }  =  \left( \begin{array}{cccc}
a & a' & \cdots & a' \\
a' & a & \cdots & a' \\
\vdots & \vdots & \ddots & \vdots \\
a' & a' & \cdots & a
\end{array} \right)
\end{equation*}
where $a:=\frac{1+(k-2)\psi}{1+(k-2)\psi-(k-1)\psi^2}>0$, $a':=\frac{-\psi}{1+(k-2)\psi-(k-1) \psi^2}<0$ and, for $\bm{\theta}:=(1,\theta_2,\ldots,\theta_{k})\in\mathbb{R}_+^k$,
\[
\bm{\theta}\bm{\chi}^{-1}\mathrm{diag}(\bm{\theta}^{-1})=\left(a+a'\sum_{i=2}^{k}\theta_i,a+a'\sum_{i=3}^{k}\frac{\theta_i}{\theta_2}+\frac{a'}{\theta_2},\ldots,
a+a'\sum_{i\neq j}^{2,k}\frac{\theta_i}{\theta_j}+\frac{a'}{\theta_j},\ldots,
a+a'\sum_{i=2}^{k-1}\frac{\theta_i}{\theta_k}+\frac{a'}{\theta_k}\right).
\]

\sidebyside{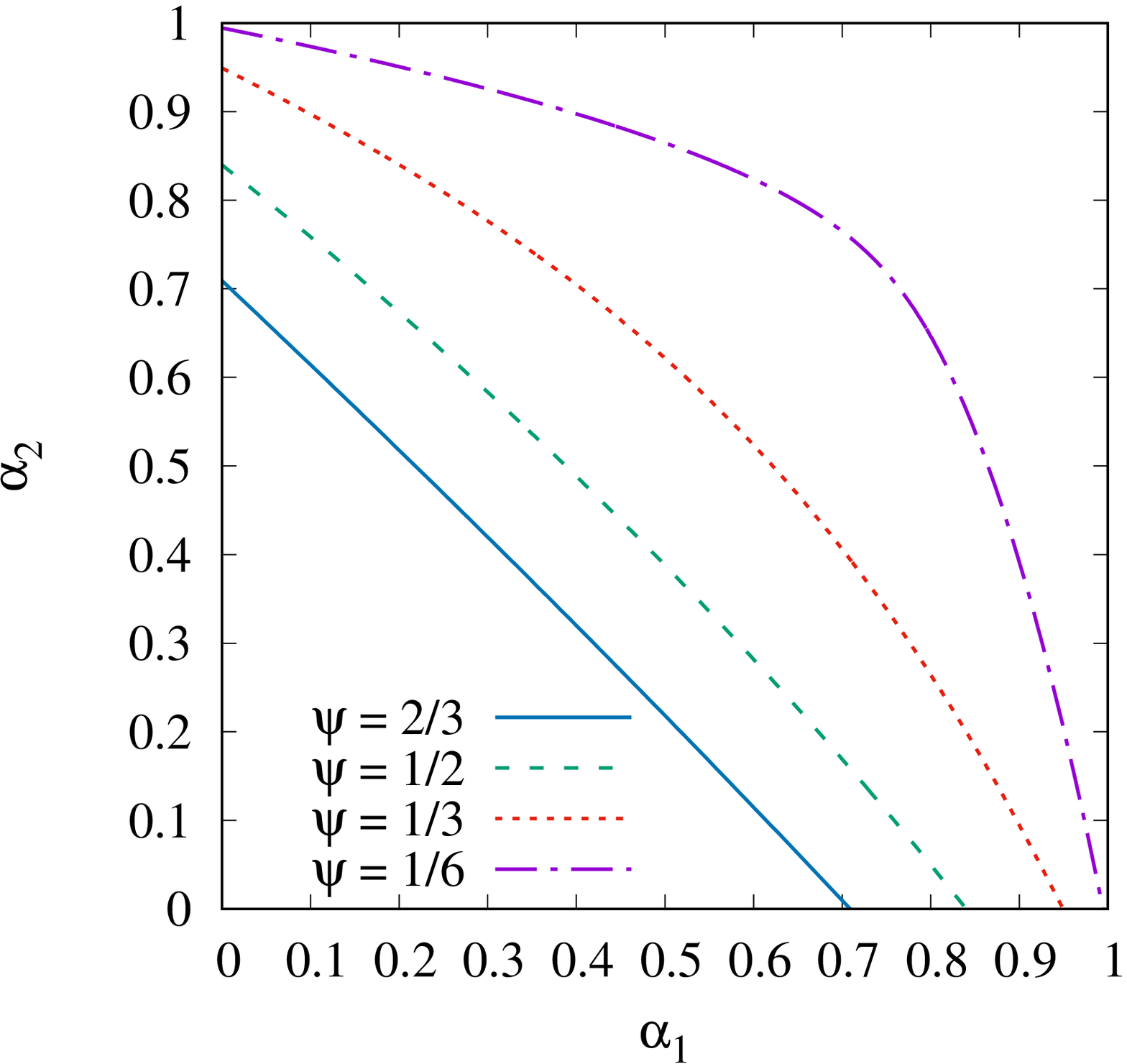}{Critical curves $\mathcal{R}_{Crit}$ for $k=r=2$ and different values of $\psi$.}
{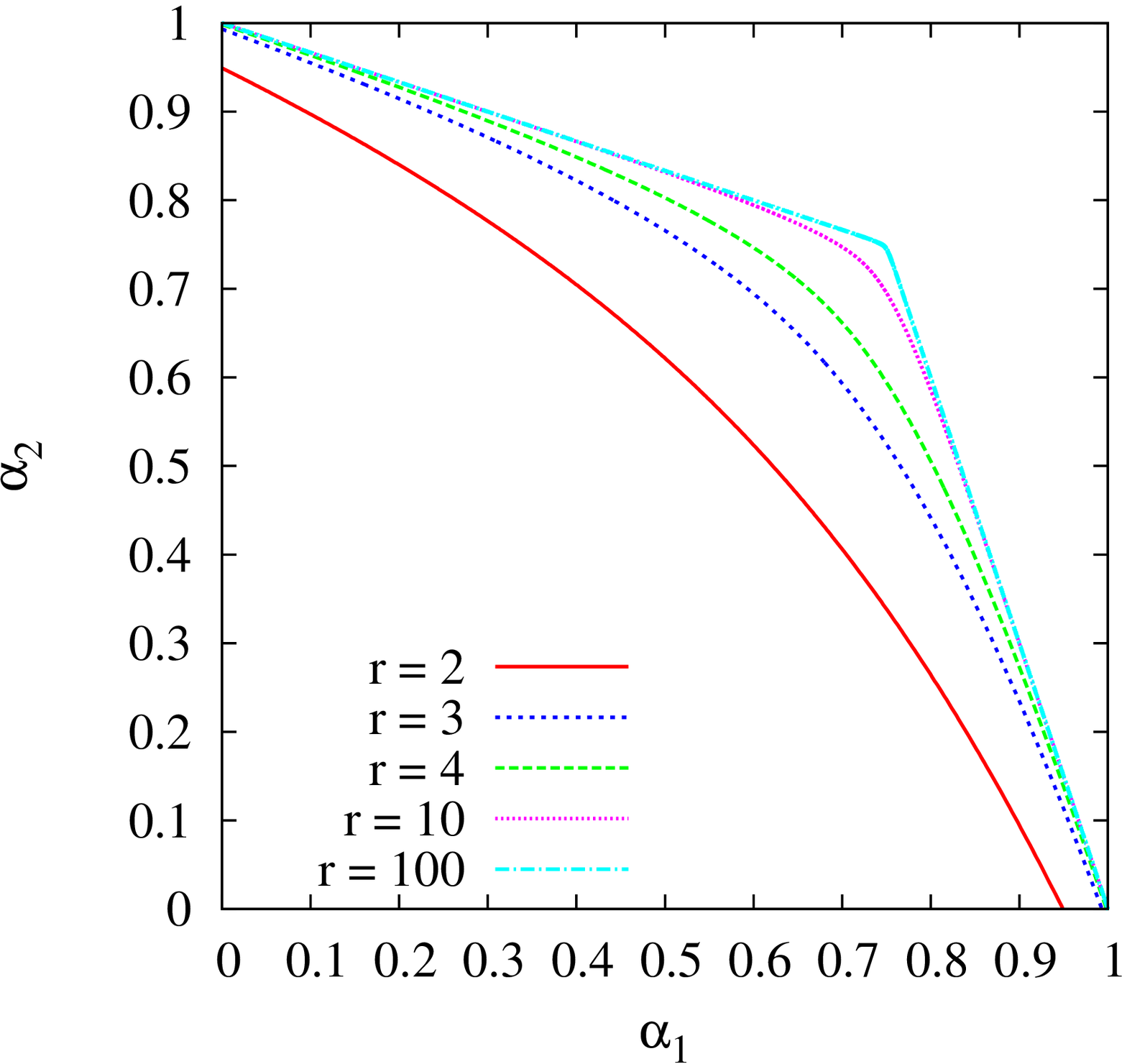}{Critical curves $\mathcal{R}_{Crit}$ for $k=2$, \mbox{$\psi=1/3$} and different values of $r$.}

In Figure \ref{fig:omega}, for fixed $k=r=2$, we plotted  curve $\mathcal{R}_{Crit}$
for different values of  parameter $\psi$. Note that, as $\psi\downarrow 0$,  sub-critical region approaches the whole square.
Indeed, this can be formally verified letting $\psi$ tend to zero in the corresponding expression of $\mathcal{R}_{Crit}$.

In Figure \ref{fig:r} for fixed $k=2$ and $\psi=1/3$, we plotted curve $\mathcal{R}_{Crit}$ for different values of parameter $r$.
Note that,
as $r\uparrow\infty$, the sub-critical region approaches the domain
\[
\{(\alpha_1,\alpha_2)\in [0,1]^2:\,\,\alpha_1+\psi\,\alpha_2<1,\,\,\alpha_2+ \psi\,\alpha_1<1\}.
\]

The above numerical illustrations show that the best way to trigger percolation
of the whole graph is to maximize the value of either $\alpha_1$ or $\alpha_2$, i.e., to put all the seeds in the same community; instead, the worst manner to trigger percolation
of the whole graph is to jointly minimize the values of $\alpha_1$ and $\alpha_2$, i.e., to equally partition  seeds between the communities.

These properties of the bootstrap percolation process on the SBM may be proved formally (for the general case of $k$ identical communities) by exploiting the convexity of the sub-critical region
and the symmetry among  communities. We remark that,
while for the worst way to trigger percolation, critical value $\bm{\alpha}$ can be computed analytically, for
the best way to trigger percolation one can only numerically estimate $\bm{\alpha}$.

More precisely, consider firstly the worst way to trigger percolation, i.e.,
assume \eqref{eq:trivial} with
\[
\alpha=\alpha_1=\ldots=\alpha_k>0.
\]
A straightforward computation shows that $(\alpha,\ldots,\alpha)\in\mathcal{R}_{Crit}$
if and only if
\[
\alpha=\left(\frac{1-\psi}{1+(k-2)\psi-(k-1)\psi^2}\right)^{r/(r-1)}=(1+(k-1)\psi)^{-r/(r-1)}.
\]
Since each community has $n/k$ nodes, for each community we have
\[
g_n^{(i)}=\nan=\left(1-\frac{1}{r}\right)\left(\frac{k (r-1)!}{n p_n^r}\right)^{1/(r-1)}.
\]
Therefore, putting the same number of seeds in each community, asymptotically,
the total number of seeds is:
\begin{equation*}
k \,\alpha \, \nan = \left(1-\frac{1}{r}\right) \left(\frac{k}{1 + (k-1) \psi}\right)^{r/(r-1)} \left(\frac{(r-1)!}{n p_n^r}\right)^{1/(r-1)}.
\end{equation*}
As expected, as $\psi\to 1$, quantity $k\alpha \nan$ tends to $\nan/k^{1/(r-1)}$, i.e., the total number of seeds in the case of one community.
Note that for $\psi < 1$, as $k$ grows large, $k\alpha \nan$ tends to
\[
\left(1-\frac{1}{r}\right)\psi^{-r/(r-1)} \left(\frac{(r-1)!}{n p_n^r}\right)^{1/(r-1)}
\]
and this quantity is asymptotically equivalent (as $n\to\infty$) to
\[
\left(1-\frac{1}{r}\right)\left(\frac{(r-1)!}{n q_n^r}\right)^{1/(r-1)},
\]
i.e, the total number of seeds in the case of one community with connection probability $q_n$.
The intuition behind this fact is the following. As $k$ grows large, most of the
neighbors of any given node belong to communities which are different from the community to which the node belongs. So the impact of neighbors belonging
to the same community of the node becomes negligible.

Now consider the best way to trigger percolation, i.e., assume $\bm{\alpha}=(\alpha,0,\ldots,0)$. In this case, applying Proposition \ref{prop:alpharegion} one has that, in general, $\alpha$
can only be numerically estimated. Indeed, $\alpha$ can be defined as the unique zero in $(0,1)$
of an algebraic equation of order $r^2-1$. In the special case of $r=2$, a closed form expression for $\alpha$ can be obtained
by using Cardano's formulas.

In Figure \ref{fig:kvary} we have plotted the total critical number of seeds in the SBM
normalized to the critical number of seeds in $G(n,p_n)$, as a function of $k$, for fixed $\psi=1/10$,
and either $r=2$ or $r=4$. The curves related to the extreme cases
are labeled \lq\lq equal-split" (when the seeds are equally divided among the communities)
and \lq\lq all-in-one" (when all the seeds are put in the same community), respectively.
We note that, although both seeds' allocations strategies
tend to perform the same as $k$ grows large, they require a quite different number of seeds
for small number of communities. The performance of any seeds' allocation strategy lies
between the two extreme cases.

\tgifeps{9}{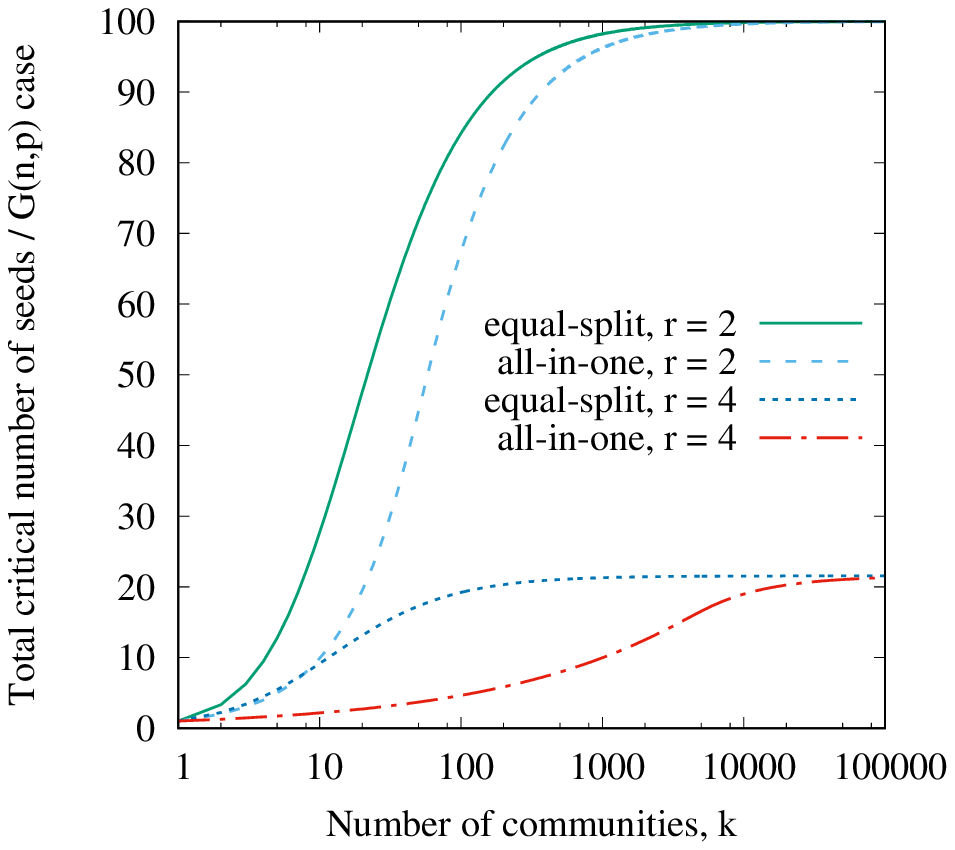}{Performance of the extreme seeds' allocation strategies as function of the number
$k$ of the communities, for $\psi=1/10$ and either $r=2$ or $r=4$.}

\section{Proofs}\label{sec:proofs}

\subsection{Preliminaries}

We start with some preliminaries. In Subsection \ref{sect:jacobian} we highlight some properties of the
Jacobian of ${\bm\oldb}$ which will be used in the proofs. In Subsection \ref{sect:addnot} we introduce some additional notation.
In Subsection \ref{cauchy-prob} we describe the Cauchy problem whose solution allows us to define the asymptotic \lq\lq fluid
strategy" previously mentioned. Finally, in Subsection \ref{sect:equiv}
we present equivalent formulations of conditions $({\mathcal Sub})$, $({\mathcal Crit})$ and $({\mathcal Sup})$, which will be exploited in the proof of the main theorems.
Throughout this section we assume model assumptions
\eqref{eq:a1a3}, \eqref{eq:hyp2bis}, \eqref{eq:hyp4}, \eqref{eq:gratiog}, \eqref{eq:trivial}, \eqref{eq:ala4} and \eqref{eq:alfaorder}.

\subsubsection{The Jacobian of $\bm\oldb$}\label{sect:jacobian}

In the following, we denote by $\overset\circ{\mathcal B}$
and $\partial{\mathcal B}$ the interior
and the boundary of a Borel set ${\mathcal B}\subseteq\mathbb{R}^k$,
respectively, and refer the reader to \cite{LT} for standard notions and results of matrix theory. Moreover, we shall not specify, when this is clear from the context,
if a vector is meant as column or row vector.

We have:
\[
\frac{\partial \oldb_i(\bold x)}{\partial x_i}=-1+(1-r^{-1})^{r-1}\left(x_i+\sum_{j\neq i}^{1,k}\chi_{ij}x_j\right)^{r-1}
\]
and
\[
\frac{\partial \oldb_i(\bold x)}{\partial x_j}=\chi_{ij}(1-r^{-1})^{r-1}\left(x_i+\sum_{h\neq i}^{1,k}\chi_{ih}x_h\right)^{r-1},\quad j\neq i.
\]
So
\[
\frac{\partial \oldb_i(\bold x)}{\partial x_i}<0\quad\text{if and only if}\quad x_i+\sum_{j\neq i}^{1,k}\chi_{ij}x_j<\frac{r}{r-1},\quad\text{and}\quad
\frac{\partial \oldb_i(\bold x)}{\partial x_j}\geq 0\quad\text{for any $j\neq i$.}
\]
Let
\[
J_{{\bm\oldb}}(\bold{x}):=\left(\begin{array}{ccc}
\frac{\partial \oldb_1}{\partial x_1}(\bold{x}) &\ldots  & \frac{\partial \oldb_1}{\partial x_k}(\bold{x})\\
\frac{\partial \oldb_2}{\partial x_1}(\bold{x}) &\ldots  & \frac{\partial \oldb_2}{\partial x_k}(\bold{x})\\
 \vdots & \vdots  & \vdots\\
\frac{\partial \oldb_k}{\partial x_1}(\bold{x}) &\ldots  & \frac{\partial \oldb_k}{\partial x_k}(\bold{x})\\
\end{array}\right),\quad\text{$\bold x\in\oD$}
\]
be the Jacobian matrix of vectorial function $\bm\oldb$ and define
\begin{equation*}
\ocalD:=\bigcup_{i=1}^{k}\left\{(x_1,\ldots, x_k)\in\oD:\,\,x_i+\sum_{j\neq i}^{1,k}\chi_{ij}x_j=\frac{r}{r-1}\right\}.
\end{equation*}
We remark that $(i)$ diagonal terms of $J_{{\bm\oldb}}(\bold{x})$
are such that $(J_{{\bm\oldb}}(\bold{x}))_{ii}\leq 0$, $i=1,\ldots,k$, and the equality holds
if and only if $\bold x\in\ocalD$, $(ii)$  terms outside the diagonal of $J_{{\bm\oldb}}(\bold{x})$ are such that
$(J_{{\bm\oldb}}(\bold{x}))_{ij}\geq 0$, $i,j=1,\ldots,k$, $i\neq j$, and the equality holds
if and only if $\bold x=\bold 0$ and/or $\chi_{ij}=0$.

Taking $\ell(\bold x)>\max_{1\leq i\leq k}|(J_{{\bm\oldb}}(\bold{x}))_{ii}|$, we rewrite the Jacobian matrix of $\bm\oldb$ as
$J_{{\bm\oldb}}(\bold{x})=\mathrm P(\bold x)-\ell(\bold x)\mathrm{I}$, where $\mathrm P(\bold x)$ is the non-negative matrix defined by
$(\mathrm P(\bold x))_{ii}:=\ell(\bold x)+(J_{{\bm\oldb}}(\bold{x}))_{ii}$, $i=1,\ldots,k$,
$(\mathrm P(\bold x))_{ij}:=(J_{{\bm\oldb}}(\bold{x}))_{ij}$, $i,j=1,\ldots,k$, $i\neq j$, and $\mathrm I$ is the identity matrix.
Note that by the irreducibility of $\bm{\chi}$ (i.e., condition \eqref{eq:ala4}) we immediately have that $\mathrm{P}(\bold x)$
is irreducible as well. Combining this with the fact that all  diagonal terms of $\mathrm{P}(\bold x)$ are strictly positive
we deduce that $\mathrm{P}(\bold x)$ is primitive. Therefore, by  Perron-Fr\"obenius theorem there exists a positive eigenvalue
$\tilde{\lambda}_{PF}(\bold x)>0$ of $\mathrm P(\bold x)$ to which it corresponds a unique positive (component-wise) eigenvector $\tilde{\bold{\ov}}_{PF}(\bold x)$. A straightforward computation shows that
\begin{equation*}
\lambda_{PF}(\bold x):=\tilde{\lambda}_{PF}(\bold x)-\ell(\bold x)
\end{equation*}
is an eigenvalue of $J_{{\bm\oldb}}(\bold x)$ to which it corresponds  (positive) eigenvector
$\bold{\ov}_{PF}(\bold x):=\tilde{\bold{\ov}}_{PF}(\bold x)$.
\textcolor{black}{Since
$\mathrm P(\bold x)$ and $J_{{\bm\oldb}}(\bold{x})$ have the same eigenvectors, it turns out that  $\bold{\ov}_{PF}(\bold x)$
is the unique positive eigenvector of $J_{{\bm\oldb}}(\bold{x})$.}

\subsubsection{Some more notation}\label{sect:addnot}

We denote by $\mathcal{K}_h$ the set of all communities at distance $h\ge 1$ from
community $G_1$ on  connected graph $\widetilde{G}$,
and by ${\color{black}\overline{d}}$ the maximum distance from community $G_1$, i.e.,
\begin{equation*}
\mathcal{K}_0:=\{1\},\quad \mathcal{K}_1:=\{i\in\{2,\ldots,k\}:\,\,(\bm{\chi}^h)_{1i}>0\},
\end{equation*}
\begin{equation}\label{eq:mathcalk}
\mathcal{K}_h:=\{i\in\{2,\ldots,k\}:\,\,(\bm{\chi}^h)_{1i}>0\,\,and\,\,(\bm{\chi}^s)_{1i}=0\,\,\forall\,\,1\leq s\leq h-1\},\quad h\geq 2
\end{equation}
where $\bm{\chi}^h$ is the $h$th power of matrix $\bm{\chi}$, and
\begin{equation}\label{eq:maxdist}
{\color{black}\overline{d}}:=\max\{h:\,\,\mathcal{K}_h\neq\emptyset\}.
\end{equation}

For later purposes (see e.g. Lemma \ref{le:initialpoint}) we note that, by the irreducibility of $\bm{\chi}$,
for any $i\in\{1,\ldots,k\}$, there exists a unique ${\color{black}d_i}\in\{0,\ldots,{\color{black}\overline{d}}\}$ such that $i\in\mathcal{K}_{{\color{black}d_i}}$.
Quantity ${\color{black}d_i}$ is nothing but that the distance,
on the connected graph $\widetilde{G}$, between  communities $G_i$ and $G_1$.

We put
\[
\bold{1}_{\mathcal{K}_h}:=(\bold{1}(1\in\mathcal{K}_h),\ldots,\bold{1}(k\in\mathcal{K}_h)),\quad h=0,\ldots,{\color{black}\overline{d}}
\]
and
\begin{equation}\label{eq:beta}
\beta_h:=\frac{r^{-1}(1-r^{-1})^{r-1}}{2}(\chi\beta_{h-1} )^r,\quad 1\leq h\leq{\color{black}\overline{d}},\quad\beta_0:=\frac{\alpha_1}{2}
\end{equation}
where
\begin{equation}\label{eq:chimin}
\chi:=\min\{\chi_{ij}:\chi_{ij}>0\}.
\end{equation}
Finally, we define  vectors
\begin{equation}\label{eq:x00}
\bold{x}_0^{(0)}:=(\alpha_1/2,0,0,\ldots,0)=\beta_0\bold{1}_{\mathcal{K}_0}
\end{equation}
and
\begin{equation}\label{eq:x0h}
\bold{x}_0^{(h)}=\sum_{s=0}^{h}\beta_{s}\bold{1}_{\mathcal{K}_{s}},\quad 1\le h \le {\color{black}\overline{d}}.
\end{equation}

\subsubsection{A related Cauchy problem}\label{cauchy-prob}

Hereon, we denote by
$\mathrm{B}_{\delta}(\bold x)$ the
closed ball of $\mathbb{R}^k$ centered at $\bold x\in\mathbb R^k$ with radius $\delta>0$
and set
\begin{equation}\label{eq:x0}
\bold{x}_{0}:=\bold{x}_0^{({\color{black}\overline{d}})}.
\end{equation}
From now on, we suppose $\oldb_i:\mathrm{B}_{\mathrm{diam}(\oD)}(\bold{x}_0)\subset\mathbb R^k\to\mathbb R$, $i=1,\ldots,k$, i.e., the functions $\oldb_i$ are defined on the closed ball
$\mathrm{B}_{\mathrm{diam}(\oD)}(\bold{x}_0)$, where
$\mathrm{diam}(\oD)$ is the diameter of $\oD$. We consider the Cauchy problem defined by:
\begin{equation}\label{eq:cauchy}
\bold{x}'(y)={\bm\oldb}(\bold{x}(y)),\quad\bold{x}(y_0):=\bold{x}_0.
\end{equation}
By the regularity properties of  functions $\oldb_i$
there exists a unique solution $\bold{x}(y)$, $y\in [y_0,y_+)$, of \eqref{eq:cauchy} (see \cite{Arnold}), where
\[
y_+:=\sup\{y>y_0:\,\text{$\bold{x}(y)\in\mathrm{B}_{\mathrm{diam}(\oD)}(\bold{x}_0)$}\}.
\]
We remark that, in general, $y_+$ can be either finite or infinite. Consider the subset of $\oD$ defined by
\[
\oE':=\{\bold x\in\oD:\,\,{\bm\oldb}(\bold{x})\geq\bold 0\}.
\]
Throughout this paper, we denote by $\oE$ the closure of the largest connected component of $\oE'$ containing $\bold{x}_0$
(as we shall explain later on, we have $\bold 0\in\oE$)
and by $\ote$ the first exit time of the solution $\bold{x}(y)$, $y\in [y_0,y_+)$, of \eqref{eq:cauchy}
from $\overset\circ{\oE}$, i.e.,
\begin{equation}\label{eq:tE}
\ote:=\inf\{y\in (y_0,y_+):\,\,\bold{x}(y)\notin\overset\circ{\oE}\},\quad \ote:=+\infty\quad\text{if $\{\cdots\}=\emptyset$.}
\end{equation}
We anticipate that by Lemma \ref{le:initialpoint} it holds $\bold{x}_0\in\overset\circ{\oE}$.
Hereafter, we shall consider  point
\begin{equation}\label{eq:x(tE)}
\bold{x}(\ote):=\lim_{y_0<y\uparrow \ote}\bold{x}(y)\in\partial\oE
\end{equation}
(note that limit exists since $\bold{x}(\cdot)$ is strictly increasing (component-wise) on $\oE$).

\begin{Remark}
The following intuitive interpretation can be  given
to the solution of the Cauchy problem \eqref{eq:cauchy}.
Since
\[
\frac{\bold{x}'(y)}{\sum_{i=1}^{k}|x_i'(y)|}=\frac{{\bm\oldb}(\bold{x}(y))}{\sum_{i=1}^{k}|\oldb_i(\bold{x}(y))|},\quad\text{$y\in [y_0,y_+)$}
\]
where $\bold{x}'(y):=(x'_1(y),\ldots,x'_k(y))$,
the quantity $\frac{x'_i(y)}{\sum_{i=1}^{k}|x_i'(y)|}$ can
be interpreted as the (normalized) rate at which the active nodes
are \lq\lq explored" in community $G_i$. By Lemma \ref{le:bas2}
\[
\frac{\rho_i(\bold{x}(y))}{\sum_{i=1}^{k}|\oldb_i(\bold{x}(y))|}\sim_e \frac{E[{\color{black}A_n^{(i)}}(t)-{\color{black}U_n^{(i)}}(t)\mid\bold{U}_n(t)=\lfloor\bold{x}(y)g_n\rfloor]}
{\sum_{i=1}^{k}(E[{\color{black}A_n^{(i)}}(t)-{\color{black}U_n^{(i)}}(t)\mid\bold{U}_n(t)=\lfloor\bold{x}(y)g_n\rfloor])},\quad \text{with }
t=\sum_i^k \lfloor x_i(y)g_n \rfloor.
 \]
Therefore, the solution of Cauchy problem \eqref{eq:cauchy}
corresponds (asymptotically) to  (normalized) trajectory $t\to\bold{U}_{n}(t)$ described by a \lq\lq policy"
according to which, given $\bold{U}_n(t-1)$,
nodes to be \lq\lq explored" at time step $t$ are chosen in community $G_i$
with probability
\[
\frac{E[{\color{black}A_n^{(i)}}(t-1)\mid\bold{U}_n(t-1)]-{\color{black}U_n^{(i)}}(t-1)}{\sum_{i=1}^{k}
(E[{\color{black}A_n^{(i)}}(t-1)\mid\bold{U}_n(t-1)]-{\color{black}U_n^{(i)}}(t-1))}.
\]
\end{Remark}

\subsubsection{Equivalent formulations of $({\mathcal Sub})$, $({\mathcal Crit})$ and $({\mathcal Sup})$}\label{sect:equiv}

The proofs of Theorems \ref{teo-blocksub} and \ref{teosupcrit} rely on an equivalent formulation of conditions $({\mathcal Sub})$, $({\mathcal Crit})$, $({\mathcal Sup})$
(denoted by $(\bold{Sub})'$, $(\bold{Crit})'$, $(\bold{Sup})'$, respectively)
expressed in terms of  properties of the solution of  Cauchy problem \eqref{eq:cauchy}.
In turn, the proof of such equivalence
makes use of a third equivalent formulation of these conditions (denoted by $(\bold{Sub})$, $(\bold{Crit})$, $(\bold{Sup})$, respectively).

We recall that, by definition, a continuous and non-decreasing curve is a mapping
$\zeta: [\oc{1},\oc{2}]\to\mathbb{R}^k$, $\oc{1},\oc{2}\in[0,+\infty]$, $\zeta(\gamma):=(\zeta_1(\gamma),\ldots,\zeta_k(\gamma))$,
which is continuous and component-wise non-decreasing. We call $\zeta([\oc{1},\oc{2}])$ the trace of the curve.

Setting $\bold{x}_0^{(-1)}:=\bold{0}$, for $1\le h \le {\color{black}\overline{d}}$, we define
\[
\mathcal{S}_h:=\{\theta\bold{x}_0^{(h)}+(1-\theta)\bold{x}_0^{(h-1)}:\,\,0\le\theta\le 1\}\quad\text{and}\quad
\mathcal{S}:=\bigcup_{1\leq h\leq{\color{black}\overline{d}}}\mathcal{S}_h.
\]
Note that $\mathcal{S}\subset\oE$, indeed $\bold{x}_0\in\mathcal{S}$, $\mathcal S$ is clearly connected
and it can be checked that functions $\oldb_i$ are non-negative on $\mathcal S$
(as a straightforward consequence of Lemma \ref{le:initialpoint}, to be stated in the next subsection).
Let $\mathcal{C}_{\oE}$ be the trace of  curve
$
\bold{x}:[y_0,\ote]\to\mathbb R^k$ (where $\bold x(\cdot)$ denotes the solution of \eqref{eq:cauchy}),
define
\[
\tilde{\mathcal{S}}:=\{\bold x\in\oD:\,\, \bold{x}=\bold{x}(\ote)+\theta\bold{\ov}_{PF}(\bold{x}(\ote)),\,\theta\geq 0\}
\]
and let $\mathcal Z$ denote the set of distinct zeros of $\bm\oldb$ within $\oD$. We shall consider conditions:
\\
\\
\noindent$(\bold{Sub})$: There exists $\bold x\in\oD$ such that ${\bm\oldb}(\bold x)<\bold 0$.\\
\noindent$(\bold{Crit})$: $\mathcal{Z}=\{\bold z\}$, $\bold{z}\in\overset\circ{\oD}$, $\mathrm{det}J_{{\bm\oldb}}(\bold z)=0$ and $\lambda_{PF}(\bold z)=0$.\\
\noindent$(\bold{Sup})$: ${\bm\oldb}({\bf x})\neq\bold 0$, $\forall$ $\bold x\in\oD$.\\
\noindent$(\bold{Sub})'$: $\ote=+\infty$, $\bold{x}(\ote)\in\mathcal{Z}\cap\overset\circ{\oD}$,
the curve with trace $\mathcal{C}:=\mathcal{S}\cup\mathcal{C}_\oE\cup\tilde{\mathcal S}$ is continuous, non-decreasing and satisfies
\begin{equation}\label{eq:intersect_k}
\text{$\mathcal C\subset\oD$, $\bold{0}\in\mathcal{C}$, $\mathcal{C}\cap\ocalD\neq\emptyset$}
\end{equation}
and
\begin{equation}\label{sub-cond_k}
\min_{\bold x\in\mathcal{C}_{\oE}\cup\tilde{\mathcal{S}}}\min_{1\leq i\leq k} \oldb_i(\bold x)<0.
\end{equation}
\noindent$(\bold{Crit})'$: $\ote=+\infty$, $\bold{x}(\ote)\in\mathcal{Z}\cap\overset\circ{\oD}$,
the curve with trace $\mathcal{C}:=\mathcal{S}\cup\mathcal{C}_\oE\cup\tilde{\mathcal S}$ is continuous, non-decreasing and satisfies
\eqref{eq:intersect_k} and
\begin{equation}\label{crit_k}
\min_{\bold x\in\mathcal{C}_{\oE}\cup\tilde{\mathcal{S}}}\min_{1\leq i\leq k} \oldb_i(\bold x)=0.
\end{equation}
\noindent$(\bold{Sup})'$: $\ote<+\infty$, $\bold{x}(\ote)\in\ocalD$, ${\bm\oldb}(\bold{x}(\ote))>\bold 0$,
the curve with trace $\mathcal{C}:=\mathcal{S}\cup\mathcal{C}_\oE$ is continuous, non-decreasing and satisfies \eqref{eq:intersect_k} and
\begin{equation}\label{sup-cond_k}
\min_{\bold x\in\mathcal{C}_{\oE}} \min_{1\leq i\leq k} \oldb_i(\bold x)>0.
\end{equation}
\\

The following theorem (whose proof is given in Subsection \ref{subsec:thm61}) holds.

\begin{Theorem}\label{prop:sub}
Assume \eqref{eq:a1a3}, \eqref{eq:hyp2bis}, \eqref{eq:hyp4}, \eqref{eq:gratiog}, \eqref{eq:trivial}, \eqref{eq:ala4}, \eqref{eq:alfaorder}. Then:\\
\noindent$(i)$ $({\mathcal Sub})$, $(\bold{Sub})$ and $(\bold{Sub})'$ are equivalent.\\
\noindent$(ii)$ $({\mathcal Sup})$, $(\bold{Sup})$ and $(\bold{Sup})'$ are equivalent.\\
\noindent$(iii)$ $({\mathcal Crit})$, $(\bold{Crit})$ and $(\bold{Crit})'$ are equivalent.
\end{Theorem}

\subsection{Proof of Theorem \ref{teo-blocksub}}

\begin{table}
 \begin{center}
\caption{Some notation used in the proof of Theorem \ref{teo-blocksub}}
 \begin{tabular}{|l|l|l|}
  \hline
\multicolumn{3}{|c|}{Special instants}\\
\hline
Symbol & Mathematical definition & Remarks\\
 \hline
$t_n^{(s)}$  & $t_n^{(s)}\in\otaun([0,1])$ \,\, (see below \eqref{eq:tau}) & $t_n^{(s)}< t_n^{(s+1)} $,  $ 0\le s \le m+1$\\
$\hat{t}_n^{(h)}$ & $\sum_{i=1}^{k}\lfloor\sum_{\ell=0}^{h}\beta_\ell\bold{1}(i\in\mathcal{K}_\ell)\nac{i}{n}\rfloor$ \,\, (see \eqref{eq:tmin})&
$h\in\{0,\ldots,{\color{black}\overline{d}}\}$\\
$T''_n$ & $\min\{1\leq t\leq\otaun(1):\,{\color{black}A_n^{(i)}}(t-1)<\ozeta_n^{(i)}(t)\}$ \,\, (see \eqref{eq:tnprimo}) & \\
$T'_n$ &  $ T''_n-1$ \,\, (see \eqref{eq:tnprimo})& \\
\hline
\multicolumn{3}{|c|}{Special points and sets}\\
\hline
 Symbol & Mathematical definition & Remarks\\
 \hline
$\gamma_{\oE}$ & $\oldz^{-1}(\bold{x}(\ote))$ \,\, (see {\it Step 2})&  \\
$\gamma_0^{(h)}$ & $\oldz^{-1}(\bold{x}_0^{(h)})$ \,\, (see \eqref{eq:I1}) &   \\
$\bar{\mathcal C}_\varepsilon$ & $\{\bold x\in\mathbb R^k:\,\,\mathrm{dist}(\bold x,\bar{\mathcal C})\leq\varepsilon\}$ \,\, (see {\it Step 3.3}) & \\
$\oTheta_n(\delta)$ & $\sup\{t:\,{\color{black}U_n^{(i)}}(t)\le\lfloor \oldz_i(\gamma_{\oE}+\delta)\nac{i}{n}\rfloor , \; \forall\, 1\le i \le k\}$
& \\
$\mathcal{F}_{\bold{v}}^{(j)} $ &  $\{\bold{w}\in \mathbb{N}^k:\,\,w_j=v_j,\,\, w_h\le v_h, \forall h\neq j\} $ & $ {\bold v} \in \mathbb N^k$ \\
$\mathcal{F}_{\bold v} $  & $ \mathcal{F}_{\bold v}:=\cup_{j=1}^{k}\mathcal{F}^{(j)}_{\bold v}$ \,\, (see {\it Step 4}) &  \\
 \hline
\multicolumn{3}{|c|}{Functions/curves}\\
\hline
 Symbol & Mathematical definition & Remarks\\
 \hline
  $ \oldz(\gamma)$ & parametrization of the curve with trace $ \mathcal C$ \,\, (see above \eqref{eq:tau}) & $\gamma\in [0,1]$ \\
  $\otaun(\gamma)$ & $ \sum_{i=1}^{k}\lfloor \oldz_i(\gamma) \nac{i}{n} \rfloor $ \,\, (see \eqref{eq:tau}) & $\gamma\in [0,1]$\\
  $\otaumn(t_n^{(s)})$ & $\inf\{ \gamma \in [0,1] : \otaun(\gamma)=t_n^{(s)}\}$ \,\, (see \eqref{eq:taumeno1}) &  \\
  $\ozeta^{(i)}_n(t)$ & (see \eqref{eq:zetaext})
  &
  non-decreasing \\
  $\bar{\mathcal C}$ & $\oldz([\gamma_0^{({\color{black}\overline{d}})},\gamma_{\oE}-\delta])$ \,\, (see {\it Step 3.3})& \\
\hline
\multicolumn{3}{|c|}{Probabilities}\\
\hline
 Symbol & Mathematical definition & Remarks\\
 \hline
  $T_{1,n}({\color{black}\overline{d}})$  & (see \eqref{eq:hug0n}) & \\
  $T_{2,n}({\color{black}\overline{d}})$  & (see \eqref{eq:uso5bisn}) &  \\
  \hline
 \end{tabular}
 \end{center}
\end{table}

In this subsection we show that \eqref{eq:limsub} holds with
\begin{equation}\label{eq:xstar}
x_*:=\sum_{i=1}^{k}x_i(\ote)(\nu_{1i}(\mu_{1i})^r)^{1/(r-1)}>0,
\end{equation}
where $\bold{x}(y)=(x_1(y),\ldots,x_k(y))$ is the solution of Cauchy problem \eqref{eq:cauchy} and $y_{\mathcal E}$
is defined by \eqref{eq:tE}. Our proof uses
Lemmas \ref{le:June15}, \ref{le:aspi} and \ref{le:initialpoint} below, which are proved in Subsection \ref{subsec:lemmas}.

\begin{Lemma}\label{le:June15}
Assume \eqref{eq:a1a3}, \eqref{eq:hyp2bis}, \eqref{eq:hyp4}, \eqref{eq:gratiog} and \eqref{eq:trivial}.
Let $\mathcal W$ be a compact set such that $\mathcal{W}\subset [\varepsilon,\infty)^k$, for some $\varepsilon>0$. Then
\begin{equation*}
\lim_{n\to\infty}\sup_{\bold{x}\in\mathcal{W}}\Big|\frac{{\color{black}\oldB_n^{(i)}}(\lfloor \bold{x}\nan\rfloor,\obq{{i}})}{\nac{i}{n}}-\oldb_i(\bold x)\Big|=0,
\quad i=1,\ldots,k.
\end{equation*}
\end{Lemma}

\begin{Lemma}\label{le:aspi}
Assume \eqref{eq:a1a3}, \eqref{eq:hyp2bis} and \eqref{eq:hyp4}. Then, for $i=1,\ldots,k$,
\begin{align*}
\oldpi(\lfloor\bold{x}\nan\rfloor,\obq{i})&=
\left(1+O\left(\sum_{j:\,\,x_j>0}(\lfloor x_j \nac{j}{n}\rfloor \oq{{ij}}+ (\lfloor x_j \nac{j}{n}\rfloor)^{-1})
\right)\right)
\left(\sum_{j=1}^{k}\lfloor x_j \nac{j}{n}\rfloor \oq{{ij}}\right)^{r}/r!.
\end{align*}
\end{Lemma}

\begin{Lemma}\label{le:initialpoint}
Assume \eqref{eq:a1a3}, \eqref{eq:hyp2bis}, \eqref{eq:hyp4}, \eqref{eq:gratiog} and \eqref{eq:trivial}. Then:\\
\noindent$(i)$ $\oldb_1(\bold{x}_0^{(h)})>\alpha_1/2$, for any $h=0,\ldots,{\color{black}\overline{d}}$.\\
\noindent$(ii)$ For a fixed $h\in\{1,\ldots,{\color{black}\overline{d}}\}$, if $i\in\bigcup_{s=1}^{h}\mathcal{K}_{s}$, then $\oldb_i(\bold{x}_0^{(h)})\geq\frac{r^{-1}(1-r^{-1})^{r-1}}{2}(\chi\beta_{{\color{black}d_i}-1})^r>0$.
\\
\noindent$(iii)$ For a fixed $h\in\{0,\ldots,{\color{black}\overline{d}}\}$, if $i\in\mathcal{K}_{h+1}$, then $\oldb_i(\bold{x}_0^{(h)})\geq r^{-1}(1-r^{-1})^{r-1}(\chi\beta_{h})^r>0$.\\
\noindent$(iv)$ For a fixed $h\in\{0,\ldots,{\color{black}\overline{d}}-1\}$, if $i\notin\bigcup_{s=0}^{h+1}\mathcal{K}_{s}$, then $\oldb_i(\bold{x}_0^{(h)})\geq 0$.\\
\noindent$(v)$ $\oldb_i(\bold{x}_0)>0$, for any $i=1,\ldots,k$.
\end{Lemma}

\noindent{\it Proof\,\,of\,\,Theorem\,\,\ref{teo-blocksub}.} We divide the proof in four main steps.\\
\noindent{\it Step\,\,1:\,\,chosing\,\,a\,\,suitable\,\,\lq\lq strategy".} Let $\zeta:[0,1]\to\mathbb{R}^k$ be a parametrization of the curve with trace $\mathcal C$
defined in $(\bold{Sub})'$  (recall that by Theorem \ref{prop:sub} $({\mathcal Sub})$ is equivalent to $(\bold{Sub})'$) and set
\begin{equation}\label{eq:tau}
\otaun(\gamma):=\sum_{i=1}^{k}\lfloor \oldz_i(\gamma) \nac{i}{n}\rfloor,\quad \gamma\in [0,1].
\end{equation}
Note that $\otaun([0,1])$ is a finite subset of $\mathbb{N}\cup \{0\}$, i.e.,
$\otaun([0,1])=\{t_n^{(0)},t_n^{(1)},\ldots,t_n^{(m+1)}\}$. Without loss of generality we assume $t_n^{(0)}:=0<t_n^{(1)}<\ldots<t_n^{(m)}<t_n^{(m+1)}:=\otaun(1)$.
We set
\begin{equation}\label{eq:taumeno1}
\otaumn(t_n^{(s)}):=\inf\{\gamma\in [0,1]:\,\,\otaun(\gamma)=t_n^{(s)}\},\quad\text{$s=0,\ldots,m+1$}
\end{equation}
and
\begin{equation}\label{eq:zetasullimag}
\ozeta_n^{(i)}(t_n^{(s)}):=\lfloor \oldz_i(\otaumn(t_n^{(s)}))\nac{i}{n}\rfloor,\quad\text{$i=1,\ldots,k$, $s=0,\ldots,m+1$.}
\end{equation}
Now we extend the definition of  $\ozeta_n^{(i)}(\cdot)$, to  set $(0,\otaun(1))\cap\mathbb N$ , by
``interpolating'' the values in $\otaun([0,1])$ as follows. First, for $s=0,\ldots,m$, we define vector $\bold{e}(t_n^{(s)},t_{n}^{(s+1)})\in\{0,1\}^k$
with components
\[
\ozeta_n^{(i)}(t_n^{(s+1)})-\ozeta_n^{(i)}(t_n^{(s)})\in\{0,1\},\quad\text{$i=1,\ldots,k$.}
\]
Then we define recursively  indexes
\begin{align}
j_1&:=\min\{i\in\{1,\ldots,k\}:\,\,(\bold{e}(t_n^{(s)},t_{n}^{(s+1)}))_i=1\}\nonumber\\
j_{\ell+1}&:=\min\{i\in\{1,\ldots,k\}\setminus\{j_1,\ldots,j_\ell\}:\,\,(\bold{e}(t_n^{(s)},t_{n}^{(s+1)}))_i=1\},\quad \ell=1,\ldots,t_n^{(s+1)}-t_{n}^{(s)}-1.\nonumber
\end{align}
Finally, for any $t\in (0,\otaun(1))\cap\mathbb N)\setminus\otaun([0,1])$, there exist $t_n^{(s)},t_n^{(s+1)}\in\otaun([0,1])$, $s=0,\ldots,m$, such that
$t_n^{(s)}<t<t_n^{(s+1)}$ and we set
\begin{equation}\label{eq:zetaext}
\ozeta_n^{(i)}(t):=\lfloor \oldz_i(\otaumn(t_n^{(s)}))\nac{i}{n}\rfloor+\bold{1}(i\in\mathcal{J}_n^{(t-t_n^{(s)})}),
\end{equation}
where $\mathcal{J}_n^{(t-t_n^{(s)})}:=\{j_1,\ldots,j_{t-t_n^{(s)}}\}$.
Note that, by construction,
\[
\sum_{i=1}^{k}\ozeta_n^{(i)}(t)=t,\quad\text{$\forall$ $t\in\{0,\ldots,\otaun(1)\}$.}
\]
We define random variables
\begin{equation}
T'_n:=T_n''-1\quad\text{where}\quad T''_n:=\min\{1\leq t\leq\otaun(1):\,\,{\color{black}A_n^{(i)}}(t-1)<\ozeta_n^{(i)}(t)
\,\,\text{for some $1\le i\le k$}\}.\label{eq:tnprimo}
\end{equation}
For $0\leq t\leq T_n'$, we consider the \lq\lq policy" defined by
\begin{equation}\label{eq:policydet}
{\color{black}C_n^{(i)}}(0):=0,\quad {\color{black}C_n^{(i)}}(t):=\ozeta_n^{(i)}(t)-\ozeta_n^{(i)}(t-1),\quad\text{$i=1,\ldots,k$, $1\leq t\leq T_n'$.}
\end{equation}
Note that ${\color{black}C_n^{(i)}}(t)$ is the indicator function of  event
\text{\{$G_i$ is selected at time $t$\}}. Indeed, by construction, for each time step $t\leq T'_n$, there exists only one index $i\in\{1,\ldots,k\}$
such that ${\color{black}C_n^{(i)}}(t)=1$ and ${\color{black}C_n^{({j})}}(t)=0$ for each $j\neq i$.
For $T'_n <t\leq T_n $, we use an arbitrary \lq\lq (feasible)  strategy",
for concreteness we chose the \lq\lq strategy" defined by \eqref{eq:strategy}.
In words, the chosen \lq\lq strategy" is deterministic and equal to \eqref{eq:policydet} as long as it is possible.
Indeed, $T''_n$ is the first time at which the deterministic \lq\lq strategy" \eqref{eq:policydet} can not be employed
because of the lack of \lq\lq usable" active nodes in the  selected community at time $T''_n$.
Note that defining   $\bm{\ozeta}_n(t):=(\ozeta_n^{(1)}(t),\ldots,\ozeta_n^{(k)}(t))$, we have
\begin{equation}\label{eq:policyadd1}
\bold{U}_n(t)=\bm{\ozeta}_n(t),\,\,\forall t\leq T_n'\quad\text{and}\quad
\bold{U}_n(T_n'')\neq\bm{\ozeta}_n(T_n'').
\end{equation}
From which, one can easily check the consistency of  our  construction,  given that  $T_n=A^*_n\geq T'_n$.
\\
\noindent{\it Step\,\,2:\,\,outline\,\,of\,\,the\,\,proof.} Let $\gamma_{\oE}\in (0,1)$ be such that $\oldz(\gamma_{\oE})=\bold{x}(\ote)$
and note that
\[
\otaun(\gamma_{\oE})=\sum_{i=1}^{k}\lfloor \oldz_i(\gamma_\oE) \nac{i}{n}\rfloor=\sum_{i=1}^{k}\lfloor x_i(\ote) \nac{i}{n}\rfloor.
\]
Throughout this proof, for an arbitrarily small $\delta>0$, we shall consider quantities $\otaun(\gamma_{\oE}\pm\delta)$.
It is easily seen that
\[
\lim_{\delta\to 0}\limsup_{n \to \infty}\otaun(\gamma_{\oE}+\delta)/\otaun(\gamma_{\oE})=1
\]
and
\[
\lim_{\delta\to 0}\liminf_{n\to \infty}\otaun(\gamma_{\oE}-\delta)/\otaun(\gamma_{\oE})=1.
\]
Therefore, for any $\varepsilon\in (0,1)$, there exist  $\delta_\varepsilon>0$ and $n_{\varepsilon}\in\mathbb N$ such that for any $\delta<\delta_\varepsilon$ and any $n\geq n_{\varepsilon}$,
it holds $\otaun(\gamma_{\oE}+\delta)<(1+\varepsilon)\otaun(\gamma_{\oE})$ and
$\otaun(\gamma_{\oE}-\delta)>(1-\varepsilon)\otaun(\gamma_{\oE})$. So, for an arbitrarily fixed $\varepsilon\in (0,1)$, any $n\ge n_\varepsilon$ and all $\delta<\delta_\varepsilon$,
\begin{align}
\{|A_n^*/\otaun(\gamma_{\oE})-1|>\varepsilon\}&=\{A_n^*>(1+\varepsilon)\otaun(\gamma_{\oE}), A_n^*\geq\otaun(\gamma_{\oE})\}\cup\{A_n^*<(1-\varepsilon)\otaun(\gamma_{\oE}),
A_n^*<\otaun(\gamma_{\oE})\}\nonumber\\
&\subseteq\{A_n^*>\otaun(\gamma_{\oE}+\delta)\}\cup\{T'_n<\otaun(\gamma_{\oE}-\delta)\}.\nonumber
\end{align}
By the definition of $\nac{i}{n}$, \eqref{eq:a1a3} and the second relation in \eqref{eq:gratiog}
we have $\otaun(\gamma_{\oE})/\nac{1}{n}\to x_*$. So the claim follows if we prove that, for $\delta$ small enough,
\begin{equation}\label{eq:Iconvinprob}
P(T'_n<\otaun(\gamma_{\oE}-\delta))=O(\mathrm{e}^{-c(\varepsilon)g_n^{(1)}})
\end{equation}
and
\begin{equation}\label{eq:IIconvinprob}
P(A_n^*>\otaun(\gamma_{\oE}+\delta))=O(\mathrm{e}^{-c(\varepsilon)g_n^{(1)}}).
\end{equation}
\noindent {\it Step\,\,3:\,\,proof\,\,of\,\,\eqref{eq:Iconvinprob}.} We divide the proof of \eqref{eq:Iconvinprob} in three further steps.
In the Step 3.1 we prove inequality
\begin{equation}\label{eq:T1T2}
P(T_n'<\otaun(\gamma_{\oE}-\delta))\leq 1-T_{1,n}({\color{black}\overline{d}})T_{2,n}({\color{black}\overline{d}}),\quad\text{for all $n$ large enough}
\end{equation}
 where
\begin{align}
T_{1,n}({\color{black}\overline{d}}):=\prod_{h=1}^{{\color{black}\overline{d}}}\prod_{s={\color{black}\hat{t}_n^{(h-1)}}}^{{\color{black}\hat{t}_n^{(h)}}-1}\Biggl(1-
\sum_{j\in\mathcal{K}_h}P(\mathrm{Bin}(\on{j}-\oan{j},\oldpi(\bm{\ozeta}_n(s),\obq{j}))<\ozeta_n^{(j)}(s+1)-\oan{j})\Biggr),\label{eq:hug0n}
\end{align}
\begin{align}
T_{2,n}({\color{black}\overline{d}}):=\prod_{s={\color{black}\hat{t}_n^{({\color{black}\overline{d}})}}}^{\otaun(\gamma_{\oE}-\delta)-1}
\Biggl(1-\sum_{j=1}^{k}P(\mathrm{Bin}(\on{j}-\oan{j},\oldpi(\bm{\ozeta}_n(s),\obq{j}))<\ozeta_n^{(j)}(s+1)-\oan{j})
\Biggr)\label{eq:uso5bisn}
\end{align}
and
\begin{equation}\label{eq:tmin}
{\color{black}\hat{t}_n^{(h)}}:=\sum_{i=1}^{k}\lfloor\sum_{\ell=0}^{h}\beta_\ell\bold{1}(i\in\mathcal{K}_\ell)\nac{i}{n}\rfloor,\quad h\in\{0,\ldots,{\color{black}\overline{d}}\}.
\end{equation}
In  Step 3.2 we prove
\begin{equation}\label{eq:T1zero}
T_{1,n}({\color{black}\overline{d}})\geq 1-O(\mathrm{e}^{-c(\varepsilon)g_n^{(1)}})
\end{equation}
and in  Step 3.3 we prove
\begin{equation}\label{eq:T2zero}
T_{2,n}({\color{black}\overline{d}})\geq 1-O(\mathrm{e}^{-c(\varepsilon)g_n^{(1)}}).
\end{equation}
\noindent {\it Step\,\,3.1:\,\,proof\,\,of\,\,\eqref{eq:T1T2}.}
\textcolor{black}{Hereon, we set $\bold{A}_n(t):=({\color{black}A_n^{(1)}}(t),\ldots,{\color{black}A_n^{(k)}}(t))$.} Since, for $t\in\{0,\ldots,\otaun(1)-1\}$,
\begin{align}
\{T'_n>t\}&=\{\bold{A}_n(s)\geq\bm{\ozeta}_n(s+1)\,\,
\forall 1\le s\le t\}=\{\bold{U}_n(s)=\bm{\ozeta}_n(s)\,\,\forall 1\le s\le t+1\},\label{eq:defTgretat}
\end{align}
by \eqref{eq:policyadd1}
we have
\begin{align}
P(T'_n>t\mid T'_n >t-1)&=P(\bold{A}_n(t)\geq\bm{\ozeta}_n(t+1)\,\mid\,
\bold{U}_n(s)=\bm{\ozeta}_n(s)\,\,\forall 1\le s\le t)\nonumber\\
&=P(\bold{A}_n(t)\geq\bm{\ozeta}_n(t+1)\,\mid\,\bold{U}_n(t)=\bm{\ozeta}_n(t)),\label{eq:july5}
\end{align}
where the last equality follows from the fact that, given $\{\bold{U}_n(t)=\bm{\ozeta}_n(t)\}$,
random vector $\bold{A}_n(t)$ is independent of $\{\bold{U}_n(s)=\bm{\ozeta}_n(s)\,\forall 1\leq s<t\}$. So
\begin{align}
P(T'_n>t)&=\prod_{s=0}^{t}P(T'_n>s\mid T'_n>s-1)
=\prod_{s=0}^{t}P(\bold{A}_n(s)\geq\bm{\ozeta}_n(s+1)\,\mid\,\bold{U}_n(s)=\bm{\ozeta}_n(s)).\nonumber
\end{align}
Consequently,
\begin{align}
&P(T_n'<\otaun(\gamma_{\oE}-\delta))=1-P(T_n'\geq\otaun(\gamma_{\oE}-\delta))\nonumber\\
&\leq 1-\prod_{s=0}^{\otaun(\gamma_{\oE}-\delta)}P(\bold{A}_n(s)\geq\bm{\ozeta}_n(s+1)\,\mid\,
\bold{U}_n(s)=\bm{\ozeta}_n(s))\nonumber\\
&\leq 1-\prod_{s=0}^{\otaun(\gamma_{\oE}-\delta)}\Biggl(1-\sum_{j=1}^{k}P({\color{black}S_n^{(j)}}(s)+\oan{j}-\ozeta_n^{(j)}(s+1)< 0\,\,
\mid
\,\bold{U}_n(s)=\bm{\ozeta}_n(s))\Biggr)\nonumber\\
&\leq 1-\prod_{h=0}^{{\color{black}\overline{d}}+1}\prod_{s={\color{black}\hat{t}_n^{(h-1)}}}^{{\color{black}\hat{t}_n^{(h)}}-1}\Biggl(1-\sum_{j=1}^{k}P({\color{black}S_n^{(j)}}(s)+\oan{j}-\ozeta_n^{(j)}(s+1)< 0\,
\mid\,\bold{U}_n(s)=\bm{\ozeta}_n(s))\Biggr),\label{eq:disprod}
\end{align}
where ${\color{black}\hat{t}_n^{(-1)}}:=0$, ${\color{black}\hat{t}_n^{({\color{black}\overline{d}}+1)}}:=\otaun(\gamma_{\oE}-\delta)$.
Since $\zeta$ is a parametrization of $\mathcal C$ and $\bold{x}_0^{(h)}=\sum_{\ell=0}^{h}\beta_{\ell}\bold{1}_{\mathcal{K}_{\ell}}\in\mathcal S\subset\mathcal C$
we have
\begin{equation}\label{eq:I1}
\oldz(\gamma_0^{(h)})=\bold{x}_0^{(h)}\quad\text{for some $\gamma_0^{(h)}\in [0,1]$}
\end{equation}
and so
\begin{align}
\otaun(\gamma_0^{(h)})=\sum_{i=1}^{k}\lfloor \oldz_i(\gamma_0^{(h)})\nac{i}{n}\rfloor=\sum_{i=1}^{k}\lfloor\sum_{\ell=0}^{h}\beta_{\ell}\bold{1}(i\in\mathcal{K}_{\ell})\nac{i}{n}\rfloor
={\color{black}\hat{t}_n^{(h)}}.\label{eq:I2}
\end{align}
For $h\in\{1,\ldots,{\color{black}\overline{d}}\}$, if $i\notin\bigcup_{\ell=0}^{h}\mathcal{K}_\ell$ then $\ozeta_n^{(i)}({\color{black}\hat{t}_n^{(h)}})=0$
and so (by monotonicity of $\ozeta_n^{(i)}(\cdot)$) $\ozeta_n^{(i)}(s)=0$ for any $s\leq {\color{black}\hat{t}_n^{(h)}}$. Consequently, for $s\leq {\color{black}\hat{t}_n^{(h)}}-1$,
\begin{align}
&\sum_{j=1}^{k}P({\color{black}S_n^{(j)}}(s)+\oan{j}-\ozeta_n^{(j)}(s+1)<0\,
\mid\,\bold{U}_n(s)=\bm{\ozeta}_n(s))\nonumber\\
&=\sum_{j\in\bigcup_{\ell=0}^{h}\mathcal{K}_\ell}P({\color{black}S_n^{(j)}}(s)+\oan{j}-\ozeta_n^{(j)}(s+1)<0\,
\mid\,\bold{U}_n(s)=\bm{\ozeta}_n(s)).\label{eq:k1}
\end{align}
Now we check that, for ${\color{black}\hat{t}_n^{(h-1)}}\leq s<{\color{black}\hat{t}_n^{(h)}}$, we have
\begin{align}
&\sum_{j\in\bigcup_{\ell=0}^{h}\mathcal{K}_\ell}P({\color{black}S_n^{(j)}}(s)+\oan{j}-\ozeta_n^{(j)}(s+1)<0\,
\mid\,\bold{U}_n(s)=\bm{\ozeta}_n(s))\nonumber\\
&=\sum_{j\in\mathcal{K}_h}P({\color{black}S_n^{(j)}}(s)+\oan{j}-\ozeta_n^{(j)}(s+1)<0\,
\mid\,\bold{U}_n(s)=\bm{\ozeta}_n(s)).\label{eq:sumsimple}
\end{align}
This in turn follows if we check that, for any $j\in\bigcup_{\ell=0}^{h-1}\mathcal{K}_\ell$, we have
\begin{equation}\label{eq:equalzeta}
\ozeta_n^{(j)}({\color{black}\hat{t}_n^{(h-1)}})=\ozeta_n^{(j)}({\color{black}\hat{t}_n^{(h)}}).
\end{equation}
Indeed this implies that quantity $\ozeta_n^{(j)}(s)$ is constant on interval $[{\color{black}\hat{t}_n^{(h-1)}},{\color{black}\hat{t}_n^{(h)}}]$, and therefore (reasoning as for the derivation of relation \eqref{eq:july5}),
for any $j\in\bigcup_{\ell=0}^{h-1}\mathcal{K}_\ell$,
\begin{align}
&P({\color{black}S_n^{(j)}}(s)+\oan{j}-\ozeta_n^{(j)}(s+1)<0\,
\mid\,\bold{U}_n(s)=\bm{\ozeta}_n(s))\nonumber\\
&=P({\color{black}A_n^{(j)}}(s)<\ozeta_n^{(j)}(s+1)\,\mid\,
\bold{U}_n(u)=\bm{\ozeta}_n(u)\,\,\forall 1\le u\le s)\nonumber\\
&=P({\color{black}A_n^{(j)}}(s)-\ozeta_n^{(j)}(s)<0\,
\mid\,\bold{U}_n(u)=\bm{\ozeta}_n(u)\,\,\forall 1\le u\le s)\nonumber\\
&=P({\color{black}A_n^{(j)}}(s)-\ozeta_n^{(j)}(s)<0\,
\mid\,\bold{A}_n(u)\geq\bm{\ozeta}_n(u+1)\,\,\forall 1\le u\le s-1)\nonumber\\
&\leq\frac{P({\color{black}A_n^{(j)}}(s)-\ozeta_n^{(j)}(s)<0,\,{\color{black}A_n^{(j)}}(s-1)\ge\ozeta_n^{(j)}(s))}{P(\bold{A}_n(u)\geq\bm{\ozeta}_n(u+1)\,\,\forall
1\le u\le s-1)}=0,\nonumber
\end{align}
where the latter equality follows by monotonicity of ${\color{black}A_n^{(j)}}(\cdot)$. Now we check \eqref{eq:equalzeta}. Let $\theta\bold{x}_0^{(h)}+(1-\theta)\bold{x}_0^{(h-1)}$, $\theta\in [0,1]$,
be a generic point of $\mathcal{S}_h$. We have
\begin{equation}\label{eq:combconv}
\theta\bold{x}_0^{(h)}+(1-\theta)\bold{x}_0^{(h-1)}=\theta\beta_h\bold{1}_{\mathcal{K}_h}+\sum_{\ell=0}^{h-1}\beta_\ell\bold{1}_{\mathcal{K}_\ell},
\end{equation}
and therefore, for any $j\in\bigcup_{\ell=0}^{h-1}\mathcal{K}_\ell$, i.e., for any $j$ such that ${\color{black}d_j}\leq h-1$ (where ${\color{black}d_j}$ is defined in Subsection \ref{sect:addnot}),
\begin{align}
\theta(\bold{x}_0^{(h)})_j+(1-\theta)(\bold{x}_0^{(h-1)})_j&=\theta\beta_h\bold{1}(j\in\mathcal{K}_h)+\sum_{\ell=0}^{h-1}\beta_\ell\bold{1}(j\in\mathcal{K}_\ell)=\beta_{{\color{black}d_j}}.\nonumber
\end{align}
Consequently,
\begin{align}
\ozeta_n^{(j)}({\color{black}\hat{t}_n^{(h)}})&=\lfloor \oldz_j(\otaumn({\color{black}\hat{t}_n^{(h)}}))\nac{j}{n}\rfloor=\lfloor \oldz_j(\gamma_0^{(h)})\nac{j}{n}\rfloor
=\lfloor\beta_{{\color{black}d_j}}\nac{j}{n}\rfloor=\ozeta_n^{(j)}({\color{black}\hat{t}_n^{(h-1)}}),\nonumber
\end{align}
and \eqref{eq:equalzeta} is checked. Then, exploiting \eqref{eq:k1} and \eqref{eq:sumsimple}, we rewrite \eqref{eq:disprod} as
\begin{align}
&P(T_n'<\otaun(\gamma_{\oE}-\delta))\nonumber\\
&\qquad\qquad\qquad
\leq 1-\prod_{h=0}^{{\color{black}\overline{d}}+1}\prod_{s={\color{black}\hat{t}_n^{(h-1)}}}^{{\color{black}\hat{t}_n^{(h)}}-1}\Biggl(1-\sum_{j\in\mathcal{K}_h}P({\color{black}S_n^{(j)}}(s)+\oan{j}-\ozeta_n^{(j)}(s+1)<0\,
\mid\,\bold{U}_n(s)=\bm{\ozeta}_n(s))\Biggr).\nonumber
\end{align}
Inequality \eqref{eq:T1T2} follows noticing that, for all $n$ large enough,
\begin{align}
&\prod_{h=0}^{{\color{black}\overline{d}}}\prod_{s={\color{black}\hat{t}_n^{(h-1)}}}^{{\color{black}\hat{t}_n^{(h)}}-1}\Biggl(1-\sum_{j\in\mathcal{K}_h}P({\color{black}S_n^{(j)}}(s)+\oan{j}-\ozeta_n^{(j)}(s+1)< 0\,
\mid\,\bold{U}_n(s)=\bm{\ozeta}_n(s))\Biggr)\nonumber\\
&=\prod_{h=0}^{{\color{black}\overline{d}}}\prod_{s={\color{black}\hat{t}_n^{(h-1)}}}^{{\color{black}\hat{t}_n^{(h)}}-1}\Biggl(1-
\sum_{j\in\mathcal{K}_h}P(\mathrm{Bin}(\on{j}-\oan{j},\oldpi(\bm{\ozeta}_n(s),\obq{j}))<\ozeta_n^{(j)}(s+1)-\oan{j})\Biggr)\label{eq:uso5}\\
&=\prod_{h=1}^{{\color{black}\overline{d}}}\prod_{s={\color{black}\hat{t}_n^{(h-1)}}}^{{\color{black}\hat{t}_n^{(h)}}-1}\Biggl(1-
\sum_{j\in\mathcal{K}_h}P(\mathrm{Bin}(\on{j}-\oan{j},\oldpi(\bm{\ozeta}_n(s),\obq{j}))<\ozeta_n^{(j)}(s+1)-\oan{j})\Biggr)\nonumber\\
&=T_{1,n}({\color{black}\overline{d}})\label{eq:hug0}
\end{align}
and
\begin{align}
&\prod_{s={\color{black}\hat{t}_n^{({\color{black}\overline{d}})}}}^{\otaun(\gamma_{\oE}-\delta)-1}\Biggl(1-\sum_{j=1}^kP({\color{black}S_n^{(j)}}(s)+\oan{j}-\ozeta_n^{(j)}(s+1)< 0\,
\mid\,\bold{U}_n(s)=\bm{\ozeta}_n(s))\Biggr)\nonumber\\
&=\prod_{s={\color{black}\hat{t}_n^{({\color{black}\overline{d}})}}}^{\otaun(\gamma_{\oE}-\delta)-1}
\Biggl(1-\sum_{j=1}^kP(\mathrm{Bin}(\on{j}-\oan{j},\oldpi(\bm{\ozeta}_n(s),\obq{j}))<\ozeta_n^{(j)}(s+1)-\oan{j})\Biggr).\label{eq:uso5bis}
\end{align}
Here in \eqref{eq:uso5} and \eqref{eq:uso5bis} we used \eqref{eq:bin}, and in \eqref{eq:hug0} we used that, for all $n$ large enough and all $0\leq s\leq {\color{black}\hat{t}_n^{(0)}}-1$,
\begin{align}
&P\Biggl(\mathrm{Bin}(\on{1}-\oan{1},\oldpi(\bm{\ozeta}_n(s),\obq{1})<\ozeta_n^{(1)}(s+1)-\oan{1}\Biggr)\nonumber\\
&\quad\quad\quad
\leq P\Biggl(\mathrm{Bin}(\on{1}-\oan{1},\oldpi(\bm{\ozeta}_n(s),\obq{1})<\lfloor(\alpha_1/2)\nac{1}{n}\rfloor-\oan{1}\Biggr)=0.\nonumber
\end{align}
\noindent {\it Step\,\,3.2:\,\,proof\,\,of\,\,\eqref{eq:T1zero}.} Let $h\in\{1,\ldots,{\color{black}\overline{d}}\}$ and $j\in\mathcal{K}_h$ be fixed. By Lemma
\ref{le:bas2}, Lemma \ref{le:initialpoint}$(iii)$, \eqref{eq:I1} and \eqref{eq:I2}, for any $\varepsilon>0$ there exists $n_\varepsilon(j,h)$ such that for all $n>n_\varepsilon(j,h)$
\begin{align}
{\color{black}\oldB_n^{(j)}}(\bm{\ozeta}_n({\color{black}\hat{t}_n^{(h-1)}}),\obq{j})&={\color{black}\oldB_n^{(j)}}(\lfloor \oldz(\gamma_0^{(h-1)})\nan\rfloor,\obq{j})\nonumber\\
&={\color{black}\oldB_n^{(j)}}(\lfloor\bold{x}_0^{(h-1)}\nan\rfloor,\obq{j})\nonumber\\
&>(\oldb_j(\bold{x}_0^{(h-1)})-\varepsilon)\nac{j}{n}\nonumber\\
&>(r^{-1}(1-r^{-1})^{r-1}(\chi\beta_{h-1})^r-\varepsilon)\nac{j}{n}.\nonumber
\end{align}
We note that, for any ${\color{black}\hat{t}_n^{(h-1)}}\leq s<{\color{black}\hat{t}_n^{(h)}}$, $\ozeta_n^{(j)}(s)\geq\ozeta_n^{(j)}({\color{black}\hat{t}_n^{(h-1)}})$ and so,
for any ${\color{black}\hat{t}_n^{(h-1)}}\leq s<{\color{black}\hat{t}_n^{(h)}}$,
\begin{align}
{\color{black}\oldB_n^{(j)}}(\bm{\ozeta}_n(s),\obq{j})&\geq {\color{black}\oldB_n^{(j)}}(\bm{\ozeta}_n({\color{black}\hat{t}_n^{(h-1)}}),\obq{j})-(\ozeta_n^{(j)}(s)-\ozeta_n^{(j)}({\color{black}\hat{t}_n^{(h-1)}}))\nonumber\\
&\geq {\color{black}\oldB_n^{(j)}}(\bm{\ozeta}_n({\color{black}\hat{t}_n^{(h-1)}}),\obq{j})-(\ozeta_n^{(j)}({\color{black}\hat{t}_n^{(h)}})-\ozeta_n^{(j)}({\color{black}\hat{t}_n^{(h-1)}}))\nonumber\\
&>(r^{-1}(1-r^{-1})^{r-1}(\chi\beta_{h-1})^r-\varepsilon)\nac{j}{n}-(\ozeta_n^{(j)}({\color{black}\hat{t}_n^{(h)}})-\ozeta_n^{(j)}({\color{black}\hat{t}_n^{(h-1)}})).\label{eq:1}
\end{align}
By \eqref{eq:combconv}
\begin{align}
\ozeta_n^{(j)}({\color{black}\hat{t}_n^{(h)}})-\ozeta_n^{(j)}({\color{black}\hat{t}_n^{(h-1)}})
&=\lfloor (\bold{x}_0^{(h)})_j\nac{j}{n}\rfloor-\lfloor (\bold{x}_0^{(h-1)})_j\nac{j}{n}\rfloor\nonumber\\
&=\lfloor\beta_h \nac{j}{n}\rfloor.\label{eq:2}
\end{align}
Combining \eqref{eq:1}, \eqref{eq:2} and \eqref{eq:beta}, for any fixed $h\in\{1,\ldots,{\color{black}\overline{d}}\}$, $j\in\mathcal{K}_h$ and ${\color{black}\hat{t}_n^{(h-1)}}\leq s< {\color{black}\hat{t}_n^{(h)}}$,
for any $\varepsilon>0$ there exists $n_\varepsilon(j,h)$ such that for all $n>n_\varepsilon(j,h)$
\[
{\color{black}\oldB_n^{(j)}}(\bm{\ozeta}_n(s),\obq{j})>(2^{-1}r^{-1}(1-r^{-1})^{r-1}(\chi\beta_{h-1})^r-\varepsilon)\nac{j}{n}.
\]
Therefore by the definition of ${\color{black}\oldB_n^{(j)}}$ and the fact that $\ozeta_n^{(j)}(s+1)\leq\ozeta_n^{(j)}(s)+1$, we have
\begin{equation} \nonumber 
(\on{j}-\oan{{j}})\oldpi(\bm{\ozeta}_n(s),\obq{j})
>(2^{-1}r^{-1}(1-r^{-1})^{r-1}(\chi\beta_{h-1})^r-\varepsilon)\nac{j}{n}+(\ozeta_n^{(j)}(s+1)-1-\oan{j}).
\end{equation}
Let
\begin{equation}\label{eq:H}
H(x):=1-x+x\log x,\quad x>0,\quad H(0)=1.
\end{equation}
By a classical concentration inequality for the binomial distribution (see e.g. Eq. $(1.6)$ of Lemma 1.1 in \cite{P}),
monotonicity in $x$ of $x/(x+s)$, $s>0$, monotonicity in $h$ of $\beta_h$ and monotonicity properties of $H$,
for any fixed $h\in\{1,\ldots,{\color{black}\overline{d}}\}$, $j\in\mathcal{K}_h$ and ${\color{black}\hat{t}_n^{(h-1)}}\leq s< {\color{black}\hat{t}_n^{(h)}}$,
for any $\varepsilon>0$ there exists $n_\varepsilon(j,h)$ such that for all $n>n_\varepsilon(j,h)$
we have
\begin{align}
&P\Biggl(\mathrm{Bin}(\on{j}-\oan{j},\oldpi(\bm{\ozeta}_n(s),\obq{j}))\leq\ozeta_n^{(j)}(s+1)-1-\oan{j}\Biggr)\nonumber\\
&\quad\quad
\leq\bold{1}\{\ozeta_n^{(j)}(s+1)-\oan{j}-1\geq 0\}
\exp\left(-(\on{j}-\oan{j})\oldpi(\bm{\ozeta}_n(s),\obq{j})H\left(\frac{\ozeta_n^{(j)}(s+1)-\oan{j}-1}{(\on{j}-\oan{{j}})\oldpi(\bm{\ozeta}_n(s),\obq{j})}\right)\right)\nonumber\\
&\quad\quad
\leq\bold{1}\{\ozeta_n^{(j)}(s+1)-\oan{j}-1\geq 0\}\nonumber\\
&\quad\quad\quad
\times\exp\Biggl(-
(2^{-1}r^{-1}(1-r^{-1})^{r-1}(\chi\beta_{h-1})^r-\varepsilon)\nac{j}{n}\nonumber\\
&\quad\quad\quad
\times H\left(\frac{\ozeta_n^{(j)}(s+1)-\oan{j}-1}{
(2^{-1}r^{-1}(1-r^{-1})^{r-1}(\chi\beta_{h-1})^r-\varepsilon)\nac{j}{n}+\ozeta_n^{(j)}(s+1)-\oan{j}-1
}\right)\Biggr)\nonumber\\
&\quad\quad
\leq\exp\Biggl(-
(2^{-1}r^{-1}(1-r^{-1})^{r-1}(\chi\beta_{h-1})^r-\varepsilon)\nac{j}{n}\nonumber\\
&\quad\quad\quad
\times H\left(\frac{\ozeta_n^{(j)}({\color{black}\hat{t}_n^{(h)}})}{
(2^{-1}r^{-1}(1-r^{-1})^{r-1}(\chi\beta_{h-1})^r-\varepsilon)\nac{j}{n}+\ozeta_n^{(j)}({\color{black}\hat{t}_n^{(h)}})
}\right)\Biggr)\nonumber\\
&\quad\quad
\leq\exp\Biggl(-
(2^{-1}r^{-1}(1-r^{-1})^{r-1}(\chi\beta_{h-1})^r-\varepsilon)\nac{j}{n}H\left(\frac{\beta_h}{
(2^{-1}r^{-1}(1-r^{-1})^{r-1}(\chi\beta_{h-1})^r-\varepsilon)+\beta_h}\right)\Biggr)\nonumber\\
&\quad\quad
=\exp\Biggl(-
(\beta_h-\varepsilon)\nac{j}{n}H\left(\frac{\beta_h}{
2\beta_h-\varepsilon}\right)\Biggr)\nonumber\\
&\quad\quad
\leq\exp\Biggl(-
(\beta_{{\color{black}\overline{d}}}-\varepsilon)\nac{j}{n}H\left(\frac{2}{3}\right)\Biggr).\nonumber
\end{align}
Therefore, for all $n$ large enough,
\begin{equation*}
\max_{1\leq h\leq{\color{black}\overline{d}}}\sup_{{\color{black}\hat{t}_n^{(h-1)}}\leq s<{\color{black}\hat{t}_n^{(h)}}}P\Biggl(\mathrm{Bin}(\on{j}-\oan{j},\oldpi(\bm{\ozeta}_n(s),\obq{j}))\leq\ozeta_n^{(j)}(s+1)-1-\oan{j}\Biggr)
\leq\mathrm{e}^{-\oc{1} \nac{1}{n}},
\end{equation*}
for some positive constant $\oc{1}>0$. By this inequality easily follows that, for all $n$ large enough,
\[
T_{1,n}({\color{black}\overline{d}})\geq (1-k\mathrm{e}^{-\oc{1} \nac{1}{n}})^{{\color{black}\hat{t}_n^{({\color{black}\overline{d}})}}},
\]
which, for all $n$ large enough, yields
\begin{equation} \nonumber 
T_{1,n}({\color{black}\overline{d}})\geq (1-k\mathrm{e}^{-\oc{1} \nac{1}{n}})^{\oc{2} \nac{1}{n}},
\end{equation}
for some positive constant $\oc{2}>0$, and therefore \eqref{eq:T1zero}.\\
\noindent {\it Step\,\,3.3:\,\,proof\,\,of\,\,\eqref{eq:T2zero}.} Recall that $\oldz(\gamma_0^{({\color{black}\overline{d}})})=\bold{x}_0$. Therefore $\bar{\mathcal C}:=\oldz([\gamma_0^{({\color{black}\overline{d}})},\gamma_{\oE}-\delta])$
is a compact set of $\mathbb R^k$
contained in $\mathcal{C}_\oE\cap\overset\circ{\oE}$. For $\varepsilon>0$, let $\bar{\mathcal C}_\varepsilon$ be the $\varepsilon$-thickening of $\bar{\mathcal C}$, i.e.,
\begin{equation*}
\bar{\mathcal C}_\varepsilon:=\{\bold x\in\mathbb R^k:\,\,\mathrm{dist}(\bold x,\bar{\mathcal C})\leq\varepsilon\},
\end{equation*}
where, for $\mathcal B\subset\mathbb R^k$,
\[
\mathrm{dist}(\bold x,\mathcal B):=\inf\{\|\bold x-\bold y\|:\,\,\bold y\in\mathcal B\}
\]
and $\|\cdot\|$ is the Euclidean norm. By the regularity properties of $\bm\rho$, we have that there exists $\varepsilon_0>0$ small enough such that
$\bar{\mathcal C}_{\varepsilon_0}\subset\overset\circ{\oE}$. We are going to show that
there exists a positive integer $\bar n$ such that
$\hat{\bm{\ozeta}}_n(s):=(\ozeta_n^{(1)}(s)/\nac{1}{n},\ldots,\ozeta_n^{(k)}(s)/\nac{k}{n})\in \bar{\mathcal C}_{\varepsilon_0}$ for any $n>\bar n$ and ${\color{black}\hat{t}_n^{({\color{black}\overline{d}})}}\leq s\leq\otaun(\gamma_{\oE}-\delta)-1$.
Indeed, for any $n\in\mathbb N$ and ${\color{black}\hat{t}_n^{({\color{black}\overline{d}})}}\leq s\leq\otaun(\gamma_{\oE}-\delta)-1$, we have
\[
\frac{\lfloor \oldz_j(\otaumn(t_n^{(s)}))\nac{j}{n}\rfloor}
{\nac{j}{n}}\leq\frac{\ozeta_n^{(j)}(s)}{\nac{j}{n}}\leq\frac{\lfloor \oldz_j(\otaumn(t_n^{(s)}))\nac{j}{n}\rfloor+1}{\nac{j}{n}},
\]
i.e.,
\[
\frac{\lfloor \oldz_j(\gamma_s)\nac{j}{n}\rfloor}
{\nac{j}{n}}\leq\frac{\ozeta_n^{(j)}(s)}{\nac{j}{n}}\leq\frac{\lfloor \oldz_j(\gamma_s)\nac{j}{n}\rfloor+1}{\nac{j}{n}}
\]
where $\gamma_s\in [\gamma_0^{({\color{black}\overline{d}})},\gamma_{\oE}-\delta]$ is such that $\otaun(\gamma_s)=t_n^{(s)}$. This relation implies
\[
\Big|\frac{\ozeta_n^{(j)}(s)}{\nac{j}{n}}-\oldz_j(\gamma_s)\Big|\leq (\nac{j}{n})^{-1}
\]
for any $n\in\mathbb N$ and ${\color{black}\hat{t}_n^{({\color{black}\overline{d}})}}\leq s\leq\otaun(\gamma_{\oE}-\delta)-1$; therefore
we can select $\bar n$ such that $\|\hat{\bm{ \ozeta}}_{n}(s)-\oldz(\gamma_s)\|\leq\varepsilon_0$, and so
$\hat{\bm{\ozeta}}_n(s)\in \bar{\mathcal C}_{\varepsilon_0}$,
for any $n>\bar n$ and ${\color{black}\hat{t}_n^{({\color{black}\overline{d}})}}\leq s\leq\otaun(\gamma_{\oE}-\delta)-1$.
Consequently, for any $\varepsilon>0$ and all $n$ sufficiently large,
\begin{align}
T_{2,n}({\color{black}\overline{d}})&\geq\prod_{s={\color{black}\hat{t}_n^{({\color{black}\overline{d}})}}}^{\otaun(\gamma_{\oE}-\delta)-1}
\Biggl(1-\sum_{j=1}^{k}P\Biggl(\mathrm{Bin}(\on{j}-\oan{j},\oldpi(\bm{\ozeta}_n(s),\obq{j})<\ozeta_n^{(j)}(s)-\oan{j}+1\Biggr)
\Biggr)\nonumber\\
&\geq
(1-\mathrm{sup}_n)^{\otaun(\gamma_{\oE}-\delta)}\geq(1-\mathrm{sup}_n)^{\oc{4} \nac{1}{n}},\label{eq:lowerbdsup}
\end{align}
where $\oc{4}>0$ is a positive constant and
\[
\mathrm{sup}_n:=
\sum_{j=1}^{k}\sup_{\bold{z}\in \bar{\mathcal C}_{\varepsilon_0}}P\Biggl(\mathrm{Bin}(\on{j}-\oan{j},\oldpi(\lfloor\bold{z}\nan\rfloor,\obq{j}))/\nac{j}{n}<
z_j-\alpha_j+
\varepsilon\Biggr).
\]
Note that \eqref{eq:T2zero} easily follows by \eqref{eq:lowerbdsup} if we prove that, for some positive constant $\oc{3}>0$,
\begin{align}
\mathrm{sup}_n=O(\mathrm{e}^{-\oc{3}\nac{1}{n}}). \nonumber 
\end{align}
For this it suffices to prove that, for any $j\in\{1,\ldots,k\}$, for all $n$ large enough it holds
\begin{align}
&\sup_{\bold{z}\in\bar{\mathcal C}_{\varepsilon_0}}P\Biggl(\mathrm{Bin}(\on{j}-\oan{j},\oldpi(\lfloor\bold{z}\nan\rfloor,\obq{j}))/\nac{j}{n}<z_j-\alpha_j+
\varepsilon\Biggr)\nonumber\\
&\qquad\qquad\qquad\qquad
=O(\mathrm{e}^{-\oc{3} \nac{1}{n}}). \nonumber 
\end{align}
We shall show this relation for $j=1$ as the other relations can be proved similarly.
Note that $\epsilon:=\min_{\bold z\in\textcolor{black}{\bar{\mathcal C}_{\varepsilon_0}}}\oldb_1(\bold z)>0$. From now on, we choose $\varepsilon<\epsilon$.
By Lemma \ref{le:June15} we have
that for any $0<\varepsilon'<\epsilon-\varepsilon$, for all $n$ large enough
\begin{align}
\frac{(\on{1}-\oan{{1}})\oldpi(\lfloor\bold{z}\nan\rfloor,\obq{1})}{\nac{1}{n}}
&>
r^{-1}(1-r^{-1})^{r-1}\left(\sum_{\ell=1}^{k}z_\ell\chi_{1\ell}\right)^r-\varepsilon'\nonumber\\
&=z_1-\alpha_1+\oldb_1(\bold z)-\varepsilon'\nonumber\\
&>z_1-\alpha_1+\epsilon-\varepsilon'\nonumber\\
&>z_1-\alpha_1+\varepsilon,\quad\text{for all $\bold z\in \bar{\mathcal C}_{\varepsilon_0}$.}\label{eq:relationunif}
\end{align}
Therefore, by a classical concentration inequality for the binomial distribution (see e.g. Eq.
$(1.6)$ of Lemma 1.1. p. 16 in \cite{P}), for all $n$ large enough, we have
\begin{align}
&\sup_{\bold z  \in\bar{\mathcal C}_{\varepsilon_0}}P\Biggl(\mathrm{Bin}(\on{1}-\oan{1},\oldpi(\lfloor\bold{z}\nan\rfloor,\obq{1}))
\leq (z_1-\alpha_1+\varepsilon)\nac{1}{n}\Biggr)\nonumber\\
&\leq\sup_{\bold z\in \bar{\mathcal C}_{\varepsilon_0}}\exp\Biggl(-(\on{1}-\oan{1})\oldpi(\lfloor\bold{z}\nan\rfloor,\obq{1})
H\left(\frac{(z_1-\alpha_1+\varepsilon)\nac{1}{n}}{(\on{1}-\oan{1})\oldpi(\lfloor\bold{z}\nan\rfloor,\obq{1})}\right)\Biggr)\nonumber\\
&\leq\sup_{\bold z\in \bar{\mathcal C}_{\varepsilon_0}}\exp\Biggl(-(z_1-\alpha_1+\varepsilon)\nac{1}{n}H\left(\frac{z_1-\alpha_1+\varepsilon}{z_1-\alpha_1+\epsilon-\varepsilon'}\right)\Biggr)\label{eq:useH01}\\
&=O(\mathrm{e}^{-\oc{3} \nac{1}{n}}),\nonumber
\end{align}
where in \eqref{eq:useH01} we used \eqref{eq:relationunif} and that $H$ decreases on $(0,1)$.
\\
{\it Step\,\,4:\,\,proof\,\,of\,\,\eqref{eq:IIconvinprob}.} Define random variable
\[
\oTheta_n(\delta):=\sup\{t:\,\,{\color{black}U_n^{(i)}}(t)\le\lfloor \oldz_i(\gamma_{\oE}+\delta)\nac{i}{n}\rfloor , \;\; \forall\, 1\le i \le k   \}
\]
and note that $\oTheta_n(\delta)\le\otaun(\gamma_{\oE}+\delta)$ and $U_{n}^{(i_0)}(\oTheta_n(\delta))=\lfloor\zeta_{i_0}(\gamma_{\oE}+\delta)g_n^{(i_0)}\rfloor$ for some $1 \le i_0\le k$,
almost surely. For $\bold v\in\mathbb N^k$, we define set
\[
\mathcal{F}_{\bold v}:=\cup_{j=1}^{k}\mathcal{F}^{(j)}_{\bold v},
\]
where
\[
\mathcal{F}_{\bold{v}}^{(j)}=\{\bold{w}\in \mathbb{N}^k:\,\,w_j=v_j,\,\, w_h\le v_h,
\forall h\neq j\}.
\]
As usual, we denote by $\lfloor \oldz(\gamma_{\oE}+\delta)\nan\rfloor$ the vector with components $\lfloor \oldz_i(\gamma_{\oE}+\delta)\nac{i}{n}\rfloor$, $i=1,\ldots,k$.
By construction,  vector $\bold{U}_n(\oTheta_n(\delta))$ (whose components are  ${\color{black}U_n^{(i)}}(\oTheta_n(\delta))$, $i=1,\ldots,k$)
satisfies $\bold{U}_n(\oTheta_n(\delta))\in\mathcal{F}_{\lfloor \oldz(\gamma_{\oE}+\delta)\nan\rfloor}$. We have
\begin{align}
\left\{A_n^*>\otaun(\gamma_{\oE}+\delta)\right\}&\subseteq\bigcap_{i=1}^{k}\bigcap_{t\leq\otaun(\gamma_{\oE}+\delta)}\left\{{\color{black}S_n^{(i)}}(t)+\oan{i}-{\color{black}U_n^{(i)}}(t)\geq 0\right\}\nonumber\\
&\subseteq\bigcap_{i=1}^{k}\{{\color{black}S_n^{(i)}}(\oTheta_n(\delta))+\oan{i}- {\color{black}U_n^{(i)}}(\oTheta_n(\delta))\ge 0\}\nonumber\\
&=\bigcup_{\bold{u}\in\mathcal{F}_{\lfloor \oldz(\gamma_{\oE}+\delta)\nan\rfloor}}
\bigcap_{i=1}^{k}\{{\color{black}S_n^{(i)}}(\oTheta_n(\delta))+\oan{i}-{\color{black}U_n^{(i)}}(\oTheta_n(\delta))\geq 0,\bold{U}_n(\oTheta_n(\delta))=\bold{u}\}.\nonumber
\end{align}
Therefore
\begin{align}
&P(A_n^*>\otaun(\gamma_{\oE}+\delta))\leq\sum_{\bold{u}\in\mathcal{F}_{\lfloor \oldz(\gamma_{\oE}+\delta)\nan\rfloor}}P\left(
\bigcap_{i=1}^{k}\{{\color{black}S_n^{(i)}}(\oTheta_n(\delta))+\oan{i}- {\color{black}U_n^{(i)}}(\oTheta_n(\delta))\geq 0\}\,\Big|\,\bold{U}_n(\oTheta_n(\delta))=\bold{u}\right)\nonumber\\
&\leq\prod_{i=1}^{k}\lfloor \oldz_i(\gamma_{\oE}+\delta)\nac{i}{n}\rfloor\times
\max_{\bold{u}\in\mathcal{F}_{\lfloor \oldz(\gamma_{\oE}+\delta)\nan\rfloor}}P\left(\bigcap_{i=1}^{k}\left\{{\color{black}S_n^{(i)}}\left(\sum_{i=1}^{k}u_i\right)+\oan{i}-u_i\ge 0\right\}
\,\Big|\,\bold{U}_n(\oTheta_n(\delta))=\bold{u}\right)\label{eq:ineqpIAN}\\
&\leq\prod_{i=1}^{k}\lfloor \oldz_i(\gamma_{\oE}+\delta)\nac{i}{n}\rfloor\nonumber\\
&\qquad\qquad\qquad
\times\max_{1\leq j\leq k}\max_{\bold{u}\in\mathcal{F}_{\lfloor \oldz(\gamma_{\oE}+\delta)\nan\rfloor}^{(j)}}P\left(\bigcap_{i=1}^{k}\left\{{\color{black}S_n^{(i)}}\left(\sum_{i=1}^{k}u_i\right)+\oan{i}-u_i\ge 0\right\}
\,\Big|\,\bold{U}_n(\oTheta_n(\delta))=\bold{u}\right).\label{eq:finlim}
\end{align}
Here, in the derivation of the inequality \eqref{eq:ineqpIAN} we exploited
relations
\[
\sum_{i=1}^{k}{\color{black}U_n^{(i)}}(\oTheta_n(\delta))=\oTheta_n(\delta)\quad\text{and}\quad
|\mathcal{F}_{\lfloor \oldz(\gamma_{\oE}+\delta)\nan\rfloor}|\leq\prod_{i=1}^{k}\lfloor \oldz_i(\gamma_{\oE}+\delta)\nac{i}{n}\rfloor.
\]
Note that, for fixed $j\in\{1,\ldots,k\}$ and
$\bold{u}\in\mathcal{F}_{\lfloor \oldz(\gamma_{\oE}+\delta)\nan\rfloor}^{(j)}$,
\begin{align}
&P\left(\bigcap_{i=1}^{k}\left\{{\color{black}S_n^{(i)}}\left(\sum_{i=1}^{k}u_i\right)+\oan{i}-u_i\ge 0\right\}
\,\Big|\,\bold{U}_n(\oTheta_n(\delta))=\bold{u}\right)\nonumber\\
&\qquad\qquad
\leq P\left({\color{black}S_n^{(j)}}\left(\sum_{i=1}^{k}u_i\right)+\oan{j}\ge\lfloor \oldz_j(\gamma_{\oE}+\delta)\nac{j}{n}\rfloor
\,\Big|\,\bold{U}_n(\oTheta_n(\delta))=\bold{u}\right)\nonumber\\
&\qquad\qquad
=P(\mathrm{Bin}(\on{j}-\oan{j},\oldpi(\bold u,\obq{j}))\geq\lfloor \oldz_j(\gamma_{\oE}+\delta)\nac{j}{n}\rfloor-\oan{j})\nonumber\\
&\qquad\qquad
\leq P(\mathrm{Bin}(\on{j}-\oan{j},\oldpi(\lfloor \oldz(\gamma_{\oE}+\delta)\nan\rfloor,\obq{j}))
\geq\lfloor \oldz_j(\gamma_{\oE}+\delta)\nac{j}{n}\rfloor-\oan{j}),\nonumber
\end{align}
where the latter inequality follows by the stochastic ordering property of the binomial distribution.
Combining this inequality with \eqref{eq:finlim} we have
\begin{align}
&P(A_n^*>\otaun(\gamma_{\oE}+ \delta))\leq\prod_{i=1}^{k}\lfloor \oldz_i(\gamma_{\oE}+\delta)\nac{i}{n}\rfloor\nonumber\\
&\qquad\qquad
\times\max_{1\leq j\leq k}P(\mathrm{Bin}(\on{j}-\oan{j},\oldpi(\lfloor \oldz(\gamma_{\oE}+\delta)\nan\rfloor,\obq{j}))
\geq\lfloor \oldz_j(\gamma_{\oE}+\delta)\nac{j}{n}\rfloor-\oan{j}).\nonumber
\end{align}
Since
\[
\prod_{i=1}^{k}\lfloor \oldz_i(\gamma_{\oE}+\delta)\nac{i}{n}\rfloor\sim (\nac{1}{n})^k,
\]
the claim then follows if we prove that, for an arbitrarily fixed $i\in\{1,\ldots,k\}$, quantity
$P(\mathrm{Bin}(\on{i}-\oan{i},\oldpi(\lfloor \oldz(\gamma_{\oE}+\delta)\nan\rfloor,\obq{i}))
\geq\lfloor \oldz_i(\gamma_{\oE}+\delta)\nac{i}{n}\rfloor-\oan{i})$ goes to zero exponentially fast with respect to $\nac{1}{n}$. For this
we exploit the same idea and similar computations as in the proof of \eqref{eq:Iconvinprob}, and so we omit many details. By Lemma \ref{le:aspi} and the definition of $\bm{\chi}$
we have
\begin{equation}\label{eq:equiv33}
(\on{i}-\oan{i})\oldpi(\lfloor \oldz(\gamma_{\oE}+\delta)\nan\rfloor,\obq{i})\sim_e
(\oldz_i(\gamma_{\oE}+\delta)-\alpha_i)\nac{i}{n}+\oldb_i(\oldz(\gamma_{\oE}+\delta))\nac{i}{n}.
\end{equation}
Therefore
\begin{align}
\frac{\lfloor \oldz_i(\gamma_{\oE}+\delta)\nac{i}{n}\rfloor-\oan{i}}{(\on{i}-\oan{i})\oldpi(\lfloor \oldz(\gamma_{\oE}+\delta)\nan\rfloor,\obq{i})}
\to\frac{\oldz_i(\gamma_{\oE}+\delta)-\alpha_i}
{\oldz_i(\gamma_{\oE}+\delta)-\alpha_i+\oldb_i(\oldz(\gamma_{\oE}+\delta))}.\label{limite}
\end{align}
By $\bold{(Sub)'}$ it follows that there exists $\delta_0>0$ such that
\[
\max_{1\leq i\leq k}\oldb_i(\oldz(\gamma_{\oE}+\delta))<0,\quad\text{for any $0<\delta \le \delta_0$.}
\]
Indeed, necessarily,
$J_{\bm\oldb}(\oldz(\gamma_{\oE}))\bold{\ov}_{PF}(\oldz(\gamma_{\oE}))=\lambda_{PF}(\oldz(\gamma_{\oE}))\bold{\ov}_{PF}(\oldz(\gamma_{\oE}))<\bold{0}$ (i.e., $\lambda_{PF}(\oldz(\gamma_{\oE}))<0$)
otherwise we would get a contradiction with \eqref{sub-cond_k}, since
the directional derivative of each function $\oldb_i$ along $\tilde{\mathcal{S}}$ is increasing.
Therefore, by a classical concentration inequality for the binomial distribution (see e.g. Eq. $(1.5)$ of Lemma 1.1 p. 16 in \cite{P}),
for all $n$ large enough, we have
\begin{align}
&P\Biggl(\mathrm{Bin}(\on{i}-\oan{i},\oldpi(\lfloor \oldz(\gamma_{\oE}+\delta)\nan\rfloor,\obq{i}))
\geq\lfloor \oldz_i(\gamma_{\oE}+\delta)\nac{i}{n}\rfloor-\oan{i}\Biggr)\nonumber\\
&\leq\exp\Biggl(-(\on{i}-\oan{i})\oldpi(\lfloor \oldz(\gamma_{\oE}+\delta)\nan\rfloor,\obq{i})\nonumber\\
&\,\,\,\,\,\,
\,\,\,\,\,\,
\,\,\,\,\,\,
\,\,\,\,\,\,
\,\,\,\,\,\,
\,\,\,\,\,\,
\,\,\,\,\,\,
\,\,\,\,\,\,
\,\,\,\,\,\,
\,\,\,\,\,\,
\,\,\,\,\,\,
\times H\left(\frac{\lfloor \oldz_i(\gamma_{\oE}+\delta)\nac{i}{n}\rfloor-\oan{i}}{(\on{i}-\oan{i})\oldpi(\lfloor \oldz(\gamma_{\oE}+\delta)\nan\rfloor,\obq{i})}\right)\Biggr),
\nonumber
\end{align}
where function $H$ is defined by \eqref{eq:H}. The claim (i.e. the exponential decreasing rate with respect to $\nac{1}{n}$) easily follows combining this inequality with
\eqref{eq:equiv33} and \eqref{limite}, and using that $H$ increases on $(1,\infty)$.
\\
\noindent$\square$

\subsection{Proof of Theorem \ref{teosupcrit}}

Let $\zeta:[0,\infty)\to\mathbb{R}^k$ be a parametrization of the curve with trace $\mathcal{C}\cup\tilde{\mathcal{S}}_{\mathrm{ext}}$, where  $\mathcal C$ is
defined in $(\bold{Sup})'$ (recall that by Theorem \ref{prop:sub} $({\mathcal Sup})$ is equivalent to $(\bold{Sup})'$) and
\[
\tilde{\mathcal{S}}_{\mathrm{ext}}:=\{\bold x\in\mathbb{R}^k\setminus\oD:\,\, \bold{x}=\bold{x}(\ote)+\theta\bold{\ov}_{PF}(\bold{x}(\ote)),\,\theta\geq 0\}
\]
is the extension of $\tilde{\mathcal{S}}$ outside $\oD$. Let $\otaun(\gamma)$ be defined by \eqref{eq:tau} but with $\gamma\in [0,\infty)$. Similarly to the proof of
Theorem \ref{teo-blocksub} we put $\otaun([0,\infty))=\{t_n^{(s)}\}_{s\in\mathbb N\cup\{0\}}$ and,
without loss of generality, we assume $t_n^{(0)}:=0<t_n^{(1)}<\ldots<t_n^{(s)}<\ldots$. In this context, we define $\otaumn(t_n^{(s)})$, $s\in\mathbb{N}\cup\{0\}$,
as in \eqref{eq:taumeno1} (clearly with $[0,\infty)$ in place of $[0,1]$), $\ozeta_n^{(i)}(t_n^{(s)})$, $i=1,\ldots,k$, $s\in\mathbb{N}\cup\{0\}$, as in \eqref{eq:zetasullimag}, and we extend the definition of
$\ozeta_n^{(i)}(\cdot)$ to any $t\in\mathbb{N}\setminus\otaun([0,\infty))$ as in \eqref{eq:zetaext}. We define the random variables $T_n'$ as in \eqref{eq:tnprimo}
(with $n$ in place of $\otaun(1)$).
As in the proof of Theorem \ref{teo-blocksub},
for $t\leq T'_n$, we consider the \lq\lq policy" ${\color{black}C_n^{(i)}}$ defined by \eqref{eq:policydet}. We proceed dividing the proof in three main steps.\\
\noindent{\it Step\,\,1:\,\,outline\,\,of\,\,the\,\,proof.} For $\varepsilon\geq 1$, we have
\[
P(n-A_n^*>\varepsilon n)=0,
\]
and the claim is obvious. Since $A_n^*\geq T_n'$, for $\varepsilon\in (0,1)$,
we have
\[
P(n-A_n^*>\varepsilon n)\leq 1-P(T_n'\geq (1-\varepsilon)n),
\]
and therefore it suffices to show that
\[
P(T_n'\geq (1-\varepsilon)n)\geq 1-O(\mathrm{e}^{-c(\varepsilon)g_n^{(1)}}).
\]
For a suitable positive constant $\oc{1}>0$ and all $n$ large enough, we have
\begin{align}
P(T_n'\geq (1-\varepsilon)n)&=P(T_n'\geq (1-\varepsilon)n\,|\,T_n'\geq \lfloor \oc{1} \nac{1}{n}\rfloor+1)P(T_n'\geq\lfloor \oc{1} \nac{1}{n}\rfloor+1\,|\,T_n'\geq {\color{black}\hat{t}_n^{({\color{black}\overline{d}})}})\nonumber\\
&\quad\quad\quad
\times P(T_n'\geq {\color{black}\hat{t}_n^{({\color{black}\overline{d}})}}), \nonumber
\end{align}
where ${\color{black}\hat{t}_n^{({\color{black}\overline{d}})}}$ is defined by \eqref{eq:tmin} with $h={\color{black}\overline{d}}$. Since $P(T_n'\geq {\color{black}\hat{t}_n^{({\color{black}\overline{d}})}})\geq T_{1,n}({\color{black}\overline{d}})$,
where $T_{1,n}({\color{black}\overline{d}})$ is defined by \eqref{eq:hug0n}, arguing as in the proof of Theorem
\ref{teo-blocksub} one has $T_{1,n}({\color{black}\overline{d}})\geq 1-O(\mathrm{e}^{-c(\varepsilon)g_n^{(1)}})$,
and so $P(T_n'\geq {\color{black}\hat{t}_n^{({\color{black}\overline{d}})}})\geq 1-O(\mathrm{e}^{-c(\varepsilon)g_n^{(1)}})$. In the next two steps we shall show that
\begin{equation}\label{eq:Ip}
P(T_n'>\lfloor (1-\varepsilon)n\rfloor\,|\,T_n'>\lfloor \oc{1} \nac{1}{n}\rfloor)\geq 1-O(\mathrm{e}^{-c(\varepsilon)g_n^{(1)}})
\end{equation}
and
\begin{equation}\label{eq:2p}
P(T_n'>\lfloor \oc{1} \nac{1}{n}\rfloor\,|\,T_n'>{\color{black}\hat{t}_n^{({\color{black}\overline{d}})}}-1)\geq 1-O(\mathrm{e}^{-c(\varepsilon)g_n^{(1)}}),
\end{equation}
concluding the proof.\\
\noindent{\it Step\,\,2:\,\,proof\,\,of\,\,\eqref{eq:Ip}.}\\
For all $n$ large enough, we have
\begin{align}
P(T_n'>\lfloor(1-\varepsilon)n\rfloor\,|\,T_n'>\lfloor \oc{1} \nac{1}{n}\rfloor)&=P(T_n'>\lfloor(1-\varepsilon)n\rfloor\,|\,T_n'>\lfloor (\op{1})^{-1}\rfloor)\nonumber\\
&\quad\quad\quad
\times P(T_n'>\lfloor (\op{1})^{-1}\rfloor\,|\,T_n'>\lfloor \oc{1} \nac{1}{n}\rfloor). \nonumber 
\end{align}
The claim clearly follows if we prove
\begin{equation}\label{eq:ianv1}
P(T_n'>\lfloor (\op{1})^{-1}\rfloor\,|\,T_n'>\lfloor \oc{1} \nac{1}{n}\rfloor)\geq 1-O(\mathrm{e}^{-c(\varepsilon)g_n^{(1)}})
\end{equation}
and
\begin{equation}\label{eq:ianv2}
P(T_n'>\lfloor(1-\varepsilon)n\rfloor\,|\,T_n'>\lfloor (\op{1})^{-1}\rfloor)\geq 1-O(\mathrm{e}^{-c(\varepsilon)g_n^{(1)}}).
\end{equation}
We divide the proof of this step in two parts. In the Step 2.1 we prove \eqref{eq:ianv1} and in the Step 2.2 we prove \eqref{eq:ianv2}.\\
\noindent{\it Step\,\,2.1:\,\,proof\,\,of\,\,\eqref{eq:ianv1}.} For $n$ large enough,
define
\[
{\color{black}l_n}:=\min\{\ell\geq\lceil \oc{1}\rceil:\,\,\op{1}{\color{black}m_n^{(\ell)}}\geq 1\},\quad\text{where ${\color{black}m_n^{(\ell)}}:=\lfloor \oc{1}^{\ell/\lceil \oc{1}\rceil}\nac{1}{n}\rfloor$}
\]
and $\lceil x\rceil$ denotes the smallest integer greater than or equal to $x\in\mathbb R$. Since ${\color{black}m_n^{({{\color{black}l_n}})}}\geq\lfloor (\op{1})^{-1}\rfloor$, we have
\[
[\lfloor \oc{1} \nac{1}{n}\rfloor,\lfloor (\op{1})^{-1}\rfloor]\cap\mathbb N\subseteq\bigcup_{\ell=\lceil \oc{1}\rceil}^{{\color{black}l_n}-1}[{\color{black}m_n^{(\ell)}},{\color{black}m_n^{({\ell+1})}}]\cap\mathbb N
\]
and so, for all $n$ large enough, using relation \eqref{eq:defTgretat} and the fact paths ${\color{black}A_n^{(i)}}(\cdot)$ and functions $\ozeta_n^{(i)}(\cdot)$
are non-decreasing, we have
\begin{align}
&P(T_n'>\lfloor (\op{1})^{-1}\rfloor\,|\,T_n'>\lfloor \oc{1} \nac{1}{n}\rfloor)
\geq\prod_{\ell=\lceil \oc{1}\rceil}^{{\color{black}l_n}-1}P(T_n'>{\color{black}m_n^{({\ell+1})}}\,|\,T_n'>{\color{black}m_n^{(\ell)}})\nonumber\\
&=\prod_{\ell=\lceil \oc{1}\rceil}^{{\color{black}l_n}-1}P(\bold{A}_n(s)\ge\bm{\ozeta}_{n}(s+1)\,\,
\forall 1\le s\le {\color{black}m_n^{({\ell+1})}}\,|\,\bold{U}_n(s)=\bm{\ozeta}_{n}(s)\,\,\forall 1\le s\le{\color{black}m_n^{({\ell})}}+1)
\nonumber\\
&=\prod_{\ell=\lceil \oc{1}\rceil}^{{\color{black}l_n}-1}P(\bold{A}_{n}(s)\ge\bm{\ozeta}_n(s+1)\,\,
\forall {\color{black}m_n^{({\ell})}}+1\le s\le {\color{black}m_n^{({\ell+1})}}\,|\,\bold{U}_{n}(s)=\bm{\ozeta}_{n}(s)\,\,\forall 1\le s\le{\color{black}m_n^{({\ell})}}+1)
\nonumber\\
&\geq\prod_{\ell=\lceil \oc{1}\rceil}^{{\color{black}l_n}-1}P(\bold{A}_n({\color{black}m_n^{({\ell})}}+1)\ge\bm{\ozeta}_{n}({\color{black}m_n^{({\ell+1})}}+1)
\,|\,\bold{U}_n(s)=\bm{\ozeta}_{n}(s)\,\,\forall 1\le s\le{\color{black}m_n^{({\ell})}}+1)\nonumber\\
&=\prod_{\ell=\lceil \oc{1}\rceil}^{{\color{black}l_n}-1}P(\bold{A}_n({\color{black}m_n^{({\ell})}}+1)\ge\bm{\ozeta}_{n}({\color{black}m_n^{({\ell+1})}}+1)
\,|\,\bold{U}_n({\color{black}m_n^{({\ell})}}+1)=\bm{\ozeta}_{n}({\color{black}m_n^{({\ell})}}+1))\label{eq:markov}\\
&\geq\prod_{\ell=\lceil \oc{1}\rceil}^{{\color{black}l_n}-1}P(\bold{S}'_{n}({\color{black}m_n^{({\ell})}}+1)\ge\bm{\ozeta}_{n}({\color{black}m_n^{({\ell+1})}}+1)
\,|\,\bold{U}_{n}({\color{black}m_n^{({\ell})}}+1)=\bm{\ozeta}_{n}({\color{black}m_n^{({\ell})}}+1))\label{eq:Sprimo}\\
&=\prod_{\ell=\lceil \oc{1}\rceil}^{{\color{black}l_n}-1}P(\mathrm{Bin}(\on{i},\oldpi(\bm{\ozeta}_{n}({\color{black}m_n^{({\ell})}}+1),\obq{i}))\geq\ozeta_n^{(i)}({\color{black}m_n^{({\ell+1})}}+1)
\,\,\forall 1\le i\le k)\label{eq:binagain}\\
&\geq\prod_{\ell=\lceil \oc{1}\rceil}^{{\color{black}l_n}-1}\left(1-\sum_{i=1}^{k}P(\mathrm{Bin}(\on{i},\oldpi(\bm{\ozeta}_{n}({\color{black}m_n^{({\ell})}}+1),\obq{i}))<\ozeta_n^{(i)}({\color{black}m_n^{({\ell+1})}}+1)
\right),\label{eq:unionbound}
\end{align}
where in \eqref{eq:markov} we used \eqref{eq:july5},
in \eqref{eq:Sprimo} we used that
\[
S_{i}^{(n)'}(s):=\sum_{v\in G_j}\bold{1}\{Y_{v}\leq s\}\leq {\color{black}S_n^{(i)}}(s)+\oan{i}={\color{black}A_n^{(i)}}(s)
\]
and we put $\bold{S}'_{n}:=(S^{(n)'}_1,\ldots,S^{(n)'}_k)$, in \eqref{eq:binagain} we used \eqref{eq:bin}, and \eqref{eq:unionbound} follows by the union bound.
Choosing $\oc{1}$ large enough and arguing as in the proof of relation $(59)$ in \cite{TGL}, for all $i$, $\ell$
and $n$ large enough we get
\begin{equation}\label{eq:disSPA}
P(\mathrm{Bin}(\on{{i}},\oldpi(\bm{\ozeta}_n({\color{black}m_n^{({\ell})}}+1),\obq{i}))<\ozeta_n^{(i)}({\color{black}m_n^{({\ell+1})}}+1))\leq\mathrm{e}^{-\oc{2} \nac{1}{n}}\mathrm{e}^{-(\ell-\lceil \oc{1}\rceil)\oc{3} \nac{1}{n}},
\end{equation}
for some positive constants $\oc{2},\oc{3}>0$. Using that
\begin{equation*}
\text{for $x_1,\ldots,x_n\in (0,1)$ it holds $\prod_{i=1}^{n}(1-x_i)>1-\sum_{i=1}^{n}x_i$}
\end{equation*}
and \eqref{eq:disSPA}, by \eqref{eq:unionbound}, for all $n$ large enough we have
\begin{align}
P(T_n'>\lfloor (\op{1})^{-1}\rfloor\,|\,T_n'>\lfloor \oc{1} \nac{1}{n}\rfloor)&\geq 1-\sum_{i=1}^{k}\sum_{\ell=\lceil \oc{1}\rceil}^{{\color{black}l_n}-1}
P(\mathrm{Bin}(\on{i},\oldpi(\bm{\ozeta}_n({\color{black}m_n^{(\ell)}}+1),\obq{i}))<\ozeta_n^{(i)}({\color{black}m_n^{({\ell+1})}}+1))\nonumber\\
&\geq 1-\oc{4}\mathrm{e}^{-\oc{2} \nac{1}{n}},\nonumber
\end{align}
for some positive constant $\oc{4}>0$, which yields \eqref{eq:ianv1}.\\
\noindent{\it Step\,\,2.2:\,\,proof\,\,of\,\,\eqref{eq:ianv2}.} Let $\oc{5}\in (0,1)$ be a small positive constant such that,
for all $n$ large enough $P(\mathrm{Bin}(\lfloor (\op{1})^{-1},\op{1})\geq r)\geq 2c _5$
(see e.g. the proof of Lemma 8.2 Case 3 p. 26 in \cite{JLTV}). For all $n$ large enough we have
\begin{align}
&P(T_n'>\lfloor(1-\varepsilon)n\rfloor\,|\,T_n'>\lfloor (\op{1})^{-1}\rfloor)\nonumber\\
&\qquad\qquad\qquad
=P(T_n'>\lfloor(1-\varepsilon)n\rfloor\,|\,T_n'>\lfloor \oc{5} n\rfloor)
P(T_n'>\lfloor \oc{5} n\rfloor\,|\,T_n'>\lfloor(\op{1})^{-1}\rfloor).\label{eq:13J0}
\end{align}
Arguing similarly to the derivation of \eqref{eq:unionbound}, we have
\begin{align}
&P(T_n'>\lfloor \oc{5} n\rfloor\,|\,T_n'>\lfloor(\op{1})^{-1}\rfloor)\nonumber\\
&
=P(\bold{A}_n(s)\ge\bm{\ozeta}_{n}(s+1)\,\,
\forall 1\le s\le\lfloor \oc{5} n\rfloor\,|\,\bold{U}_n(s)=\bm{\ozeta}_{n}(s)\,\,\forall 1\le s\le\lfloor(\op{1})^{-1}\rfloor +1)\nonumber\\
&
=P(\bold{A}_n(s)\ge\bm{\ozeta}_{n}(s+1)\,\,
\forall \lfloor(\op{1})^{-1}\rfloor+1\le s\le\lfloor \oc{5} n\rfloor\,|\,\bold{U}_{n}(s)=\bm{\ozeta}_{n}(s)\,\,
\forall 1\le s\le\lfloor(\op{1})^{-1}\rfloor +1)\nonumber\\
&
\geq P(\bold{A}_n(\lfloor(\op{1})^{-1}\rfloor+1)\geq\bm{\ozeta}_n(\lfloor \oc{5} n\rfloor+1)
\,|\,\bold{U}_{n}(s)=\bm{\ozeta}_{n}(s)\,\,\forall 1\le s\le\lfloor(\op{1})^{-1}\rfloor+1)\nonumber\\
&
=P(\bold{A}_{n}(\lfloor(\op{1})^{-1}\rfloor+1)\geq\bm{\ozeta}_n(\lfloor \oc{5} n\rfloor+1)
\,|\,\bold{U}_{n}(\lfloor(\op{1})^{-1}\rfloor+1)=\bm{\ozeta}_{n}(\lfloor(\op{1})^{-1}\rfloor+1))\nonumber\\
&
\geq 1-\sum_{i=1}^{k}P(\mathrm{Bin}(\on{i},\oldpi(\bm{\ozeta}_n(\lfloor(\op{1})^{-1}\rfloor+1),\obq{i}))<\ozeta_n^{(i)}(\lfloor \oc{5} n\rfloor+1)).\label{eq:13J1}
\end{align}
Similarly, we get
\begin{align}
&P(T_n'>\lfloor(1-\varepsilon)n\rfloor\,|\,T_n'>\lfloor \oc{5} n\rfloor)\nonumber\\
&\qquad\qquad\qquad
\geq 1-\sum_{i=1}^{k}P(\mathrm{Bin}(\on{i},\oldpi(\bm{\ozeta}_n(\lfloor \oc{5} n\rfloor+1),\obq{i}))<\ozeta_n^{(i)}(\lfloor(1-\varepsilon)n\rfloor+1)).\label{eq:13J3}
\end{align}
The following inequalities are proved in Step 4 of the proof of Proposition 4.1 in \cite{TGL} and
hold for any $i$ and all $n$ large enough:
\begin{align}
&P(\mathrm{Bin}(\on{i},\oldpi(\bm{\ozeta}_n(\lfloor(\op{1})^{-1}\rfloor+1),\obq{i}))<\ozeta_n^{(i)}(\lfloor \oc{5} n\rfloor+1))
\leq \oc{6}\mathrm{e}^{-\oc{7} n},\quad\text{for some constants $\oc{6},\oc{7}>0$}\nonumber\\
&P(\mathrm{Bin}(\on{i},\oldpi(\bm{\ozeta}_n(\lfloor \oc{5} n\rfloor+1),\obq{i}))<\ozeta_n^{(i)}(\lfloor(1-\varepsilon)n\rfloor+1))
\leq\mathrm{e}^{-\oc{8} \nac{1}{n}},\quad\text{for some constant $\oc{8}>0$.}\nonumber
\end{align}
These relations and \eqref{eq:13J0}, \eqref{eq:13J1} and \eqref{eq:13J3} clearly imply \eqref{eq:ianv2}.\\
\noindent{\it Step\,\,3:\,\,proof\,\,of\,\,\eqref{eq:2p}.}\\
Let $\gamma'\in (0,\infty)$ be such that
\[
\otaun(\gamma')=\sum_{i=1}^{k}\lfloor \oldz_i(\gamma')\nac{i}{n}\rfloor\geq\lfloor \oc{1} \nac{1}{n}\rfloor+1.
\]
For all $n$ large enough, we have
\begin{align}
&P(T_n'>\lfloor \oc{1} \nac{1}{n}\rfloor\,|\,T_n'>{\color{black}\hat{t}_n^{({\color{black}\overline{d}})}}-1)\nonumber\\
&\geq\prod_{s={\color{black}\hat{t}_n^{({\color{black}\overline{d}})}}}^{\lfloor \oc{1} \nac{1}{n}\rfloor-1}
\Biggl(1-\sum_{i=1}^{k}P(\mathrm{Bin}(\on{i}-\oan{i},\oldpi(\bm{\ozeta}_n(s),\obq{i}))<\ozeta_n^{(i)}(s+1)-\oan{i})
\Biggr)\nonumber\\
&\geq\prod_{s={\color{black}\hat{t}_n^{({\color{black}\overline{d}})}}}^{\otaun(\gamma')-1}
\Biggl(1-\sum_{i=1}^{k}P(\mathrm{Bin}(\on{i}-\oan{i},\oldpi(\bm{\ozeta}_n(s),\obq{i}))<\ozeta_n^{(i)}(s+1)-\oan{i})
\Biggr)=:\widetilde{T}_{2,n}({\color{black}\overline{d}}).\nonumber
\end{align}
The claim follows if we prove that $\widetilde{T}_{2,n}({\color{black}\overline{d}})\geq 1-O(\mathrm{e}^{-c(\varepsilon)g_n^{(1)}})$. This
can be verified along the same lines as in  the proof of relation $T_{2,n}({\color{black}\overline{d}})\geq 1-O(\mathrm{e}^{-c(\varepsilon)g_n^{(1)}})$, provided within  the proof of Theorem \ref{teo-blocksub}.
Hereon, we skip many details and highlight the main differences. For $\varepsilon>0$, let $\bar{\mathcal C}_\varepsilon$ be the $\varepsilon$-thickening
of compact set $\bar{\mathcal C}:=\oldz([\gamma_0^{{\color{black}\overline{d}}},\gamma'])$. As in the proof of Theorem \ref{teo-blocksub}, one has that
there exists $\varepsilon_0>0$ small enough so that $\bar{\mathcal C}_{\varepsilon_0}\subset\overset\circ{\oE}$ and it can be shown that
there exists $\bar n$ such that $\hat{\bm{\ozeta}}_n(s)\in\bar{\mathcal C}_{\varepsilon_0}$ for any $n>\bar n$ and any ${\color{black}\hat{t}_n^{({\color{black}\overline{d}})}}\leq s\leq\otaun(\gamma')-1$. Then, for any $\varepsilon>0$
and all $n$ large enough we get $\widetilde{T}_{2,n}({\color{black}\overline{d}})\geq (1-\mathrm{sup}_n)^{\oc{9}\nac{1}{n}}$, where $\mathrm{sup}_n$ is defined as in the proof of
Theorem \ref{teo-blocksub} and $\oc{9}>0$ is a positive constant. By assumption (${\mathcal Sup}$) (which is equivalent to $(\bold{Sup})'$) we have $\min_{\bold z\in \bar{\mathcal C}_{\varepsilon_0}}\oldb_1(\bold z)>0$.
Then one can show that $\mathrm{sup}_n=O(\mathrm{e}^{-\oc{10}\nac{1}{n}})$, for some constant $\oc{10}>0$, exactly as in the proof of Theorem \ref{teo-blocksub}, and the proof is completed.

\subsection{Proof of Theorem \ref{prop:sub}}\label{subsec:thm61}

First, we  introduce  the following additional conditions:
\\
\\
\noindent$(\bold{Sub})''$: Either $(\bold{Sub})(i)$ or $(\bold{Sub})(ii)$ or $(\bold{Sub})(iii)$ holds.\\
\noindent$(\bold{Sub})'''$: There exists $\bold x\in\oD$ such that ${\bm\oldb}(\bold x)\leq\bold 0$ and ${\bm\oldb}(\bold x)\neq\bold 0$;
\\
\\
where
\\
\\
\noindent$(\bold{Sub})(i)$: $|\mathcal Z|\geq 2$.\\
\noindent$(\bold{Sub})(ii)$: $\mathcal{Z}=\{\bold z\}$, $\bold{z}\in\overset\circ{\oD}$ and $\mathrm{det}J_{{\bm\oldb}}(\bold z)\neq 0$.\\
\noindent$(\bold{Sub})(iii)$: $\mathcal{Z}=\{\bold z\}$, $\bold{z}\in\overset\circ{\oD}$, $\mathrm{det}J_{{\bm\oldb}}(\bold z)=0$ and $\lambda_{PF}(\bold z)\neq 0$.
\\
\\
and:
\\
\\
\noindent$(\bold{I})$: $\mathcal{Z}=\{\bold z\}$ and $\bold{z}\in\ocalD$.
\\
\\
For reader's convenience, we summarized in Table \ref{table1} all the different \lq\lq critical" conditions that have been introduced in this paper.
The proof of Theorem \ref{prop:sub} is based on the following lemmas, which will be proved later on.

\begin{table}[t]
\begin{center}
\caption{Main \lq\lq critical" assumptions}\label{table1}
\begin{tabular}{|c|c|c|c|}
\hline
$({\mathcal Sub})$ & \multicolumn{2}{c|}{ $\min_{\bold x\in\oD_{\bm\oldb}}\oldb_1(\bold{x})<0$}  & \multirow{ 3}{*}{$\oD_{{\bm\oldb}}:=\{\bold x\in\oD:\,\,\oldb_1(\bold x)=\ldots=\oldb_k(\bold x)\}$ } \\
\cline{1-3}
$({\mathcal Crit})$ & \multicolumn{2}{c|}{$\min_{\bold x\in\oD_{\bm\oldb}}\oldb_1(\bold{x})=0$}  &  \\
\cline{1-3}
$({\mathcal Sup})$ &   \multicolumn{2}{c|}{$\min_{\bold x\in\oD_{\bm\oldb}}\oldb_1(\bold x)>0$}   &  \\
\hline
\multirow{ 3}{*}{$(\bold{Sub}')$} & \multirow{ 3}{*}{  $\min_{\bold x\in\mathcal{C}_{\mathcal E}\cup \tilde{\mathcal{S}  }}\min_{1\leq i\leq k} \oldb_i(\bold x)<0$ }&
\multirow{ 2}{*}{ $\ote=+\infty$ }    &  \multirow{ 2}{*}{$\mathcal S:= \cup_{1\le h\le {\color{black}\overline{d}}} \mathcal S_h$}\\
 & & &  \multirow{ 3}{*}{$\mathcal S_h:= \{ \bold{x}_0^{(h)} + (1-\theta) \bold{x}_0^{(h-1)}, \; \theta \in [0,1]\}$ } \\
 & & \multirow{ 2}{*}{$\bold{x}(\ote)\in\mathcal{Z}\cap\overset\circ{\oD}$ } & \\
 \cline{1-2}
 \multirow{ 3}{*}{$(\bold{Crit}')$} &   \multirow{ 3}{*}{$\min_{\bold x\in\mathcal{C}_{\mathcal E}\cup \tilde{\mathcal{S}  }}\min_{1\leq i\leq k} \oldb_i(\bold x)=0$} &   &
 \multirow{ 4}{*}{$\tilde{\mathcal{S}}:=\{ \bold{x}(\ote)+\theta\bold{\ov}_{PF}(\bold{x}(\ote))\in \oD,\,\theta\geq 0\}$}\\
 & &  \multirow{ 2}{*}{$\mathcal{C}:=\mathcal{S}\cup\mathcal{C}_\oE\cup\tilde{\mathcal S}$}  &  \\
 & & & \\
 \cline{1-3}
\multirow{ 4}{*}{$(\bold{Sup}')$} & \multirow{ 4}{*}{ $\min_{\bold x\in\mathcal{C}_{\mathcal E}}\min_{1\leq i\leq k} \oldb_i(\bold x)>0$}& $\ote<\infty$
& \multirow{ 4}{*}{ $\mathcal C\subset\oD$, $\bold{0}\in\mathcal{C}$, $\mathcal{C}\cap\ocalD\neq\emptyset$}\\
 & & $\bold{x}(\ote)\in\ocalD$,   &  \\
 & & ${\bm\oldb}(\bold{x}(\ote))>\bold 0$& \\
& &  $\mathcal{C}:=\mathcal{S}\cup\mathcal{C}_\oE$ & \\
 \hline
 $(\bold{Sub})$    &\multicolumn{3}{c|}{ $\exists \bold{x}\in \oD: {\bm\oldb}( \bold{x} )< \bold 0$ } \\
 \hline
 $(\bold{Crit})$ & \multicolumn{3}{c|}{$\mathcal{Z}=\{\bold z\}$,   $\bold{z}\in\overset\circ{\oD}$,    $\mathrm{det}J_{{\bm\oldb}}(\bold z)=0$,  $\lambda_{PF}(\bold z)=0$} \\
  \hline
 $(\bold{Sup})$ & \multicolumn{3}{c|}{$\mathcal{Z}=\emptyset$}\\
 \hline
  $(\bold{Sub})(i)$ &  \multicolumn{3}{c|}{$|\mathcal Z|\geq 2$}\\
 \hline
 $(\bold{Sub})(ii)$  & \multicolumn{3}{c|}{ $\mathcal{Z}=\{\bold z\}$, $\bold{z}\in\overset\circ{\oD}$ and $\mathrm{det}J_{{\bm\oldb}}(\bold z)\neq 0$}\\
 \hline
 \noindent$(\bold{Sub})(iii)$ &\multicolumn{3}{c|}{ $\mathcal{Z}=\{\bold z\}$, $\bold{z}\in\overset\circ{\oD}$, $\mathrm{det}J_{{\bm\oldb}}(\bold z)=0$ and $\lambda_{PF}(\bold z)\neq 0$}\\
 \hline
 $(\bold{I})$ & \multicolumn{3}{c|}{ $\mathcal{Z}=\{\bold z\}$ and $\bold{z}\in\ocalD$}\\
 \hline
 $(\bold{Sub})''$ & \multicolumn{3}{c|}{Either $(\bold{Sub})(i)$ or $(\bold{Sub})(ii)$ or $(\bold{Sub})(iii)$ holds }\\
\hline
 $(\bold{Sub})'''$ & \multicolumn{3}{c|}{ $\exists \bold{x}\in \oD: {\bm\oldb}( \bold{x} )\le \bold 0, \;\, {\bm\oldb}( \bold{x} )\neq \bold 0$ }\\
 \hline
 \end{tabular}
\end{center}
\end{table}

\begin{Lemma}\label{le:Sup}
Under the assumptions of Theorem \ref{prop:sub}, we have that  conditions $(\bold{Sup})$ and $(\bold{Sup})'$ are equivalent.
\end{Lemma}

\begin{Lemma}\label{le:esistecurva}
Under the assumptions of Theorem \ref{prop:sub}, if moreover $\mathcal{Z}\neq\emptyset$,
then $\ote=+\infty$ and $\bold{x}(\ote)\in\mathcal{Z}\cap\oE$.
\end{Lemma}

\begin{Lemma}\label{le:subzeri}
Under the assumptions of Theorem \ref{prop:sub}, we have $(\bold{Sub})\Rightarrow\mathcal Z\neq\emptyset$.
\end{Lemma}

\begin{Lemma}\label{le:subsub'}
Under the assumptions of Theorem \ref{prop:sub}, we have $({\mathcal Sub})\Rightarrow (\bold{Sub}) \Rightarrow(\bold{Sub})'$.
\end{Lemma}

\begin{Lemma}\label{le:critcrit'}
Under the assumptions of Theorem \ref{prop:sub}, we have $(\bold{Crit})\Rightarrow (\bold{Crit})'$.
\end{Lemma}

\begin{Lemma}\label{le:subequiv}
Under the assumptions of Theorem \ref{prop:sub}, we have that conditions $(\bold{Sub})$, $(\bold{Sub})''$ and $(\bold{Sub})'''$ are equivalent.
\end{Lemma}

\begin{Lemma}\label{le:Ithensub}
Under the assumptions of Theorem \ref{prop:sub}, we have $(\bold I)\Rightarrow (\bold{Sub})$.
\end{Lemma}

\noindent{\it Proof\,\,of\,\,Theorem\,\,\ref{prop:sub}.}\\
\noindent{\it Proof\,\,of\,\,$(i)$.} Due to Lemma \ref{le:subsub'} we only have to show
$(\bold{Sub})'\Rightarrow (\bold{Sub})\Rightarrow ({\mathcal Sub})$. We start proving
$(\bold{Sub})'\Rightarrow (\bold{Sub})$. We reason by contradiction and suppose that the claim is false.
Then by Lemma \ref{le:subequiv} we have that $(\bold{Sub})''$ does not hold. So either $(\bold{I})$ or $(\bold{Crit})$ holds
(indeed $(\bold{Sup})$ can not hold by Lemma \ref{le:Sup} since we are assuming $(\bold{Sub})'$). If $(\bold{I})$ holds, then $\bold{x}(\ote)\in\ocalD$ by Lemma \ref{le:esistecurva},
and this contradicts $(\bold{Sub})'$. If $(\bold{Crit})$ holds, then, by Lemma \ref{le:critcrit'}, $(\bold{Crit})'$ holds, and this again contradicts $(\bold{Sub})'$. We now prove
$(\bold{Sub})\Rightarrow ({\mathcal Sub})$. We reason again by contradiction and suppose that either $({\mathcal Sup})$ or $({\mathcal Crit})$ holds. By Lemma
\ref{le:subzeri} we have $\mathcal Z\neq\emptyset$ and so $\min_{\bold{x}\in\oD_{\bm\oldb}}\oldb_1(\bold x)\leq 0$. This excludes $({\mathcal Sup})$.
We shall show in part $(iii)$ of this proof that $({\mathcal Crit})\Rightarrow (\bold{Crit})$. We get a contradiction noticing that
conditions $(\bold{Crit})$ and $(\bold{Sub})''$ are not compatible and
by Lemma \ref{le:subequiv} we have that $(\bold{Sub})''$ holds.\\
\noindent{\it Proof\,\,of\,\,$(ii)$.} Due to Lemma \ref{le:Sup} it suffices to show that $({\mathcal Sup})$ and $(\bold{Sup})$ are equivalent. We start proving
$({\mathcal Sup})\Rightarrow (\bold{Sup})$. By assumption $\min_{\bold x\in\oD_{\bm\oldb}}\oldb_1(\bold{x})>0$, so
there are no zeros of $\bm{\rho}$ in $\oD_{\bm\oldb}$, and therefore, there are no zeros of $\bm{\rho}$ in $\oD$. We now prove
$(\bold{Sup})\Rightarrow ({\mathcal Sup})$. Reasoning by contradiction assume that $({\mathcal Sup})$ does not hold. Then either $({\mathcal Sub})$ or $({\mathcal Crit})$ holds.
If $({\mathcal Sub})$ holds, then by Lemma \ref{le:subsub'} $(\bold{Sub})$ holds. By Lemma \ref{le:subequiv} we then get $(\bold{Sub})''$, which clearly contradicts $(\bold{Sup})$.
If $({\mathcal Crit})$ holds, then $\bm{\rho}$ has zeros on $\oD$. Therefore, either $(\bold{Sub})''$ or $(\bold{Crit})$ or $(\bold I)$ holds. Due to
Lemmas \ref{le:subequiv} and \ref{le:Ithensub} and the fact that $(\bold{Sub})''$ contradicts $(\bold{Sup})$, neither $(\bold{Sub})''$ nor $(\bold{I})$ can hold.
The claim follows noticing that also $(\bold{Crit})$ contradicts $(\bold{Sup})$.\\
\noindent{\it Proof\,\,of\,\,$(iii)$.} Due to Lemma \ref{le:critcrit'} we need to prove $(\bold{Crit})'\Rightarrow (\bold{Crit})$ and $(\bold{Crit})\Leftrightarrow ({\mathcal Crit})$.
We start proving $(\bold{Crit})'\Rightarrow (\bold{Crit})$.
Reasoning by contradiction assume that $(\bold{Crit})$ does not hold. Then we distinguish two cases: either $\mathcal{Z}=\emptyset$ or $\mathcal{Z}\neq\emptyset$.
The first case implies $(\bold{Sup})$, i.e., by Lemma \ref{le:Sup} $(\bold{Sup})'$. This clearly contradicts $(\bold{Crit})'$. In the second case we have that one of the following four conditions holds:
$(\bold{Sub})(i)$, $(\bold{Sub})(ii)$, $(\bold{Sub})(iii)$, $(\bold{I})$. If one of the first three conditions holds, then by Lemmas \ref{le:subsub'} and \ref{le:subequiv}, we have
that condition $(\bold{Sub})'$ holds, and this contradicts $(\bold{Crit})'$. If $(\bold{I})$ holds, then
$\bold{x}(\ote)\in\ocalD$ by Lemma \ref{le:esistecurva}, and this again contradicts $(\bold{Crit})'$. We now prove
$(\bold{Crit})\Rightarrow ({\mathcal Crit})$.
We reason again by contradiction and assume that $({\mathcal Crit})$ does not hold. Then either $({\mathcal Sub})$ or $({\mathcal Sup})$ holds. If $({\mathcal Sub})$ holds, then by Lemmas \ref{le:subsub'} and \ref{le:subequiv}
$(\bold{Sub})''$ holds, which contradicts $(\bold{Crit})$. If $({\mathcal Sup})$ holds, then $\bm{\rho}$ has no zeros on $\oD$,
which again contradicts $(\bold{Crit})$. Finally we prove $({\mathcal Crit})\Rightarrow (\bold{Crit})$.
Clearly $({\mathcal Crit})$ implies that $\bm{\rho}$ has zeros on $\oD$. Therefore either $(\bold{Sub})''$
or $(\bold{I})$ or $(\bold{Crit})$ holds. So, by Lemmas \ref{le:subequiv} and \ref{le:Ithensub}, either $(\bold{Sub})$ or $(\bold{Crit})$ holds.
We show that if $(\bold{Sub})$ holds, then we get a contradiction. Hereon, we explicit the dependence of ${\bm\oldb}$ and $\oD_{\bm\oldb}$ on $\bm{\alpha}:=(\alpha_1,\ldots,\alpha_k)$
writing ${\bm\oldb}(\bold x,\bm{\alpha})$ and $\oD_{{\bm\oldb},\bm{\alpha}}$ in place of ${\bm\oldb}(\bold x)$ and $\oD_{\bm\oldb}$, respectively. On one hand,
if $(\bold{Sub})$ holds, then there exists $\bold x\in\oD$ such that ${\bm\oldb}(\bold x,\bm{\alpha})<\bold 0$. So there exists $\varepsilon>0$ such that
${\bm\oldb}(\bold x,\bm{\alpha}+\bm{\varepsilon})={\bm\oldb}(\bold x,\bm{\alpha})+\bm{\varepsilon}<\bold 0$ (component-wise), where
$\bm{\varepsilon}:=(\varepsilon,\ldots,\varepsilon)$ and we used the linearity of ${\bm\oldb}(\bold x,\cdot)$. Therefore, the SBM with parameter $\bm{\alpha}+\bm{\varepsilon}$
satisfies $(\bold{Sub})$ and so, by Lemma \ref{le:subsub'}, it satisfies $(\bold{Sub})'$. On the other hand, the linearity of ${\bm\oldb}(\bold x,\cdot)$ ensures
$\oD_{{\bm\oldb},\bm{\alpha}}\equiv\oD_{{\bm\oldb},\bm{\alpha}+\bm{\varepsilon}}$. Therefore, by $({\mathcal Crit})$
\begin{equation}\label{contr2}
\min_{\bold y\in\oD_{{\bm\oldb},\bm{\alpha}+\bm{\varepsilon}}}\oldb_1(\bold y,\bm{\alpha}+\bm{\varepsilon})
=\min_{\bold y\in\oD_{{\bm\oldb},\bm{\alpha}}}\oldb_1(\bold{y},\bm{\alpha}+\bm{\varepsilon})=\bm{\varepsilon}+\min_{\bold{y}\in\oD_{{\bm\oldb},\bm{\alpha}}}\oldb_1(\bold{y},\bm{\alpha})=\bm{\varepsilon}>\bold 0.
\end{equation}
Since $({\mathcal Sup})$ implies $(\bold{Sup})$ (see proof of part $(ii)$ of this theorem), by Lemma \ref{le:Sup} and \eqref{contr2} we have that
the SBM with parameter $\bm{\alpha}+\bm{\varepsilon}$ satisfies $(\bold{Sup})'$, yielding a contradiction.
\\
\noindent$\square$

\noindent {\it Proof\,\,of\,\,Lemma\,\,\ref{le:Sup}.} We start proving $(\bold{Sup})\Rightarrow (\bold{Sup})'$.
We shall check later on that $\ote<\infty$, $\bold{x}(\ote)\in\ocalD$ and ${\bm\oldb}(\bold{x}(\ote))>\bold 0$.
Then, consider the curve with trace $\mathcal C:=\mathcal S\cup\mathcal{C}_\oE$.
It is easily realized that it is continuous, non-decreasing and $\mathcal C$ satisfies
\eqref{eq:intersect_k}. Curve $\mathcal C$ satisfies also \eqref{sup-cond_k}. Indeed, since $\oldb_i(\bold x)>0$
for any $i$, $\bold x\in\mathcal C \cap\overset\circ{\oE}$ and
${\bm\oldb}(\bold{x}(\ote))>\bold 0$ (component-wise), we have
\[
\min_{\bold x\in\mathcal{C}_{\oE}}\min_{1\leq i\leq k}\oldb_i(\bold x)>0.
\]
\noindent${\it Proof\,\,of\,\, \ote<\infty.}$ Since by assumption $\bm{\rho}\neq\bold 0$ over compact $\oD$, we have
\[
m_{\oD}^{({\bm\oldb})}:=\min_{\bold{x}\in \oD}\sum_{i=1}^{k}|\oldb_i(\bold{x})|\in\mathbb{R}_+
\]
and so
\begin{align}
\infty>\sum_{i=1}^{k}x_i(\ote)&=\sum_{i=1}^{k}\left ((\bold{x}_0)_i+\lim_{y_0<y\uparrow \ote}\int_{y_0}^{y}x_i'(s)\,\mathrm d s\right)\nonumber\\
&=\sum_{i=1}^{k}(\bold{x}_0)_i+\lim_{y_0<y\uparrow \ote}\int_{y_0}^{y}\left(\sum_{i=1}^{k}x_i'(s)\right)\,\mathrm d s\nonumber\\
&>\lim_{y_0<y\uparrow \ote}\int_{y_0}^{y}\left(\sum_{i=1}^{k}\oldb_i(\bold{x}(s))\right)\,\mathrm d s\nonumber\\
&\geq m_{\oD}^{({\bm\oldb})}\lim_{y_0<y\uparrow \ote}\int_{y_0}^{y}\,\mathrm d s =m_{\oD}^{({\bm\oldb})}(\ote-y_0).\label{eq:mDfin}
\end{align}
Therefore $\ote<\infty$.\\
\noindent${Proof\,\,of\,\,\bold{x}(\ote)\in\ocalD\,\, and \,\,{\bm\oldb}(\bold{x}(\ote))>\bold 0}.$
We have either $\bold{x}(\ote)\in\partial\oE\cap\ocalD^c$ or
$\bold{x}(\ote)\in\partial\oE\cap\ocalD$, where $\ocalD^c$ stands for the complement of $\ocalD$.
Reasoning by contradiction, suppose $\bold{x}(\ote)\in\partial\oE\cap\ocalD^c$. Define
\[
\mathcal{I}:=\{i:\,\,\oldb_i(\bold{x}(\ote))=0\}\quad\text{and}\quad\mathcal{J}:=\{j:\,\,\oldb_j(\bold{x}(\ote))>0\},
\]
and note that $\mathcal{I}\neq\emptyset$
(since in particular $\bold{x}(\ote)\in\partial\oE\setminus\bigcup_{i=1}^{k}\{\bold{x}:\,\,x_i=0\}$, and so there exists
$i\in\{1,\ldots,k\}$ such that $\oldb_i(\bold{x}(\ote))=0$) and $\mathcal{J}\neq\emptyset$ (since $\bold{x}(\ote)\in\oE\setminus\mathcal{Z}$).
By the irreducibility of $\bm{\chi}$
\begin{equation}\label{eq:iojo}
\text{There exist $i_0\in\mathcal I$ and $j_0\in\mathcal{J}$ such that $\frac{\partial \oldb_{i_0}}{\partial x_{j_0}}(\bold{x}(\ote))>0$.}
\end{equation}
Since $\oldb_{i_0}(\bold{x}(y))>0$ for any $y\in (y_0,\ote)$ and $\oldb_{i_0}(\bold{x}(\ote))=0$,
we have
\begin{equation}\label{deriv1}
\frac{\diff \oldb_{i_0}(\bold{x}(y))}{\diff y}|_{y=\ote}\le 0.
\end{equation}
On the other hand, for any $y\in (y_0,y_+)$,
\begin{align*}
\frac{\diff \oldb_{i_0}(\bold{x}(y))}{\diff y }&=\sum_{j=1}^{k}\frac{\partial \oldb_{i_0}}{\partial x_j}(\bold{x}(y))x_j'(y)\\
&=\sum_{j=1}^{k}\frac{\partial \oldb_{i_0}}{\partial x_j}(\bold{x}(y))\oldb_j(\bold{x}(y))\\
&=\sum_{j\in\mathcal{I}}\frac{\partial \oldb_{i_0}}{\partial x_{j}}(\bold{x}(y))\oldb_{j}(\bold{x}(y))+
\sum_{j\in\mathcal{J}}\frac{\partial \oldb_{i_0}}{\partial x_j}(\bold{x}(y))\oldb_j(\bold{x}(y)).
\end{align*}
Computing this quantity at $y=\ote$, by \eqref{eq:iojo}, the definition of the set of indexes $\mathcal I$ and $\mathcal J$ and the fact that all non-diagonal terms of
$J_{\bm\oldb}(\bold x)$ are non-negative, we have
\begin{equation}\label{eq:derpos}
\frac{\diff \oldb_{i_0}(\bold{x}(y))}{\diff y }|_{y=\ote}=
\sum_{j\in\mathcal J}\frac{\partial \oldb_{i_0}}{\partial x_j}(\bold{x}(\ote))\oldb_j(\bold{x}(\ote))\geq
\frac{\partial \oldb_{i_0}}{\partial x_{j_0}}(\bold{x}(\ote))\oldb_{j_0}(\bold{x}(\ote))>0,
\end{equation}
which contradicts \eqref{deriv1}. Therefore $\bold{x}(\ote)\in\ocalD$.
To check ${\bm\oldb}(\bold{x}(\ote))>\bold 0$, we start noticing that we clearly have ${\bm\oldb}(\bold{x}(\ote))\geq\bold 0$.
Reasoning by contradiction, assume that
$\oldb_{i_0}(\bold{x}(\ote))=0$ for some $i_0\in\{1,\ldots,k\}$. Then, arguing as above we get both \eqref{deriv1} and \eqref{eq:derpos},
and the claim follows. Finally we prove $(\bold{Sup})'\Rightarrow (\bold{Sup})$.
Reasoning by contradiction, suppose $\mathcal Z\neq\emptyset$ and consider set $\hat{\mathcal{F}}_{\bold z}$, $\bold{z}\in\mathcal Z$, defined by
\begin{equation}\label{eq:setF}
\hat{\mathcal{F}}_{\bold z}:=\cup_{j=1}^{k}\hat{\mathcal{F}}^{(j)}_{\bold z},
\end{equation}
where
\[
\hat{\mathcal{F}}^{(j)}_{\bold z}:=\{\bold w\in\oD:\,\,w_j=z_j \text{ and }  w_h\le z_h,\;\;   \forall h\neq j\}.
\]
Since $\bold z$ is a zero of $\oldb_i$ for any $i$, we have $\bold z\neq\bold 0$ with all components $z_i\neq 0$.
By the monotonicity properties of the
$\oldb_i$'s, for any $j\in\{1,\ldots,k\}$ and any $\bold x\in\hat{\mathcal{F}}^{(j)}_\bold z$, it holds $\oldb_j(\bold{x})\leq \oldb_j(\bold z)=0$, and so
\[
\min_{1\leq i\leq k} \oldb_i(\bold{x})\leq 0,\quad\text{$\forall$ $\bold{x}\in\hat{\mathcal{F}}_{\bold z}$.}
\]
By the continuity and the monotonicity of  curve $\mathcal C$ and condition \eqref{eq:intersect_k}, we necessarily have
\[
\mathcal{C}\cap\hat{\mathcal{F}}_{\bold z}=(\mathcal{S}\cup \mathcal{C}_{\oE})\cap\hat{\mathcal{F}}_{\bold z}\neq\emptyset.
\]
As immediate consequence of Lemma \ref{le:initialpoint}, we have $\mathcal{S} \cap\hat{\mathcal{F}}_{\bold z}=\emptyset$, and so
$\mathcal{C}_{\oE}\cap\hat{\mathcal{F}}_{\bold z}\neq\emptyset$. Taking $\bold w\in \mathcal{C}_{\oE}\cap \hat{\mathcal{F}}_{\bold z}$,  we have
\[
\min_{\bold{x}\in\mathcal{C}_{\oE}}\min_{1\leq i\leq k}\oldb_{i}(\bold{x})\leq\min_{1\leq i\leq k}\oldb_i(\bold{w})\leq 0,
\]
which contradicts \eqref{sup-cond_k}.
\\
\noindent$\square$

\noindent {\it Proof\,\,of\,\,Lemma\,\,\ref{le:esistecurva}.} Firstly, we show
\[
\bold{x}(\ote)\in\mathcal Z\cap\partial\oE\quad\text{if and only if}\quad \ote=+\infty,
\]
secondly, we prove $\bold{x}(\ote)\in\mathcal{Z}\cap\partial\oE$. The implication
$\bold{x}(\ote)\in\mathcal{Z}\cap\partial\oE\Rightarrow \ote=+\infty$
follows by the uniqueness of the solution of Cauchy problem \eqref{eq:cauchy}.We now prove
\noindent $\ote=+\infty\Rightarrow\bold{x}(\ote)\in\mathcal{Z}\cap\partial\oE$.
We reason by contradiction. For $\varepsilon>0$ arbitrarily small, let $\oD^{(\varepsilon)}$ denote the closed subset
of $\oD$ obtained by erasing from $\oD$ the open balls with radius $\varepsilon>0$ centered at zeros of $b$ in $\oD$.
We have
\[
m_{\oD^{(\varepsilon)}}^{({\bm\oldb})}:=\min_{\bold{x}\in\oD^{(\varepsilon)}}\sum_{i=1}^{k}|\oldb_i(\bold{x})|\in\mathbb{R}_+
\]
and, since $\bold{x}(\ote)\notin\mathcal Z\cap\partial\oE$,
curve $\zeta_{\oE}(y):=\bold{x}(y)$,
$y_0\leq y\leq \ote$,
is contained in $\oD^{(\varepsilon)}$. So arguing as for \eqref{eq:mDfin}, we have $m_{\oD^{(\varepsilon)}}^{({\bm\oldb})}(\ote-y_0)<\infty$,
which proves the finiteness of $\ote$ and therefore yields a contradiction. We finally prove
$\bold{x}(\ote)\in\mathcal Z\cap \partial \oE$.
Since $\bold{x}(\ote)\in\partial\oE$ 
we have $\bold{x}(\ote)\in (\mathcal{Z}\cap\partial\oE)\cup(\mathcal Z^c\cap \partial\oE)$.
We shall prove that $\bold{x}(\ote)\notin\mathcal Z^c\cap\partial\oE$. Note that
\[
\mathcal Z^c\cap\partial\oE=(\partial\oE\cap\ocalD\cap\mathcal Z^c)\cup(\partial\oE\cap\ocalD^c\cap\mathcal Z^c).
\]
We shall show that $\bold{x}(\ote)\notin\partial\oE\cap\ocalD$ and
$\bold{x}(\ote)\notin\partial\oE\cap\ocalD^c\cap\mathcal{Z}^c$.
Reasoning by contradiction, assume $\bold{x}(\ote)\in\partial\oE\cap\ocalD$. For an arbitrarily fixed $\bold z\in\mathcal Z$, there exists
$\bold{x}\in\hat{\mathcal{F}}_{\bold z}\cap\mathcal{C}{_\oE}\cap\overset\circ{\oD}$, where $\hat{\mathcal F}_\bold z$ is defined by \eqref{eq:setF}.
So by the monotonicity properties of the
$\oldb_i$'s we have $\oldb_j(\bold{x})\leq \oldb_j(\bold z)= 0$, for some $j$.
Therefore $\bold{x}\notin\oE$, which is impossible.
Reasoning again by contradiction, assume $\bold{x}(\ote)\in\partial\oE\cap\ocalD^c\cap\mathcal{Z}^c$.
Arguing as in the proof of Lemma \ref{le:Sup}
we deduce both \eqref{deriv1} and \eqref{eq:derpos}, getting in this way a contradiction.
\\
\noindent$\square$

\noindent {\it Proof\,\,of\,\,Lemma\,\,\ref{le:subzeri}.} Let $\bold x\in\overset\circ{\oD}$
be such that ${\bm\oldb}(\bold x)<\bold 0$ (component-wise) and define
\begin{equation}\label{eq:setG}
\hat{\mathcal{G}}_\bold{x}:=\{\bold w\in\oD:\,\,w_j\leq x_j\,\,\forall j=1,\ldots,k\}.
\end{equation}
It is easily seen that
\begin{equation}\label{eq:GxFx}
\text{$\oE\subset\hat{\mathcal{G}}_\bold x$\quad and\quad $\hat{\mathcal{G}}_\bold{x}\cap\ocalD=\emptyset$.}
\end{equation}
Therefore, condition $(\bold{Sub})$ implies $\oE\cap\ocalD=\emptyset$. On the other hand, we have that condition $(\bold{Sup})'$ implies
$\oE\cap\ocalD\neq\emptyset$ (indeed $\bold{x}(\ote)\in\oE\cap\ocalD$). Consequently, $(\bold{Sub})$ implies that $(\bold{Sup})'$ does not hold,
and so, by Lemma \ref{le:Sup} condition $(\bold{Sub})$ implies $\mathcal Z\neq\emptyset$.
\\
\noindent$\square$

\noindent {\it Proof\,\,of\,\,Lemma\,\,\ref{le:subsub'}.} Since the implication
$({\mathcal Sub})\Rightarrow(\bold{Sub})$ is obvious, we only need to prove $(\bold{Sub})\Rightarrow (\bold{Sub}')$.
By Lemma \ref{le:subzeri} we have $\mathcal{Z}\neq\emptyset$ and so by Lemma \ref{le:esistecurva} it holds
$\ote=+\infty$ and $\bold{x}(\ote)\in\mathcal{Z}$. Let $\bold x\in\overset\circ{\oD}$
be such that ${\bm\oldb}(\bold x)<\bold 0$ and consider sets $\hat{\mathcal{F}}_{\bold x}$
and $\hat{\mathcal{G}}_{\bold x}$ defined by \eqref{eq:setF} (with $\bold x$ in place of $\bold z$)
and \eqref{eq:setG}. By \eqref{eq:GxFx} we have
$\bold{x}(\ote)\in\overset\circ{\oD}$, and so $\bold{x}(\ote)\in\mathcal Z\cap\overset\circ{\oD}$.
Note that the curve with trace $\mathcal C:=\mathcal{S}\cup\mathcal{C}_{\oE}\cup\tilde{\mathcal S}$
is continuous, non-decreasing, satisfies \eqref{eq:intersect_k} and it is such that $\mathcal{C}\cap\hat{\mathcal{F}}_{\bold x}\neq\emptyset$.
Therefore there exists an $i$  such that $\mathcal{C} \cap \hat{\mathcal{F}}_{\bold x}^{(i)} \neq \emptyset$. Let $\bold{x}_{*}\in \mathcal{C}\cap  \hat{\mathcal{F}}^{(i)}_{\bold x}$.
As an immediate consequence of Lemma \ref{le:initialpoint}, we have $\mathcal{S}\cap\hat{\mathcal{F}}_{\bold x}=\emptyset$, and therefore
$\bold{x}_{*}\in \mathcal{C}_{\oE}\cup\tilde{\mathcal{S}}\cap \hat{\mathcal{F}}_{\bold x}^{(i)}$.
Finally, since $\oldb_i$ is non-decreasing with respect to the $j$th variable, $j\neq i$, we have
\[
\min_{1\leq j\leq k}\oldb_j(\bold x_*)\leq \oldb_i(\bold x_*)\leq\sup_{\bold{w}\in\hat{\mathcal{F}}_{\bold x}^{(i)}} \oldb_i(\bold w)=\oldb_i(\bold x)<0,
\]
which yields \eqref{sub-cond_k} and completes the proof.
\\
\noindent$\square$

\noindent {\it Proof\,\,of\,\,Lemma\,\,\ref{le:critcrit'}.}
By Lemma \ref{le:esistecurva} we have $\ote=+\infty$ and $\bold{x}(\ote)\in\mathcal Z$.
Similarly to the proof of Lemma \ref{le:subsub'} it is checked that
$\bold{x}(\ote)\in\overset\circ{\oD}$. The curve with trace $\mathcal C:=\mathcal{S}\cup\mathcal{C}_{\oE}\cup\tilde{\mathcal{S}}$ is clearly continuous, non-decreasing and satisfies
\eqref{eq:intersect_k}. Moreover, it satisfies \eqref{crit_k} since
\[
\min_{\bold x\in\mathcal C}\oldb_i(\bold x)=\oldb_i(\bold{x}(\ote))=0,\quad\text{$\forall$ $i=1,\ldots,k$.}
\]
Indeed, each $\oldb_i$ is non-negative on $\mathcal S\cup\mathcal C_{\oE}$, strictly positive
on $\tilde{\mathcal{S}}\setminus\{\bold{x}(\ote)\}$ and $\oldb_i(\bold{x}(\ote))=0$, for any $i=1,\ldots,k$. In particular, the strict positivity
of  $\oldb_i$'s on $\tilde{\mathcal{S}}\setminus\{\bold{x}(\ote)\}$ follows by the fact that $\bold{x}(\ote)$ is the unique point of minimum of
each $\oldb_i$ on $\tilde{\mathcal S}$ and $\oldb_i(\bold{x}(\ote))=0$ (indeed, this is a consequence of the convexity of $\tilde{\mathcal{S}}$, the convexity
of the $\oldb_i$'s and the fact that $\bold{x}(\ote)$ is the unique critical point of  $\oldb_i$'s on $\tilde{\mathcal S}$ since $\lambda_{PF}(\bold{x}(\ote))=0$).
\\
\noindent$\square$

\noindent {\it Proof\,\,of\,\,Lemma\,\,\ref{le:subequiv}.} We start noticing that the implication $(\bold{Sub})\Rightarrow (\bold{Sub})'''$ is obvious. We now prove
$(\bold{Sub})'''\Rightarrow (\bold{Sub}).$ Let $\bold x\in\oD$ be such that ${\bm\oldb}(\bold x)\leq\bold 0$ and ${\bm\oldb}(\bold x)\neq\bold 0$. We have
\[
\mathcal{N}_0:=\{i\in\{1,\ldots,k\}:\,\,\oldb_i(\bold{x})<0\}\neq\emptyset.
\]
Let $\mathcal{N}_0^c:=\{1,\ldots,k\}\setminus\mathcal{N}_0=\{i\in\{1,\ldots,k\}:\,\,\oldb_i(\bold{x})=0\}$ be the complement of $\mathcal{N}_0$
and define
$x^0:=\min_{i\in \mathcal{N}_0}x_i$. Note that  $x^0>0$ (indeed if there exists $j\in\mathcal{N}_0$ such that $x_j=0$ then $\oldb_j(\bold x)\ge 0$, which is impossible).
Define
\[
\mathcal{N}_1:=\Biggl\{i\in\mathcal{N}_0^c: \exists h\in\mathcal{N}_0\,\,\text{with}\,\,\frac{\partial \oldb_i}{\partial x_h}(\bold x)>0\Biggr\}.
\]
The irreducibility of $\bm{\chi}$ implies $\mathcal{N}_1\neq\emptyset$. Let $\bold{I}^{(0)}$ be the vector with components $I_i^{(0)}:=\bold{1}(i\in\mathcal{N}_0)$, $i=1,\ldots,k$,
and put $\bold{x}_1:=\bold{x}-\varepsilon\frac{x^0}{2}\bold{I}^{(0)}$ for some $\varepsilon\in (0,1]$ sufficiently small.
By construction $\bold{x}_1\in\oD$ and $\oldb_j(\bold{x}_1)<0$ for any $j\in\mathcal{N}_0\cup\mathcal{N}_1$.
Iterating the procedure we can find $(i)$ A finite family of sets $\{\mathcal{N}_{h}\}_{1\leq h\leq\ell}$ with $1\leq\ell\leq k$, $\mathcal{N}_h\neq\emptyset$, for any $1\leq h\leq\ell$,
and $\bigcup_{h=1}^{\ell}\mathcal{N}_h=\{1,\ldots,k\}$ $(ii)$ A point $\bold{x}_\ell\in\oD$ such that $\oldb_{j}(\bold{x}_\ell)<\bold 0$ for any
$j\in\bigcup_{h=1}^{\ell}\mathcal{N}_h$. The claim is proved. We now show $(\bold{Sub})'' \Rightarrow (\bold{Sub})$. For this, we shall prove  implications:
$(\bold{Sub})(i)\Rightarrow (\bold{Sub})$, $(\bold{Sub})(ii)\Rightarrow (\bold{Sub})$ and
$(\bold{Sub})(iii)\Rightarrow (\bold{Sub})$. We start proving $(\bold{Sub})(i)\Rightarrow (\bold{Sub})$.
Let $\bold{z}^{(1)},\bold{z}^{(2)}\in\oD$ be two distinct zeros of $\bm{\rho}$. Since $\oD$ is convex we have $\bold{x}_{\theta}:=\theta\bold{z}^{(1)}+(1-\theta)\bold{z}^{(2)}\in\oD$,
for any $\theta\in [0,1]$. By the convexity of $\oldb_i$ we have $\oldb_i(\bold x_\theta)\leq 0$ for any $i=1,\ldots,k$ and $\theta\in [0,1]$.
Since $\bold{z}^{(1)}\neq\bold{z}^{(2)}$ there exists $i_0\in\{1,\ldots,k\}$ such that $z_{i_0}^{(1)}\neq z_{i_0}^{(2)}$. Consider function
$\varphi_{i_0}(\theta):=\oldb_{i_0}(\bold{x}_\theta)$. A straightforward computation shows $\varphi_{i_0}''(\theta)>0$ for any $\theta\in (0,1)$. Therefore,
\[
\oldb_{i_0}(\bold{x}_\theta)<\theta\oldb_{i_0}(\bold{z}^{(1)})+(1-\theta)\oldb_{i_0}(\bold{z}^{(2)})=0,\quad\text{$\forall$ $\theta\in (0,1)$.}
\]
This implies $(\bold{Sub})'''$ and so (by one of the previous steps) $(\bold{Sub})$. We continue by proving $(\bold{Sub})(ii)\Rightarrow (\bold{Sub})$.
Since $\mathrm{det}J_{{\bm\oldb}}(\bold{z})\neq 0$, the linear map $J_{{\bm\oldb}}(\bold{z}):\mathbb R^k\to\mathbb R^k$ is surjective and therefore
there exists a vector $\bold{w}\in\mathbb{R}^k$
such that $J_{{\bm\oldb}}(\bold{z})\bold{w}<\bold 0$ (component-wise). Set $\bold{x}_\delta:=\bold{z}+\delta\bold{w}$, with $\delta>0$ small enough so that $\bold x_\delta\in\oD$.
Let $j\in\{1,\ldots,k\}$ be arbitrarily fixed. The directional derivative of $\oldb_j$ along vector $\bold w$ computed at $\bold x_\delta$ is equal to
\[
\frac{\partial \oldb_j}{\partial\bold w}(\bold x_\delta)=\sum_{i=1}^{k}\frac{\partial \oldb_j}{\partial x_i}(\bold z+\delta\bold w)w_i.
\]
For $\delta>0$ small enough, this quantity is strictly less than zero since $J_{{\bm\oldb}}(\bold{z})\bold{w}<\bold 0$.
Therefore, for $\delta>0$ small enough,
$\oldb_j(\bold{x}_\delta)<\oldb_j(\bold{z})=0$, and the proof of the claim is completed.
We now prove $(\bold{Sub})(iii)\Rightarrow (\bold{Sub})$. If $\lambda_{PF}(\bold z)>0$, then $J_{{\bm\oldb}}(\bold{z})\bold{\ov}_{PF}(\bold z)=\lambda_{PF}(\bold z)\bold{\ov}_{PF}(\bold z)>\bold 0$, where the latter inequality
holds (component-wise) due the positivity of  eigenvector $\bold{\ov}_{PF}(\bold z)$. Since $\bold z\in\overset\circ{\oD}$, $\bold{x}_\delta:=\bold{z}+\delta(-\bold{\ov}_{PF}(z))\in
\overset\circ{\oD}$ for a sufficiently small $\delta>0$. Let $j\in\{1,\ldots,k\}$ be arbitrarily fixed.
Since $J_{{\bm\oldb}}(\bold{z})\bold{\ov}_{PF}(\bold z)>\bold 0$, for $\delta>0$ small enough,
the directional derivative of $\oldb_j$ along vector $-\bold{\ov}_{PF}(\bold z)$ computed at $\bold{x}_\delta$ is strictly less than zero, indeed
\[
\frac{\partial \oldb_j}{\partial(-\bold{\ov}_{PF}(\bold z))}(\bold x_\delta)=-\sum_{i=1}^{k}\frac{\partial \oldb_j}{\partial x_i}(\bold z+\delta(-\bold{\ov}_{PF}(\bold z)))(\bold{\ov}_{PF}(\bold z))_{i}.
\]
Therefore, for $\delta>0$ small enough, $\oldb_j(\bold{x}_\delta)<\oldb_j(\bold z)=0$, and the claim is proved.
If $\lambda_{PF}(\bold z)<0$, by a similar argument one has that, for $\delta>0$ small enough,
$\oldb_j(\bold{x}_\delta)<0$ for any $j$, where $\bold{x}_\delta:=\bold{z}+\delta\bold{\ov}_{PF}(z)\in
\overset\circ{\oD}$, which concludes the proof of the claim. Finally, we prove $(\bold{Sub})\Rightarrow (\bold{Sub})''$.
By Lemma \ref{le:subzeri} we have that $(\bold{Sub})$ implies either $(\bold{Sub})''$ or $(\bold{I})$ or $(\bold{Crit})$.
By Lemma \ref{le:subsub'} condition $(\bold{Sub})$ implies $(\bold{Sub})'$ which, in turn, guarantees that
neither $(\bold I)$ nor $(\bold{Crit})'$ holds. The claim follows by Lemma \ref{le:critcrit'}.
\\
\noindent$\square$

\noindent {\it Proof\,\,of\,\,Lemma\,\,\ref{le:Ithensub}.}
By assumption function $\bm{\rho}$ has a unique zero $\bold z$ which lies in $\ocalD$.
We preliminary note that: $1)$ at least one term $(J_{{\bm\oldb}}(\bold z))_{ii}$, $i\in\{1,\ldots,k\}$,
on the diagonal of $J_{{\bm\oldb}}(\bold z)$ is equal to zero and all non-zero terms on the diagonal of $J_{{\bm\oldb}}(\bold z)$ are strictly negative;
$2)$ all terms $(J_{{\bm\oldb}}(\bold z))_{ij}$, $i,j\in\{1,\ldots,k\}$, $i\neq j$, outside the diagonal of $J_{{\bm\oldb}}(\bold z)$ are non-negative.
Without loss of generality we assume $(J_{{\bm\oldb}}(\bold z))_{11}=0$. For ease of notation, for $i=2,\ldots,k$, we set
\[
\bold{1}_{1i}:=\bold{1}((J_{\bm\oldb}(\bold z))_{1i}\in\mathbb{R}_+).
\]
Let $\bold{v}^{(\delta)}:=(-1,-\delta\bold{1}_{12},\ldots,-\delta\bold{1}_{1k})\in\mathbb R^k$, for some $\delta>0$.
For any $i=1,\ldots,k$, we have
\[
(J_{{\bm\oldb}}(\bold z)\bold{v}^{(\delta)})_i=\sum_{h=1}^{k}(J_{{\bm\oldb}}(\bold z))_{ih}v_h^{(\delta)}
=(J_{{\bm\oldb}}(\bold z))_{ii}v_i^{(\delta)}+\sum_{h\neq i}^{1,k}(J_{{\bm\oldb}}(\bold z))_{ih}v_h^{(\delta)}.
\]
In particular,
\begin{equation}\label{eq:jb1neg}
(J_{{\bm\oldb}}(\bold z)\bold{v}^{(\delta)})_1=-\delta\sum_{h=2}^{k}(J_{{\bm\oldb}}(\bold z))_{1h}\bold{1}_{1h}<0,
\end{equation}
where the latter inequality follows noticing that by the irreducibility of matrix $\bm{\chi}$ we have $\bold{1}_{1h}=1$ for some $h\in\{2,\ldots,k\}$.
Define $\bold{x}_\varepsilon^{(\delta)}:=\bold z+\varepsilon\bold v^{(\delta)}$, $\varepsilon>0$. Since $\bold{z}\in\ocalD$, for $\varepsilon>0$
small enough, we have $\bold x_\varepsilon^{(\delta)}\in\overset\circ{\oD}$. Due to \eqref{eq:jb1neg},
arguing as in the proof of Lemma \ref{le:subequiv}
we deduce that the directional derivative of $\oldb_1$, along vector $\bold v^{(\delta)}$ and computed at $\bold x_\varepsilon^{(\delta)}$, is strictly negative
and so
\begin{equation}\label{eq:b1neg}
\oldb_1(\bold x_{\varepsilon}^{(\delta)})<\oldb_1(\bold z)=0.
\end{equation}
Now, take $i_0\neq 1$. We have that for sufficiently small $\delta$:
\[
(J_{{\bm\oldb}}(\bold z)\bold{v}^{(\delta)})_{i_0}=-\delta(J_{{\bm\oldb}}(\bold z))_{i_0i_0}\bold{1}_{1i_0}-(J_{{\bm\oldb}}(\bold z))_{i_01}-\delta\sum_{h\neq i_0}^{2,k}(J_{{\bm\oldb}}(\bold z))_{i_0 h}\bold{1}_{1h}\leq 0.
\]
If $(J_{{\bm\oldb}}(\bold z)\bold{v}^{(\delta)})_{i_0}<0$, then studying as usual the directional derivative along $\bold v^{(\delta)}$ we deduce
\begin{equation}\label{eq:bi0neg1}
\text{$\oldb_{i_0}(\bold{x}_\varepsilon^{(\delta)})<0$ (for $\varepsilon$ and $\delta$ small enough).}
\end{equation}
If $(J_{{\bm\oldb}}(\bold z)\bold{v}^{(\delta)})_{i_0}=0$, then
$\chi_{i_01}=0$ and $\chi_{i_0 h}=0$ for any $h\in\{2,\ldots,k\}\setminus\{i_0\}$ such that $\chi_{1h}>0$ (indeed, for $i\neq j$,
$(J_{\bm\oldb}(\bold z))_{ij}=0$ if and only if $\chi_{ij}=0$). Therefore (for $\varepsilon$ and $\delta$ small enough)
\begin{align}
\oldb_{i_0}(\bold{x}_{\varepsilon}^{(\delta)})&=\alpha_{i_0}-z_{i_0}+\varepsilon\delta\bold{1}_{1i_0}
+r^{-1}(1-r^{-1})^{r-1}\left(\sum_{h\in\{2,\ldots,k\}\setminus\{i_0\}:\chi_{1h}=0}\chi_{i_0 h}(z_h-\varepsilon\delta\bold{1}_{1h})\right)^r\nonumber\\
&=\oldb_{i_0}(\bold z)=0.\label{eq:bi00}
\end{align}
Collecting \eqref{eq:b1neg}, \eqref{eq:bi0neg1} and \eqref{eq:bi00} we get $(\bold{Sub})$.
\\
\noindent$\square$

\subsection{Proof of Proposition \ref{prop:equiv}}\label{subsec:prop21}

In this proof we simply write $\mathcal A$, $\mathcal U$ and $T$ in place of $\mathcal A_n$, $\mathcal U_n$ and $T_n$.
We first prove
\begin{equation}\label{eq:firstinclusion}
\cG\subseteq\mathcal{A}(T).
\end{equation}
This is equivalent to prove $\bigcup_{h=0}^{H}\cG_h\subseteq\mathcal{A}(T)$, for any $H\in\N\cup\{0\}$.
We show this claim by induction over $H$. We clearly have $\cG_0=\mathcal{A}(0)\subseteq\mathcal{A}(T)$. Assume
\[
\bigcup_{h=0}^H\cG_h\subseteq \cA(T),\quad\text{for some $H\in\N$.}
\]
The inclusion \eqref{eq:firstinclusion} follows if we check $\cG_{H+1}\subseteq\mathcal{A}(T)$.
Take $v\in\cG_{H+1}$ and, reasoning by contradiction, suppose $v\notin\mathcal{A}(T)$.
By the definition of the bootstrap percolation process $v\in\cG_{H+1}$ has at least $r$ neighbors in
$\bigcup_{h=0}^H\cG_h$. This set of active nodes is contained in $\mathcal{U}(T)$ due to the inductive hypothesis and relation $\mathcal{A}(T)=\mathcal{U}(T)$.
Consequently, by \eqref{eq:Mv},  $M_v(T)\geq r$ and so, by \eqref{eq:AMv}, $v\in\mathcal{A}(T)$, which is a contradiction.
We now prove
\begin{equation} \nonumber 
\cG\supseteq\mathcal{A}(T).
\end{equation}
For this it suffices to prove that $\mathcal{A}(t)\subseteq\cG$ for any $0\leq t\leq T$. We denote by $\Delta\mathcal{A}(t)$ the set of nodes that become active exactly at time $t$.
Reasoning by contradiction, assume that there exists at least a $v\in\Delta\mathcal{A}(t)$, $r\leq t\leq T$, such that $v\notin\cG_h$, for any $h\in\N\cup\{0\}$.
Then there must exist a minimum time $t_0$ with $r\leq t_0\leq T$ such that $\Delta \mathcal{A}(t_0)\not\subseteq\cG$.
Since $v\in\Delta\cA(t_0)$, it has at least $r$ neighbors in $\cU(t_0)$. By construction we have
$\cU(t_0)\subseteq\cA(t_0-1)$ and $\cA(t_0-1)\subseteq\cG$. So $v$ has $r$ neighbors in $\cG$. Therefore
$v\in\cG$, which is a contradiction.

\subsection{Proofs of Lemmas \ref{le:bas2}, \ref{le:June15}, \ref{le:aspi}, \ref{le:initialpoint}}\label{subsec:lemmas}

Since the proofs of Lemmas \ref{le:bas2} and \ref{le:June15} exploit Lemma \ref{le:aspi}, the proof of Lemma \ref{le:initialpoint} is based on Lemma \ref{le:bas2} and the proof
of Lemma \ref{le:aspi} is self-contained, we first prove Lemma \ref{le:aspi}, then we prove
Lemmas \ref{le:bas2} and \ref{le:June15}, and finally we prove Lemma \ref{le:initialpoint}.
\\
\noindent{\it Proof\,\,of\,\,Lemma\,\,\ref{le:aspi}.} Note that
\[
\oldpi(\lfloor\bold{x}\nan\rfloor,\oq{{ij}})=P\left(\sum_{j\in\{1,\ldots,k\}:\,\,x_j>0}\mathrm{Bin}(\lfloor x_j \nac{j}{n}\rfloor,\oq{{ij}})\geq r\right).
\]
Without loss of generality we denote by $x_1,\ldots,x_h$, $1\leq h\leq k$, the strictly positive $x$'s.
We proceed by induction on $h\in\mathbb N$. By the third relation in
\eqref{eq:ptcto0bm} and formula (8.1) in \cite{JLTV}, for any $i\in\{1,\ldots,k\}$, $j\in\{1,\ldots,h\}$ and $m\in\mathbb N$,
\begin{align}
P\left(\mathrm{Bin}(\lfloor x_j \nac{j}{n}\rfloor,\oq{{ij}})\geq m\right)
=\frac{(\lfloor x_j \nac{j}{n}\rfloor \oq{{ij}})^m}{m!}(1+O(\lfloor x_j \nac{j}{n}\rfloor \oq{{ij}}+(\lfloor x_j \nac{j}{n}\rfloor)^{-1})).\label{eq:apprbin}
\end{align}
The case $h=1$ is an immediate consequence of \eqref{eq:apprbin}.
Now, assume $h\geq 2$ and suppose that the claim holds for $h-1$.
By the independence of binomial random variables and \eqref{eq:apprbin}, we have
\begin{align}
&\oldpi(\lfloor\bold{x}\nan\rfloor,\obq{i})=P\left(\sum_{j=1}^{h-1}\mathrm{Bin}(\lfloor x_j\nac{j}{n}\rfloor,\oq{{ij}})
+\mathrm{Bin}(\lfloor x_h \nac{h}{n}\rfloor,\oq{{ih}})\geq r\right)\nonumber\\
&\,\,\,\,\,\,
=\sum_{r_h=0}^{r}P\left(\sum_{j=1}^{h-1}\mathrm{Bin}(\lfloor x_j\nac{j}{n}\rfloor,\oq{{ij}})\geq r-r_h\right)P\left(\mathrm{Bin}(\lfloor x_h \nac{h}{n}\rfloor,\oq{{ih}})=r_h\right)\nonumber\\
&\,\,\,\,\,\,
+P\left(\mathrm{Bin}(\lfloor x_h \nac{h}{n}\rfloor,\oq{{ih}})\geq r+1\right).\label{eq:pi1}
\end{align}
By the inductive hypothesis we have
\begin{align}
&P\left(\sum_{j=1}^{h-1}\mathrm{Bin}(\lfloor x_j \nac{j}{n}\rfloor,\oq{{ij}})\geq m\right)\nonumber\\
&\qquad\qquad
=\left(1+O_{h-1}\right)\left(\sum_{j=1}^{h-1}\lfloor x_j \nac{j}{n}\rfloor \oq{{ij}}\right)^{m}/m!\nonumber\\
&\qquad\qquad
=\left(1+O_{h-1}\right)
\sum_{(r_1,\ldots,r_{h-1}):\,\,0\leq r_j\leq m,\,\,\sum_{j}r_j=m}\prod_{j=1}^{h-1}\frac{(\lfloor x_j \nac{j}{n}\rfloor \oq{{ij}})^{r_j}}{r_j!},\nonumber
\end{align}
where for ease of notation we set
\begin{equation}\label{eq:O}
O_{\ell}:=O\left(\sum_{j=1}^{\ell}(\lfloor x_j \nac{j}{n}\rfloor \oq{{ij}}+(\lfloor x_j \nac{j}{n}\rfloor)^{-1})\right),\quad 1\leq\ell\leq h.
\end{equation}
Combining this relation with \eqref{eq:apprbin} and \eqref{eq:pi1}, we have
\begin{align}
\oldpi(\lfloor\bold{x} \nan\rfloor,\obq{i})
&=(1+O(\lfloor x_h \nac{h}{n}\rfloor \oq{{ih}}+(\lfloor x_h \nac{h}{n}\rfloor)^{-1}))(1+O_{h-1})\nonumber\\
&\quad\quad\quad
\times
\sum_{r_h=0}^{r}\sum_{(r_1,\ldots,r_{h-1}):\,\,0\leq r_j\leq r-r_h,\,\,\sum_{j}r_j=r-r_h}
\prod_{j=1}^{h}\frac{(\lfloor x_j \nac{j}{n}\rfloor \oq{{ij}})^{r_j}}{r_j!}\nonumber\\
&+\frac{(\lfloor x_h \nac{h}{n}\rfloor \oq{{ih}})^{r+1}}{(r+1)!}(1+O(\lfloor x_h \nac{h}{n}\rfloor \oq{{ih}}+(\lfloor x_h \nac{h}{n}\rfloor)^{-1}))
\nonumber\\
&=(1+O_h)\Biggl(\sum_{r_h=0}^{r}\sum_{(r_1,\ldots,r_{h-1}):\,\,0\leq r_j\leq r-r_h,\,\,\sum_{j}r_j=r-r_h}
\prod_{j=1}^{h}\frac{(\lfloor x_j \nac{j}{n}\rfloor \oq{{ij}})^{r_j}}{r_j!}\nonumber\\
&\qquad\qquad\qquad\qquad\qquad\qquad
+\frac{(\lfloor x_h \nac{h}{n}\rfloor \oq{{ih}})^{r+1}}{(r+1)!}\Biggr)\nonumber\\
&=(1+O_h)\Biggl(\sum_{(r_1,\ldots,r_{h-1},r_h):\,\,0\leq r_j\leq r,\,\,\sum_{j}r_j=r}
\prod_{j=1}^{h}\frac{(\lfloor x_j \nac{j}{n}\rfloor \oq{{ij}})^{r_j}}{r_j!}\nonumber\\
&\qquad\qquad\qquad\qquad\qquad\qquad
+\frac{(\lfloor x_h \nac{h}{n}\rfloor \oq{{ih}})^{r+1}}{(r+1)!}\Biggr)\nonumber\\
&=(1+O_h)\Biggl(\left(\sum_{j=1}^{h}\lfloor x_j \nac{j}{n}\rfloor \oq{{ij}}\right)^{r}/r!
+\frac{(\lfloor x_h \nac{h}{n}\rfloor \oq{{ih}})^{r+1}}{(r+1)!}\Biggr)\nonumber\\
&=(1+O_h)\Biggl(\left(\sum_{j=1}^{h}\lfloor x_j \nac{j}{n}\rfloor \oq{{ij}}\right)^{r}/r!\Biggr),\nonumber
\end{align}
which concludes the proof.
\\
\noindent$\square$

\noindent{\it Proof\,\,of\,\,Lemma\,\,\ref{le:bas2}.} By the definition of ${\color{black}\oldB_n^{(i)}}$, we have
\begin{align}
\frac{{\color{black}\oldB_n^{(i)}}(\lfloor\bold{x}\nan\rfloor,\obq{{i}})}{\nac{i}{n}}&=\frac{1}{\nac{i}{n}}(\oan{{i}}+(\on{i}-\oan{{i}})\oldpi(\lfloor\bold{x}\nan\rfloor,\obq{{i}})
-\lfloor x_i \nac{i}{n}\rfloor)\nonumber\\
&=\left(\frac{\oan{i}}{\nac{i}{n}}-\frac{\lfloor x_i \nac{i}{n}\rfloor}{\nac{i}{n}}\right)+\frac{\on{i}-\oan{i}}{\nac{i}{n}}
\oldpi(\lfloor\bold{x}\nan\rfloor,\obq{{i}}).\label{eq:Brel}
\end{align}
By \eqref{eq:trivial} it follows
\[
\frac{\oan{i}}{\nac{i}{n}}-\frac{\lfloor x_i \nac{i}{n}\rfloor}{\nac{i}{n}}\to\alpha_i-x_i.
\]
By the second relation in \eqref{eq:ptcto0bm}, Lemma \ref{le:aspi} and the definition of $\nac{i}{n}$,
we have
\begin{align}
\frac{\on{i}-\oan{i}}{\nac{i}{n}}
\oldpi(\lfloor\bold{x}\nan\rfloor,\obq{{i}})&\sim_e\frac{\on{i}}{\nac{i}{n}}\left(\sum_{j=1}^{k}x_j \nac{j}{n}\oq{{ij}}\right)^r/r!\nonumber\\
&=\frac{\on{i}(\nac{i}{n}\op{i})^r}{(r!)\nac{i}{n}}\left(\sum_{j=1}^{k}x_j\frac{\nac{j}{n}\oq{{ij}}}{\nac{i}{n}\op{i}}\right)^r\nonumber\\
&=\frac{\on{i}(\op{i})^{r}(\nac{i}{n})^{r-1}}{(r!)}\left(\sum_{j=1}^{k}x_j\frac{\nac{j}{n}\oq{{ij}}}{\nac{i}{n}\op{i}}\right)^r\nonumber\\
&=\frac{(1-r^{-1})^{r-1}}{r}\left(\sum_{j=1}^{k}x_j\frac{\nac{j}{n}\oq{{ij}}}{\nac{i}{n}\op{i}}\right)^r\to
r^{-1}(1-r^{-1})^{r-1}\left(\sum_{j=1}^{k}x_j\chi_{ij}\right)^r.\label{eq:limitdefbas}
\end{align}
The claim then follows by taking the limit as $n\to\infty$ in \eqref{eq:Brel}.
\\
\noindent$\square$

\noindent{\it Proof\,\,of\,\,Lemma\,\,\ref{le:June15}.}
We have
\begin{align}
&\sup_{\bold{x}\in \mathcal{W}}\Big|\frac{{\color{black}\oldB_n^{(i)}}(\lfloor\bold{x}\nan\rfloor,\obq{{i}})}{\nac{i}{n}}-\oldb_i(\bold x)\Big|\nonumber\\
&\leq\Big|\frac{\oan{i}}{\nac{i}{n}}-\alpha_i\Big|+\sup_{\bold{x}\in \mathcal{W}}\Big|\frac{\lfloor x_i \nac{i}{n}\rfloor}{\nac{i}{n}}-x_i\Big|\nonumber\\
&\qquad\qquad
+\sup_{\bold{x}\in \mathcal{W}}\Big|\frac{\on{i}-\oan{i}}{\nac{i}{n}}\oldpi(\lfloor\bold{x}\nan\rfloor,\obq{{i}})
-r^{-1}(1-r^{-1})^{r-1}\left(\sum_{j=1}^{k}x_j\chi_{ij}\right)^r\Big|\nonumber\\
&\leq\Big|\frac{\oan{i}}{\nac{i}{n}}-\alpha_i\Big|+\frac{1}{\nac{i}{n}}\nonumber\\
&\qquad\qquad
+\sup_{\bold{x}\in \mathcal{W}}\Big|\frac{\on{i}-\oan{i}}{\nac{i}{n}}\oldpi(\lfloor\bold{x}\nan\rfloor,\obq{{i}})
-r^{-1}(1-r^{-1})^{r-1}\left(\sum_{j=1}^{k}x_j^{(j)}\chi_{ij}\right)^r\Big|,\label{eq:supremumin}
\end{align}
and so we only need to prove that the supremum in \eqref{eq:supremumin} tends to zero as $n\to\infty$. For $\bold x\in \mathcal{W}$, we
consider quantity $O_k$ defined by \eqref{eq:O} with $\ell=k$
(note all $x_j$'s in $O_k$ are positive since $\bold x\in \mathcal{W}$). By Lemma \ref{le:aspi}
we have
\begin{align}
&\sup_{\bold x\in \mathcal{W}}\Big|\frac{\on{i}-\oan{i}}{\nac{i}{n}}\oldpi(\lfloor\bold{x}\nan\rfloor,\obq{{i}})
-r^{-1}(1-r^{-1})^{r-1}\left(\sum_{j=1}^{k}x_j\chi_{ij}\right)^r\Big|\nonumber\\
&\qquad\qquad
=\sup_{\bold x\in \mathcal{W}}\Big|\frac{\on{i}-\oan{i}}{\nac{i}{n}}(1+O_{k})
\left(\sum_{j=1}^{k}\lfloor x_j \nac{j}{n}\rfloor \oq{{ij}}\right)^{r}/r!
-r^{-1}(1-r^{-1})^{r-1}\left(\sum_{j=1}^{k}x_j\chi_{ij}\right)^r\Big|\nonumber\\
&\qquad\qquad
\leq\sup_{\bold x\in \mathcal{W}}\Big|\frac{\on{i}-\oan{i}}{\nac{i}{n}}O_{k}
\left(\sum_{j=1}^{k}x_j \nac{j}{n}\oq{{ij}}\right)^{r}/r!\Big|\nonumber\\
&\qquad\qquad\qquad
+\sup_{\bold x\in \mathcal{W}}\Big|\frac{\on{i}-\oan{i}}{\nac{i}{n}}
\left(\sum_{j=1}^{k}\lfloor x_j \nac{j}{n}\rfloor \oq{{ij}}\right)^{r}/r!
-r^{-1}(1-r^{-1})^{r-1}\left(\sum_{j=1}^{k}x_j\chi_{ij}\right)^r\Big|.\label{eq:suprema}
\end{align}
We start considering the first supremum in \eqref{eq:suprema}. For some positive constant $c_1>0$, we have
\begin{align}
&\sup_{\bold x\in \mathcal{W}}\Big|\frac{\on{i}-\oan{i}}{\nac{i}{n}}O_{k}
\left(\sum_{j=1}^{k}x_j \nac{j}{n}\oq{{ij}}\right)^{r}/r!\Big|\nonumber\\
&\qquad\qquad
\leq c_1\left(\sum_{j=1}^{k}(\nac{j}{n}\oq{{ij}}+(\nac{j}{n})^{-1})\right)
\frac{\on{i}-\oan{i}}{\nac{i}{n}}\left(\sum_{j=1}^{k}\nac{j}{n}\oq{{ij}}\right)^{r}\to 0,\nonumber
\end{align}
where the limit follows by \eqref{eq:ptcto0bm} and \eqref{eq:limitdefbas}.
As far as the second supremum in \eqref{eq:suprema} is concerned, by the definition of $\nac{i}{n}$,
we have, for some positive constants $c_2,c_3>0$,
\begin{align}
&\sup_{\bold x\in \mathcal{W}}\Big|\frac{\on{i}-\oan{i}}{\nac{i}{n}}
(r!)^{-1}\left(\sum_{j=1}^{k}\lfloor x_j \nac{j}{n}\rfloor \oq{{ij}}\right)^{r}
-r^{-1}(1-r^{-1})^{r-1}\left(\sum_{j=1}^{k}x_j\chi_{ij}\right)^r\Big|\nonumber\\
&\qquad\qquad
\leq\frac{\on{i}-\oan{i}}{\nac{i}{n}}
(r!)^{-1}\sup_{\bold x\in \mathcal{W}}\Big|\left(\sum_{j=1}^{k}\lfloor x_j \nac{j}{n}\rfloor \oq{{ij}}\right)^{r}-
\left(\sum_{j=1}^{k}x_j \nac{j}{n}\oq{{ij}}\right)^{r}\Big|\nonumber\\
&\qquad\qquad
+\sup_{\bold x\in \mathcal{W}}\Big|
(r!)^{-1}\frac{\on{i}-\oan{i}}{\nac{i}{n}}\left(\sum_{j=1}^{k}x_j \nac{j}{n}\oq{{ij}}\right)^{r}
-r^{-1}(1-r^{-1})^{r-1}\left(\sum_{j=1}^{k}x_j\chi_{ij}\right)^r\Big|\nonumber\\
&\qquad\qquad
\leq c_2\frac{\on{i}-\oan{i}}{\on{i}}
\sup_{\bold x\in \mathcal{W}}\Big|\left(\sum_{j=1}^{k}\frac{\lfloor x_j \nac{j}{n}\rfloor \oq{{ij}}}{\nac{i}{n}\op{i}}\right)^{r}-
\left(\sum_{j=1}^{k}\frac{x_j \nac{j}{n}\oq{{ij}}}{\nac{i}{n}\op{i}}\right)^{r}\Big|\nonumber\\
&\qquad\qquad
+r^{-1}(1-r^{-1})^{r-1}\sup_{\bold x\in \mathcal{W}}\Big|
\frac{\on{i}-\oan{i}}{\on{i}}\left(\sum_{j=1}^{k}x_j\frac{\nac{j}{n}\oq{{ij}}}{\nac{i}{n}\op{i}}\right)^{r}
-\left(\sum_{j=1}^{k}x_j\chi_{ij}\right)^r\Big|\nonumber\\
&\qquad\qquad
\leq c_2\frac{\on{i}-\oan{i}}{\on{i}}
\sup_{\bold x\in \mathcal{W}}\Big|\left(\sum_{j=1}^{k}\frac{\lfloor x_j \nac{j}{n}\rfloor \oq{{ij}}}{\nac{i}{n}\op{i}}+\sum_{j=1}^{k}\frac{\oq{{ij}}}{\nac{i}{n}\op{i}}\right)^{r}-
\left(\sum_{j=1}^{k}\frac{
\lfloor x_j\nac{j}{n}\rfloor \oq{{ij}}}{\nac{i}{n}\op{i}}\right)^{r}
\Big|\nonumber\\
&\qquad\qquad\qquad\qquad
+c_3\sup_{\bold x\in \mathcal{W}}\Big|
\frac{\on{i}-\oan{i}}{\on{i}}\left(\sum_{j=1}^{k}x_j\frac{\nac{j}{n}\oq{{ij}}}{\nac{i}{n}\op{i}}\right)^{r}
-\left(\sum_{j=1}^{k}x_j\chi_{ij}\right)^r\Big|\nonumber\\
&\qquad\qquad
\leq c_2\frac{\on{i}-\oan{i}}{\on{i}}
\sup_{\bold x\in \mathcal{W}}\Big|\left(\sum_{j=1}^{k}\frac{\lfloor x_j \nac{j}{n}\rfloor \oq{{ij}}}{\nac{i}{n}\op{i}}+\sum_{j=1}^{k}\frac{\oq{{ij}}}{\nac{i}{n}\op{i}}\right)^{r}-
\left(\sum_{j=1}^{k}\frac{
\lfloor x_j\nac{j}{n}\rfloor \oq{{ij}}}{\nac{i}{n}\op{i}}\right)^{r}
\Big|\nonumber\\
&\qquad\qquad\qquad\qquad
+c_3\frac{\on{i}-\oan{i}}{\on{i}}\sup_{\bold x\in \mathcal{W}}\Big|
\left(\sum_{j=1}^{k}x_j\frac{\nac{j}{n}\oq{{ij}}}{\nac{i}{n}\op{i}}\right)^{r}
-\left(\sum_{j=1}^{k}x_j\chi_{ij}\right)^r\Big|\nonumber\\
&\qquad\qquad\qquad\qquad
+c_3\sup_{\bold x\in \mathcal{W}}\Big|\frac{\on{i}-\nac{i}{n}}{\on{i}}\left(\sum_{j=1}^{k}x_j\chi_{ij}\right)^r-\left(\sum_{j=1}^{k}x_j\chi_{ij}\right)^r\Big|.\label{eq:supremabis}
\end{align}
Note that the latter supremum in the right-hand side of \eqref{eq:supremabis} goes to zero as $n\to\infty$. As far as the other two suprema in the right-hand side of \eqref{eq:supremabis},
note that, for any $\delta>0$ there exists $n_\delta$ such that for any $n\geq n_\delta$
\[
\chi_{ij}-\delta/k<\frac{\nac{j}{n}\oq{{ij}}}{\nac{i}{n}\op{i}}<\chi_{ij}+\delta/k
\]
and
\[
\sum_{j=1}^{k}\frac{\oq{{ij}}}{\nac{i}{n}\op{i}}<\delta.
\]
Therefore, for all $n$ large enough,
\begin{align}
&\sup_{\bold x\in \mathcal{W}}\Big|
\left(\sum_{j=1}^{k}x_j\frac{\nac{j}{n}\oq{{ij}}}{\nac{i}{n}\op{i}}\right)^{r}
-\left(\sum_{j=1}^{k}x_j\chi_{ij}\right)^r\Big|\nonumber\\
&\leq\sup_{\bold x\in \mathcal{W}}\bold{1}\left\{\sum_{j=1}^{k}x_j\frac{\nac{j}{n}\oq{{ij}}}{\nac{i}{n}\op{i}}\geq\sum_{j=1}^{k}x_j\chi_{ij}\right\}
\left[\left(\sum_{j=1}^{k}x_j\frac{\nac{j}{n}\oq{{ij}}}{\nac{i}{n}\op{i}}\right)^{r}
-\left(\sum_{j=1}^{k}x_j\chi_{ij}\right)^r\right]\nonumber\\
&+\sup_{\bold x\in \mathcal{W}}\bold{1}\left\{\sum_{j=1}^{k}x_j\frac{\nac{j}{n}\oq{{ij}}}{\nac{i}{n}\op{i}}<\sum_{j=1}^{k}x_j\chi_{ij}\right\}
\left[\left(\sum_{j=1}^{k}x_j\chi_{ij}\right)^r-\left(\sum_{j=1}^{k}x_j\frac{\nac{j}{n}\oq{{ij}}}{\nac{i}{n}\op{i}}\right)^{r}
\right]\nonumber\\
&\leq\sup_{\bold x\in \mathcal{W}}\Big|
\left(\sum_{j=1}^{k}x_j\chi_{ij}+\delta\right)^{r}
-\left(\sum_{j=1}^{k}x_j\chi_{ij}\right)^r\Big|+
\sup_{\bold x\in \mathcal{W}}\Big|\left(\sum_{j=1}^{k}x_j\chi_{ij}\right)^r-
\left(\sum_{j=1}^{k}x_j\chi_{ij}-\delta\right)^{r}
\Big|\nonumber
\end{align}
If $\sum_{j=1}^{k}\chi_{ij}=0$, then by the arbitrariness of $\delta$ we clearly have
\begin{equation}\label{eq:sup1zero}
\sup_{\bold x\in \mathcal{W}}\Big|
\left(\sum_{j=1}^{k}x_j\frac{\nac{j}{n}\oq{{ij}}}{\nac{i}{n}\op{i}}\right)^{r}
-\left(\sum_{j=1}^{k}x_j\chi_{ij}\right)^r\Big|\to 0.
\end{equation}
If $\sum_{j=1}^{k}\chi_{ij}>0$, then for all $n$ large enough
\begin{align}
\sup_{\bold x\in \mathcal{W}}\Big|
\left(\sum_{j=1}^{k}x_j\chi_{ij}\pm\delta\right)^{r}
-\left(\sum_{j=1}^{k}x_j\chi_{ij}\right)^r\Big|&\leq(\sup \mathcal{W})^r\left(\sum_{j=1}^{k}\chi_{ij}\right)^r\Big|
\left(1\pm\frac{\delta}{\varepsilon\sum_{j=1}^{k}\chi_{ij}}\right)^{r}-1\Big|,\nonumber
\end{align}
where $\sup \mathcal{W}:=\max_{1\leq i\leq k}\max_{\bold{x}\in \mathcal{W}}x_i$.
Therefore again \eqref{eq:sup1zero} follows by the arbitrariness of $\delta$. For all $n$ large enough, we also have
\begin{align}
&\sup_{\bold x\in \mathcal{W}}
\Big|\left(\sum_{j=1}^{k}\frac{\lfloor x_j \nac{j}{n}\rfloor \oq{{ij}}}{\nac{i}{n}\op{i}}+\sum_{j=1}^{k}\frac{\oq{{ij}}}{\nac{i}{n}\op{i}}\right)^{r}-
\left(\sum_{j=1}^{k}\frac{
\lfloor x_j\nac{j}{n}\rfloor \oq{{ij}}}{\nac{i}{n}\op{i}}\right)^{r}
\Big|\nonumber\\
&\qquad\qquad\qquad
\leq\sup_{\bold x\in \mathcal{W}}
\Big|\left(\sum_{j=1}^{k}\frac{\lfloor x_j \nac{j}{n}\rfloor \oq{{ij}}}{\nac{i}{n}\op{i}}+\delta\right)^{r}-
\left(\sum_{j=1}^{k}\frac{
\lfloor x_j\nac{j}{n}\rfloor \oq{{ij}}}{\nac{i}{n}\op{i}}\right)^{r}
\Big|\nonumber\\
&\qquad\qquad\qquad
\leq\sup_{\bold x\in \mathcal{W}}
\Big|\left(\sum_{j=1}^{k}\frac{\lfloor x_j \nac{j}{n}\rfloor \oq{{ij}}}{\nac{i}{n}\op{i}}+\delta\right)^{r}-
\left(\sum_{j=1}^{k}x_j\chi_{ij}\right)^{r}
\Big|\nonumber\\
&\qquad\qquad\qquad
+\sup_{\bold x\in \mathcal{W}}
\Big|
\left(\sum_{j=1}^{k}\frac{
\lfloor x_j\nac{j}{n}\rfloor \oq{{ij}}}{\nac{i}{n}\op{i}}\right)^{r}-\left(\sum_{j=1}^{k}x_j\chi_{ij}\right)^{r}
\Big|,\nonumber
\end{align}
and one can check that these two latter suprema go to zero as $n\to\infty$ arguing as in the proof of relation \eqref{eq:sup1zero}.
\\
\noindent$\square$

\noindent{\it Proof\,\,of\,\,Lemma\,\,\ref{le:initialpoint}.}\\
\noindent{\it Proof\,\,of\,\,$(i)$.} For $i=1,\ldots,k$ and $h=0,\ldots,{\color{black}\overline{d}}$, we have
\begin{align}
\oldb_i(\bold{x}_0^{(h)})&=\alpha_i-(\bold{x}_0^{(h)})_i+r^{-1}(1-r^{-1})^{r-1}\left(\sum_{j=1}^{k}(\bold{x}_0^{(h)})_j\chi_{ij}\right)^r\nonumber\\
&=\alpha_i-\sum_{s=0}^{h}\beta_s\bold{1}(i\in\mathcal{K}_s)+r^{-1}(1-r^{-1})^{r-1}\left(\sum_{j=1}^{k}\sum_{s=0}^{h}\beta_s\bold{1}(j\in\mathcal{K}_s)\chi_{ij}\right)^r\nonumber\\
&\geq\alpha_i-\sum_{s=0}^{h}\beta_s\bold{1}(i\in\mathcal{K}_s)+r^{-1}(1-r^{-1})^{r-1}\chi^r\left(\sum_{j=1}^{k}\sum_{s=0}^{h}\beta_s\bold{1}(j\in\mathcal{K}_s)\bold{1}\{\chi_{ij}>0\}\right)^r.
\label{eq:bi16Lug}
\end{align}
Therefore
\begin{align}
\oldb_1(\bold{x}_0^{(h)})
&\geq\alpha_1-\beta_0+r^{-1}(1-r^{-1})^{r-1}\chi^r\left(\sum_{j=1}^{k}\sum_{s=0}^{h}\beta_s\bold{1}(j\in\mathcal{K}_s)\bold{1}(\chi_{1j}>0)\right)^r\nonumber\\
&\geq\alpha_1-\beta_0=\alpha_1/2.\nonumber
\end{align}
\noindent{\it Proof\,\,of\,\,$(ii)$.} Since $i\in\bigcup_{s=1}^{h}\mathcal{K}_s$, then
there exists $j_0\in\{1,\ldots,k\}$ such that ${\color{black}d_{j_0}}={\color{black}d_i}-1$ and $\chi_{ij_0}>0$
(where  quantities ${\color{black}d_i}$ are defined in Subsection \ref{sect:addnot}). So by
\eqref{eq:bi16Lug} we have
\begin{align}
\oldb_i(\bold{x}_0^{(h)})&\geq\alpha_i-\beta_{{\color{black}d_i}}+r^{-1}(1-r^{-1})^{r-1}(\chi\beta_{{\color{black}d_i}-1})^r\nonumber\\
&\geq\frac{r^{-1}(1-r^{-1})^{r-1}(\chi\beta_{{\color{black}d_i}-1})^r}{2}>0,\nonumber
\end{align}
where the second inequality follows by \eqref{eq:beta}.\\
\noindent{\it Proof\,\,of\,\,$(iii)$.} The proof is similar to the proof of part $(ii)$ and therefore omitted.\\
\noindent{\it Proof\,\,of\,\,$(iv)$.} Since $i\notin\bigcup_{s=1}^{h+1}\mathcal{K}_s$, we have ${\color{black}d_i}\geq h+2$.
Therefore, $\bold{1}(i\in\mathcal{K}_s)=0$ for any $0\leq s\leq h$ and so
$\sum_{s=0}^{h}\beta_s\bold{1}(i\in\mathcal{K}_s)=0$. Moreover,
\[
\bold{1}(\chi_{ij}>0)\sum_{s=0}^{h}\beta_s\bold{1}(j\in\mathcal{K}_s)=0,\quad\text{$\forall$ $j=1,\ldots,k$}
\]
(otherwise one can find a $j_0\in\mathcal{K}_s$, $0\leq s\leq h$, with $\chi_{ij_0}>0$ and so ${\color{black}d_i}\leq h+1$, which is a contradiction).
Therefore, by \eqref{eq:bi16Lug} we have $\oldb_i(\bold{x}_0^{(h)})\geq\alpha_i\geq 0$.\\
\noindent{\it Proof\,\,of\,\,$(v)$.} It is an obvious consequence of $(i)$ and $(ii)$.
\\
\noindent$\square$

\subsection{Proof of Propositions \ref{prop:alphagreaterthan1} and \ref{prop:alpharegion}}

\noindent{\it Proof\,\,of\,\,Proposition\,\,\ref{prop:alphagreaterthan1}.} If $\alpha_1>1$, then $\oldb_1$ has no zeros in $[0,\frac{r}{r-1}]^k$, and therefore
$({\mathcal Sup})$ holds. If $\alpha_1=1$, then $\oldb_1$ has a unique zero at $\bold{z}:=(r/(r-1),0,\ldots,0)$. Since $\bm{\chi}$ is irreducible,
there exists $j\in\{2,\ldots,k\}$ such that $\oldb_j(\bold{z})>0$. So $\bold{z}\notin\oD_{\bm\oldb}$, and therefore $({\mathcal Sup})$ holds.
\\
\noindent$\square$

The proof of Proposition \ref{prop:alpharegion} exploits the following lemma which will be proved later on.

\begin{Lemma}  \label{zerounico-crit}
Let $\bm{\alpha}\ge\bold 0$. If there exists $\bold{z}\in\overset\circ{\mathcal D}$ such that ${\bm\oldb}(\bold{z},\bm{\alpha})=0$, $\mathrm{det}J_{\bm\oldb}(\bold z)=0$ and $\lambda_{PF}(\bold z)=0$,
then $\bold{(Crit)}$ holds, i.e., $\mathcal Z=\{\bold z\}$.
\end{Lemma}

\noindent{\it Proof\,\,of\,\,Proposition\,\,\ref{prop:alpharegion}.} We divide the proof in two steps. In the first step we prove that, for $\bm{\alpha}_*\geq\bold 0$, the following statements are equivalent:\\
\noindent$(i)$ $\bm{\alpha}_*\in\mathcal{R}_{Crit}$;\\
\noindent$(ii)$ There exists $\bold{z}_*\in\overset\circ{\mathcal D}$ such that
$\oldb_1(\bold{z}_*,\bm{\alpha}_*)=0$, $\bold{z}_*\in\mathcal{D}_{{\bm\oldb},\bm{\alpha_*}}$ and
\begin{equation}\label{eq:autosoluzione}
\bold{v}(\bold z_*)J_{{\bm\oldb},\bm{\alpha}_*}(\bold z_*)=\bold 0,\quad\text{for some $\bold{v}(\bold z_*):=(v_1(\bold z_*),\ldots,v_k(\bold z_*))>\bold 0$.}
\end{equation}
\\
In the second step we use this equivalence to determine  region $\mathcal{R}_{Crit}$ (and consequently, $\mathcal{R}_{Sub}$ and $\mathcal{R}_{Sup}$).\\
\noindent{\it Step\,\,1.} We start proving $(i)\Rightarrow(ii)$. By $(\bold{Crit})$ we have $\mathcal{Z}=\{\bold z_*\}$, for some
$\bold{z}_*\in\overset\circ{\mathcal D}$, $\mathrm{det}J_{{\bm\oldb},\bm{\alpha}_*}(\bold z_*)=0$ and $\lambda_{PF}(\bold z_*)=0$.
Since $\mathrm{det}J_{{\bm\oldb},\bm{\alpha}_*}(\bold z_*)=0$, we have
$\bold{v}J_{{\bm\oldb},\bm{\alpha}_*}(\bold z_*)=\bold 0$ for some $\bold{v}:=(v_1,\ldots,v_k)\neq\bold 0$. Since $\lambda_{PF}(\bold z_*)=0$ we have $\bold{v}=\bm{\phi}_{PF}(\bold z_*)$,
and so $\bold v$ is (component-wise) positive. We now prove $(ii)\Rightarrow(i)$. By $(ii)$ it follows that
$\mathrm{det}J_{{\bm\oldb},\bm{\alpha}_*}(\bold z_*)=0$ and  that $\bold{v}(\bold z_*)=\phi_{PF}(\bold z_*)$
(since $\bold{v}(\bold z_*)$ is (component-wise) positive). Therefore $\lambda_{PF}(\bold z_*)=0$ and
by Lemma \ref{zerounico-crit} we have that $(\bold{Crit})$ holds, and so $\bm{\alpha}_*\in\mathcal{R}_{Crit}$.\\
\noindent{\it Step\,\,2.} By \eqref{eq:autosoluzione} we have, for some $\bm{\theta}:=(1,\theta_2,\ldots,\theta_k)$, $\theta_i>0$, $i=2,\ldots,k$,
\[
\sum_{i=1}^{k}\theta_i\frac{\partial \oldb_i(\bold x,\alpha_i^*)}{\partial x_h}\Big|_{\bold{x}=\bold{z}_*}=0,\quad\text{$h=1,\ldots,k$}
\]
i.e.,
\begin{align}
&
\sum_{i=1}^{k}\theta_i\frac{\partial \oldb_i(\bold x,\alpha_i^*)}{\partial x_1}\Big|_{\bold{x}=\bold{z}_*}\nonumber\\
&\qquad\qquad
=-1+(1-r^{-1})^{r-1}\left(\sum_{j=1}^{k}\chi_{1j}z_j^*\right)^{r-1}
+(1-r^{-1})^{r-1}\sum_{i=2}^{k}\theta_i\left(\sum_{j=1}^{k}\chi_{ij}z_j^*\right)^{r-1}\chi_{i1}=0,\label{eq:L1}
\end{align}
\begin{align}
\sum_{i=1}^{k}\theta_i\frac{\partial \oldb_i(\bold x,\alpha_i^*)}{\partial x_h}\Big|_{\bold{x}=\bold{z}_*}
&
=(1-r^{-1})^{r-1}\left(\sum_{j=1}^{k}\chi_{1j}z_j^*\right)^{r-1}\chi_{1h}\nonumber\\
&\qquad
+\sum_{i=2}^{k}\theta_i\left(-\bold{1}_{\{i=h\}}+(1-r^{-1})^{r-1}\left(\sum_{j=1}^{k}\chi_{ij}z_j^*\right)^{r-1}\chi_{ih}\right)=0,\quad\text{$h=2,\ldots,k$.}
\label{eq:L2}
\end{align}
Setting
\begin{equation}\label{eq:xprime}
x_i':=\left[(1-r^{-1})\sum_{j=1}^{k}\chi_{ij}z_j^*\right]^{r-1},\quad i=1,\ldots,k
\end{equation}
\eqref{eq:L1} and \eqref{eq:L2} reduce to
\[
x_1'+\sum_{i=2}^{k}\theta_i\chi_{i1}x_{i}'=1,\nonumber
\]
\[
x_1'\chi_{1h}+\sum_{i=2}^{k}\theta_i\chi_{ih}x_{i}'=\theta_{h},\quad\text{$h=2,\ldots,k$.}
\]
In matrix form these relations can be rewritten as
\begin{equation*}
\bold{x}'\mathrm{diag}(\bm{\theta})\bm{\chi}=\bm{\theta}
\end{equation*}
where $\bold x':=(x_1',\ldots,x_k')$,
and $\mathrm{diag}(\bm{\theta})$ is the $k\times k$ diagonal matrix with diagonal elements $(1,\theta_2,\ldots,\theta_k)$.
Since $\bm{\chi}$ is invertible,
we have
\begin{equation*}
\bold x'=\bm{\theta}\bm{\chi}^{-1}\mathrm{diag}(\bm{\theta}^{-1}).
\end{equation*}
By this latter relation and \eqref{eq:xprime} we get
\begin{equation}\label{eq:regcrit1}
\sum_{j=1}^{k}\chi_{ij}z_j^*=(\bold{x}_{\bm{\theta}})_i,
\quad\text{$i=1,\ldots,k$}
\end{equation}
and so
\begin{equation}\label{eq:regcrit2}
z_i^*=(\bold{x}_{\bm{\theta}}(\bm{\chi}^{-1})^t)_i,\quad\text{$i=1,\ldots,k$.}
\end{equation}
The claim follows inserting expressions \eqref{eq:regcrit1} and \eqref{eq:regcrit2} into $\oldb_i(\bold{z}_*,\alpha_i^*)$, $i=1,\ldots,k$, and then solving
equations $\oldb_i(\bold{z}_*,\alpha_i^*)=0$ with respect to $\alpha_i^*$, $i=1,\ldots,k$.
\\
\noindent$\square$

\noindent{\it Proof\,\,of\,\,Lemma\,\, \ref{zerounico-crit}.} By the strict convexity of $\oldb_i$, $i=1,\ldots,k$, and the fact that $J_{\bm\oldb}(\bold z)\bold\ov_{PF}(\bold z)=\bold 0$,
for any $\gamma\in\mathbb R$ such that $\bold{z}+\gamma\bold{\ov}_{PF}(\bold z)\in\oD$, we have
\begin{equation}\label{posbi}
\oldb_i(\bold{z}+\gamma\bold{\ov}_{PF}(\bold z))\ge 0,\quad i=1,\ldots,k.
\end{equation}
Reasoning by contradiction, assume that there exists $\bold{z}_1\in\mathcal{Z}$, $\bold z_1\neq\bold z$. Then, by the strict convexity of $\oldb_i$,
$i=1,\ldots,k$, we have ${\bm\oldb}((1-\theta)\bold{z}+\theta\bold{z}_1)<\bold 0$, for any $\theta\in (0,1)$. For $\theta_0\in (0,1)$ sufficiently small,
define $\bold w= (1-\theta_0)\bold{z}+\theta_0\bold{z}_1$ and consider  sets $\hat{\mathcal{F}}_{\bold w}$ and $\hat{\mathcal{F}}^{(i)}_{\bold w}$, $i=1,\ldots,k$, defined by
\eqref{eq:setF}.
Then, for any $i=1,\ldots,k$,
\begin{equation}\label{negbi}
\oldb_i(\bold y)<0,\quad\text{$\forall$ $\bold{y}\in\hat{\mathcal{F}}^{(i)}_{\bold w}$.}
\end{equation}
Since $\bold \ov_{PF}(\bold z)>0$, we have $\bold{z}+\gamma_0\bold{\ov}_{PF}(\bold z)\in
\hat{\mathcal{F}}_{\bold w}$, for some $\gamma_0\in\mathbb R$, and therefore $\bold{z}+\gamma_0\bold{\ov}_{PF}(\bold z)\in
\hat{\mathcal{F}}_{\bold w}^{(i_0)}$, for some $i_0\in\{1,\ldots,k\}$. Combining this with
\eqref{posbi} and \eqref{negbi} we get a contradiction.
\\
\noindent$\square$


\begin{thebibliography}{10}

\bibitem{abbe15}
E. Abbe and C. Sandon.
\newblock Community detection in general stochastic block models: fundamental limits and efficient algorithms for recovery.
\newblock {\em Proceedings of the IEEE 56th Annual Symposium on Foundations of Computer Science}, 670--688, 2015.

\bibitem{abbe16}
E. Abbe, A.S. Bandeira and G. Hall.
\newblock Exact recovery in the stochastic block model.
\newblock {\em IEEE Transactions on Information Theory}, 62: 471--487, 2016.

\bibitem{fountolakis2}
M. Abdullah and N. Fountoulakis.
\newblock A phase transition in the evolution of bootstrap percolation processes on preferential attachment graphs.
\newblock {\em Random Structures and Algorithms}, 52: 379-418, 2018.

\bibitem{airoldi08}
E.M. Airoldi, D.M. Blei, S.E. Fienberg and E.P. Xing.
\newblock Mixed membership stochastic block models.
\newblock {\em  Journal of Machine Learning Research}, 9: 1981--2014, 2008.

\bibitem{amini1}
H. Amini.
\newblock Bootstrap percolation and diffusion in random graphs with given vertex degrees.
\newblock {\em Electronic Journal of Combinatorics}, 17: 1--20, 2010.

\bibitem{amini2}
H. Amini and N. Fountoulakis.
\newblock Bootstrap percolation in power-law random graphs.
\newblock {\em Journal of Statistical Physics}, 155: 72--92, 2014.

\bibitem{angel}
O. Angel and B. Kolesnik.
\newblock Large deviations for subcritical bootstrap percolation on the random graph.
\newblock {\em ArXiv}: 1705.06815v2, 2018.

\bibitem{Arnold}
V.I. Arnold.
\newblock {\em Ordinary Differential Equations}.
\newblock Springer-Verlag, Berlin, 1992.

\bibitem{bollobas}
J. Balogh and B. Bollob\'{a}s.
\newblock Bootstrap percolation on the hypercube.
\newblock {\em Probability Theory and Related Fields}, 134: 624--648, 2006.

\bibitem{BPP}
J. Balogh, Y. Peres and  G. Pete.
\newblock Bootstrap percolation on infinite trees and non-amenable groups.
\newblock {\em Combinatorics, Probability and Computing}, 15: 715--730, 2006.

\bibitem{pittel}
J. Balogh and  B.G. Pittel.
\newblock Bootstrap percolation on the random regular graph.
\newblock {\em Random Structures and Algorithms}, 30: 257--286, 2007.

\bibitem{galton}
B. Bollob\'{a}s, K. Gunderson, C. Holmgren, S. Janson and M. Przykucki.
\newblock Bootstrap percolation on Galton-Watson trees.
\newblock {\em Electronic Journal of Probability}, 19: 1--27, 2014.

\bibitem{belajigsaw}
B. Bollob\'{a}s, O. Riordan, E. Slivken and P. Smith.
\newblock The threshold for jigsaw percolation on random graphs.
\newblock {\em Electronic Journal of Combinatorics}, 24, Paper \#P2.36, 2017.

\bibitem{rgg}
M. Bradonji\'{c} and I. Saniee.
\newblock Bootstrap percolation on random geometric graphs.
\newblock {\em Probability in the Engineering and Informational Sciences,} 28: 169--181, 2014.

\bibitem{jigsaw}
C.D. Brummitt, S. Chatterjee, P.S. Dey and D. Sivakoff.
\newblock Jigsaw percolation: what social networks can collaboratively solve a puzzle?.
\newblock {\em The Annals of Applied Probability}, 25: 2013--2038, 2015.

\bibitem{chalupa}
J. Chalupa, P.L. Leath and G.R. Reich.
\newblock Bootstrap percolation on a Bethe lattice.
\newblock {\em Journal of Physics C}, 12: 31--35, 1979.

\bibitem{amin10}
A. Coja-Oghlan and A. Lanka.
\newblock Finding planted partitions in random graphs with general degree distributions.
\newblock {\em SIAM Journal on Discrete Mathematics}, 23: 1682--1714, 2010.


\bibitem{feige}
U. Feige, M. Krivelevich and D. Reichman.
\newblock Contagious sets in random graphs.
\newblock {\em The Annals of Applied Probability}, 27: 2675--2697, 2016.

\bibitem{fortunato}
S. Fortunato.
\newblock Community detection in graphs.
\newblock {\em Physics Reports}, 486(3--5), 75--174, 2010.

\bibitem{fountolakis} N. Fountoulakis, M. Kang, C. Koch, T. Makai.
\newblock A phase transition regarding the evolution of bootstrap processes in inhomogeneous random graphs.
\newblock {\em The Annals of Applied Probability}, 28: 990--1051, 2018.

\bibitem{girvanpnas02}
M. Girvan and M.E.J. Newman.
\newblock Community structure in social and biological networks.
\newblock {\em Proceedings of the National Academy of Sciences}, 99: 7821--7826, 2002.

\bibitem{gopalan13}
P.K. Gopalan and D.M. Blei.
\newblock Efficient discovery of overlapping communities in massive networks.
\newblock {\em Proceedings of the National Academy of Sciences}, 110: 14534--14539, 2013.

\bibitem{cecilia17}
C. Holmgren, T. Ju\v{s}kevi\v{c}ius, N. Kettle.
\newblock Majority bootstrap percolation on $G_{n,p}$.
\newblock {\em Electronic Journal of Combinatorics}, 24, Paper \#P1.1, 2017.

\bibitem{JLTV}
S. Janson, T. Luczak, T. Turova and T. Vallier.
\newblock Bootstrap percolation on the random graph $G_{n,p}$.
\newblock {\em The Annals of Applied Probability}, 22: 1989--2047, 2012.

\bibitem{kozma16}
S. Janson, R. Kozma, M. R\'{o}bert, M. Ruszink\'{o} and Y. Sokolov.
\newblock Bootstrap percolation on a random graph coupled with a lattice.
\newblock {\em Electronic Journal of Combinatorics}, in press.

\bibitem{newman11}
B Karrer, M.E.J. Newman.
\newblock Stochastic block models and community structure in networks.
\newblock {\em Physical Review E}, 83: 016107, 2011.

\bibitem{kempe}
D. Kempe, J. Kleinberg and E. Tardos.
\newblock Maximizing the spread of influence through a social network.
\newblock {\em Proceedings of the 9th ACM SIGKDD international conference on knowledge discovery and data mining}, 10 pp., 2003.

\bibitem{LT}
P. Lancaster and M. Tismenetski.
\newblock {\em The Theory of Matrices}.
\newblock Academic Press, New York, 1985.

\bibitem{mass14}
L. Massouli\'{e}.
\newblock Community detection thresholds and the weak Ramanujan property.
\newblock {\em Proceedings of the 46th annual ACM symposium on theory of computing}, 694--703, 2014.

\bibitem{P}
M.D. Penrose.
\newblock Random geometric graphs.
\newblock {\em Oxford University Press}, New York, USA, 2003.

\bibitem{Scalia}
G.P. Scalia-Tomba.
\newblock Asymptotic final-size distribution for some chain-binomial processes.
\newblock {\em Advances in Applied Probability}, 17: 477--495, 1985.

\bibitem{munik}
M. Shrestha and C. Moore.
\newblock Message-passing approach for threshold models of behavior in networks.
\newblock{\em Physical Review E}, 89: 022805, 2014.

\bibitem{TGL} G.L. Torrisi, M. Garetto and E. Leonardi.
\newblock A large deviation approach to super-critical bootstrap percolation on the random graph $G_{n,p}$.
\newblock{\em Stochastic Processes and their Applications}, to appear, 2018.

\bibitem{turova15}
T.S. Turova and T. Vallier.
\newblock Bootstrap percolation on a graph with random and local connections.
\newblock {\em Journal of Statistical Physics}, 160: 1249--1276, 2015.

\bibitem{watts}
D. Watts.
\newblock A simple model of global cascades in random networks.
\newblock {\em Proceedings of the National Academy of Sciences}, 99: 5766--5771, 2002.
\end{thebibliography}
\end{document}